\RequirePackage{import}
\documentclass[english,ruled, 10pt]{article}
\newif\ifsubmission
\submissionfalse
\usepackage[T1]{fontenc}
\usepackage[utf8]{inputenc}
\usepackage{amsfonts, amsmath, amssymb, amsthm}
\usepackage{geometry}
\usepackage{algorithm, algpseudocode}
\usepackage{nicematrix}
\usepackage{bm}
\usepackage{xcolor}
\usepackage{optidef}
\usepackage{pstricks, pst-tree}
\usepackage{graphicx}
\usepackage[square,numbers]{natbib}
\usepackage{mathabx}
\usepackage{mathtools}
\usepackage{listings}
\usepackage{float}
\usepackage{caption}
\usepackage{subcaption}
\usepackage{pifont}
\usepackage{babel}
\usepackage{hyperref}
\usepackage{array}
\usepackage{multirow}
\hypersetup{
	colorlinks=true,
	linkcolor=red,
	citecolor = blue
}
\usepackage{changes}
\usepackage{enumitem}

\theoremstyle{plain}
\newtheorem*{assumption*}{\protect\assumptionname}

\theoremstyle{remark}
\newtheorem*{remark*}{\protect\remarkname}

\theoremstyle{plain}

\theoremstyle{plain}
\newtheorem{theorem}{Theorem}[section]

\theoremstyle{definition}

\theoremstyle{plain}
\newtheorem{assumption}{\protect\assumptionname}

\theoremstyle{plain}
\newtheorem{proposition}{\protect\propositionname}

\newtheorem{lemma}[theorem]{Lemma}

\theoremstyle{plain}
\newtheorem{corollary}{\protect\corollaryname}

\providecommand{\assumptionname}{Assumption}
\providecommand{\corollaryname}{Corollary}
\providecommand{\definitionname}{Definition}
\providecommand{\propositionname}{Proposition}
\providecommand{\remarkname}{Remark}

\providecommand{\Halmos}{\ensuremath{\blacksquare}}

\newenvironment{sharedproof}[1][Proof]
{%
    \begin{proof}[#1]%
}
{%
    \unskip\nobreak\hfill\Halmos
    \end{proof}%
}

\numberwithin{equation}{section}

\newif\ifsiam

\siamfalse 

\title{First-order Methods for  Variational Inequality Problems with Semi-Infinite Constraints
}
\author{
Dureksha Thilakasiri\thanks{dthilakasiri@smu.edu, Operations Research and Engineering Management, Southern Methodist University}
\hspace{6em}
Digvijay Boob\thanks{dboob@smu.edu, Operations Research and Engineering Management, Southern Methodist University}
}
\date{}

\begin{document}
\allowdisplaybreaks
\global\long\def\vertiii#1{\left\vert \kern-0.25ex  \left\vert \kern-0.25ex  \left\vert #1\right\vert \kern-0.25ex  \right\vert \kern-0.25ex  \right\vert }%
\global\long\def\matr#1{\bm{#1}}%
\global\long\def\til#1{\tilde{#1}}%
\global\long\def\wt#1{\widetilde{#1}}%
\global\long\def\wh#1{\widehat{#1}}%
\global\long\def\wb#1{\widebar{#1}}%
\global\long\def\mcal#1{\mathcal{#1}}%
\global\long\def\mbb#1{\mathbb{#1}}%
\global\long\def\mtt#1{\mathtt{#1}}%
\global\long\def\ttt#1{\texttt{#1}}%
\global\long\def\inner#1#2{\langle#1,#2\rangle}%
\global\long\def\binner#1#2{\big\langle#1,#2\big\rangle}%
\global\long\def\Binner#1#2{\Big\langle#1,#2\Big\rangle}%
\global\long\def\br#1{\left(#1\right)}%
\global\long\def\bignorm#1{\bigl\Vert#1\bigr\Vert}%
\global\long\def\Bignorm#1{\Bigl\Vert#1\Bigr\Vert}%
\global\long\def\setnorm#1{\Vert#1\Vert_{-}}%
\global\long\def\rmn#1#2{\mathbb{R}^{#1\times#2}}%
\global\long\def\deri#1#2{\frac{d#1}{d#2}}%
\global\long\def\pderi#1#2{\frac{\partial#1}{\partial#2}}%
\global\long\def\onebf{\mathbf{1}}%
\global\long\def\zero{\mathbf{0}}%

\global\long\def\norm#1{\lVert#1\rVert}%
\global\long\def\bnorm#1{\big\Vert#1\big\Vert}%
\global\long\def\Bnorm#1{\Big\Vert#1\Big\Vert}%

\global\long\def\brbra#1{\big(#1\big)}%
\global\long\def\Brbra#1{\Big(#1\Big)}%
\global\long\def\rbra#1{(#1)}%
\global\long\def\sbra#1{[#1]}%
\global\long\def\bsbra#1{\big[#1\big]}%
\global\long\def\Bsbra#1{\Big[#1\Big]}%
\global\long\def\cbra#1{\{#1\}}%
\global\long\def\bcbra#1{\big\{#1\big\}}%
\global\long\def\Bcbra#1{\Big\{#1\Big\}}%

\global\long\def\grad{\nabla}%
\global\long\def\Expe{\mathbb{E}}%
\global\long\def\rank{\text{rank}}%
\global\long\def\range{\text{range}}%
\global\long\def\diam{\text{diam}}%
\global\long\def\epi{\text{epi }}%
\global\long\def\inte{\operatornamewithlimits{int}}%
\global\long\def\cov{\text{Cov}}%
\global\long\def\argmin{\operatornamewithlimits{argmin}}%
\global\long\def\argmax{\operatornamewithlimits{argmax}}%
\global\long\def\tr{\operatornamewithlimits{tr}}%
\global\long\def\dis{\operatornamewithlimits{dist}}%
\global\long\def\sign{\operatornamewithlimits{sign}}%
\global\long\def\prob{\mathbb{P}}%
\global\long\def\st{\operatornamewithlimits{s.t.}}%
\global\long\def\dom{\text{dom}}%
\global\long\def\diag{\text{diag}}%
\global\long\def\and{\text{and}}%
\global\long\def\st{\text{s.t.}}%
\global\long\def\Var{\operatornamewithlimits{Var}}%
\global\long\def\raw{\rightarrow}%
\global\long\def\law{\leftarrow}%
\global\long\def\Raw{\Rightarrow}%
\global\long\def\Law{\Leftarrow}%
\global\long\def\vep{\varepsilon}%
\global\long\def\dom{\operatornamewithlimits{dom}}%

\global\long\def\Lbf{\mathbf{L}}%

\global\long\def\Ffrak{\mathfrak{F}}%
\global\long\def\Gfrak{\mathfrak{G}}%
\global\long\def\vfrak{\mathfrak{v}}%
\global\long\def\LGfrak{{\lambda\mathfrak{G}}}%
\global\long\def\gfrak{\mathfrak{g}}%
\global\long\def\sfrak{\mathfrak{s}}%
\global\long\def\xibar{\bar{\xi}}%
\global\long\def\Cbb{\mathbb{C}}%
\global\long\def\Ebb{\mathbb{E}}%
\global\long\def\Fbb{\mathbb{F}}%
\global\long\def\Nbb{\mathbb{N}}%
\global\long\def\Rbb{\mathbb{R}}%
\global\long\def\extR{\widebar{\mathbb{R}}}%
\global\long\def\Pbb{\mathbb{P}}%
\global\long\def\Acal{\mathcal{A}}%
\global\long\def\Bcal{\mathcal{B}}%
\global\long\def\Ccal{\mathcal{C}}%
\global\long\def\Dcal{\mathcal{D}}%
\global\long\def\Fcal{\mathcal{F}}%
\global\long\def\Gcal{\mathcal{G}}%
\global\long\def\Hcal{\mathcal{H}}%
\global\long\def\Ical{\mathcal{I}}%
\global\long\def\Kcal{\mathcal{K}}%
\global\long\def\Lcal{\mathcal{L}}%
\global\long\def\Mcal{\mathcal{M}}%
\global\long\def\Ncal{\mathcal{N}}%
\global\long\def\Ocal{\mathcal{O}}%
\global\long\def\Pcal{\mathcal{P}}%
\global\long\def\Qcal{\mathcal{Q}}%
\global\long\def\Ucal{\mathcal{U}}%
\global\long\def\Scal{\mathcal{S}}%
\global\long\def\Tcal{\mathcal{T}}%
\global\long\def\Xcal{\mathcal{X}}%
\global\long\def\Ycal{\mathcal{Y}}%
\global\long\def\Vcal{\mathcal{V}}%
\global\long\def\Ubf{\mathbf{U}}%
\global\long\def\Pbf{\mathbf{P}}%
\global\long\def\Ibf{\mathbf{I}}%
\global\long\def\Ebf{\mathbf{E}}%
\global\long\def\Abs{\boldsymbol{A}}%
\global\long\def\Qbs{\boldsymbol{Q}}%
\global\long\def\Lbs{\boldsymbol{L}}%
\global\long\def\Pbs{\boldsymbol{P}}%
\global\long\def\i{i}%
\global\long\def\Ibb{\mathbb{I}}

\DeclarePairedDelimiterX{\inprod}[2]{\langle}{\rangle}{#1, #2}
\DeclarePairedDelimiter\abs{\lvert}{\rvert}
\DeclarePairedDelimiter{\bracket}{ [ }{ ] }
\DeclarePairedDelimiter{\paran}{(}{)}
\DeclarePairedDelimiter{\braces}{\lbrace}{\rbrace}
\DeclarePairedDelimiterX{\gnorm}[3]{\|}{\|_{#2}^{#3}}{#1}
\DeclarePairedDelimiter{\floor}{\lfloor}{\rfloor}
\DeclarePairedDelimiter{\ceil}{\lceil}{\rceil}

\global\long\def\tsum{{\textstyle {\sum}}}

\newcommand{\opex}{\texttt{OpEx}}
\newcommand{\opconex}{\texttt{OpConEx}}
\newcommand{\pdconex}{\texttt{PD-ConEx}}
\newcommand{\sopconex}{\texttt{P-OpEx}}
\newcommand{\stochsopconex}{\texttt{S-P-OpEx}}
\newcommand{\fstochsopconex}{\texttt{FS-P-OpEx}}
\newcommand{\adopconex}{\texttt{AdOpEx}}
\newcommand{\gconex}{\texttt{GradConEx}}
\newcommand{\stopconex}{\texttt{StOpConEx}}
\newcommand{\fstopconex}{\texttt{F-StOpConEx}}

\maketitle

\begin{abstract}
This paper studies first-order methods for monotone variational inequalities with semi-infinite constraints and presents this model as a framework covering constrained saddle-point problems, robust optimization, distributionally robust function-constrained optimization, and generalized Nash equilibrium formulations. Since feasibility must hold over infinitely many indexed constraints, direct projection onto the full feasible region is generally difficult. Therefore, we propose a single-loop primal-dual constraint extrapolation method (PD-ConEx) that updates primal variables, Lagrange multipliers, and worst-case constraint indices using first-order information, linearization, extrapolation, and projections only onto simple sets. We establish convergence guarantees for smooth and nonsmooth regimes, showing how the achievable rates depend on strong monotonicity and strong concavity of the index block. For coupled convex-concave saddle-point problems, strong convexity in the minimization block can accelerate its smooth terms even without strong concavity in the maximization block. We also develop an adaptive smooth variant with $O(1/K)$ VI-gap and feasibility guarantees without an a priori bound on generated multipliers. Moment-based and Wasserstein distributionally robust models fit within the framework, and numerical experiments on robust quadratic instances and distributionally robust learning illustrate the method's optimality and feasibility behavior.
\end{abstract}
\section{Introduction}

Variational inequality (VI) problems provide a general mathematical framework for modeling constrained optimization, equilibrium, game-theoretic, and saddle-point problems. For a given closed convex set $\widetilde X \subseteq \mathbb{R}^n$ and a {\em monotone} operator $F:\mathbb{R}^n \to \mathbb{R}^n$, the standard VI problem seeks
\begin{equation}\label{eq:vi}
x^\ast \in \wt{X}\ \st\ \inner{F(x^\ast)}{x - x^\ast} \ge 0,\quad \forall\, x \in \wt{X}.
\end{equation}
The VI problem \eqref{eq:vi} has various applications in finance, economics,
transportation, communication, and many others.  We refer the reader to the monologue \cite{FacchineiPang2003} for classic applications and mathematical foundations of the VI area. The VI area has witnessed significant algorithmic research due to its wide applicability in equilibrium problems, adversarial machine learning, and reinforcement learning.

In this paper, we study a monotone VI problem in which the feasible set is defined by semi-infinite function constraints. Specifically, the feasible region is given by
\begin{equation}\label{eq:Xtilde}
\wt{X}
=
X \cap
\cbra{
x:\ g_i(x,y_i)\le 0,\ \forall\,y_i\in Y_i,\ i=1,\dots,m
}.
\end{equation}
Here, $X\subset\Rbb^n$ is a simple convex compact set that admits efficient projection, each set $Y_i$ indexes a family of constraints, and the function $g_i$ is convex with respect to $x$ and concave with respect to $y_i$ for all $i  =1, \dots, m$. The term semi-infinite refers to the fact that there are finitely many constraint families, but each family contains infinitely many inequalities indexed by $y_i\in Y_i$. We say that problem \eqref{eq:vi} is smooth if the operator $F$ and the gradients of the constraints $g_i$, for all $i = 1, \dots, m$, are Lipschitz continuous. If any of these assumptions are not satisfied, then we call it a {\em composite} nonsmooth problem. We use composite to highlight that in a general nonsmooth VI problem, we can have both $F$ and $g$ as sum of smooth and nonsmooth components. In such a setup, we would like to accelerate convergence for the smooth components even when the nonsmooth part is present. In the composite nonsmooth case, we assume that the nonsmooth component of $F$ and the subgradient of the nonsmooth component of $g_i$, for all $i =1, \dots, m$, remain bounded over the sets $X$ and $Y_i$, respectively.

The main computational difficulty in \eqref{eq:vi}--\eqref{eq:Xtilde} is that projection onto $\wt{X}$ is generally hard. Although projection onto the simple set $X$ is assumed to be easy, the additional semi-infinite constraints may destroy this projection-friendly structure. As a result, standard first-order VI methods that rely on exact projection onto the feasible region cannot be directly applied. This difficulty is closely related to function-constrained VI problems \cite{BoobDengKhalafi2024}, but the semi-infinite structure introduces an additional challenge because feasibility must hold uniformly over the index sets $Y_i$.

Semi-infinite constraints arise naturally in several important models. In robust optimization, constraints must hold for all realizations of uncertain parameters in an uncertainty set \cite{BenTalNemirovski1998, BenTalElGhaouiNemirovski2009}. In distributionally robust optimization, objectives or constraints are often defined through worst-case expectations over ambiguity sets of probability distributions \cite{Scarf1958, RahimianMehrotra2019, EsfahaniKuhn2018}. Similar structures also appear in fairness literature where fairness is imposed by bounding a function whose tractable formulation involves an index variable similar to $y_i$ in \eqref{eq:Xtilde}. In order to limit the violation of fairness, it is common to impose hard bounds on fairness measuring functions, which naturally leads to a semi-infinite form. These applications motivate the development of scalable first-order methods that can handle semi-infinite feasibility without requiring projection onto the feasible set.

\subsection{Literature Review} 
The classical VI literature mainly studied structural properties and algorithms such as Newton-type, quasi-Newton, and interior-point methods \cite{FacchineiPang2003}. These methods can have strong local convergence properties under suitable regularity assumptions, but they are often less suitable for large-scale problems where scalable first-order methods are needed.

This limitation has motivated a large body of work on first-order methods for monotone VIs. Korpelevich's extragradient method \cite{Korpelevich1976} and Nemirovski's Mirror-Prox method \cite{Nemirovski2004} are foundational methods for smooth monotone VIs and obtain the optimal $O(1/K)$ rate. Nesterov's dual extrapolation method \cite{Nesterov2007} provides a related optimal first-order framework, and faster convergence is possible under strong monotonicity \cite{NesterovScrimali2011}. Other developments include optimistic gradient methods \cite{Popov1980}, projection-reflected and golden-ratio methods \cite{Malitsky2015,Malitsky2020}, and operator extrapolation methods \cite{KotsalisLanLi2022}. These works form the foundation of modern first-order VI algorithms. However, most of them assume that the feasible set is either unconstrained or simple enough that exact projection can be computed efficiently.

The requirement of efficient projection onto the feasible set becomes restrictive when the set is described by function constraints. In such problems, projecting onto the full feasible region may require solving a difficult auxiliary optimization problem. Recent works have addressed this issue for function-constrained VIs using primal-dual, interior-point, and constraint-extrapolation ideas. Boob, Deng, and Khalafi \cite{BoobDengKhalafi2024} developed first-order methods for function-constrained VIs, including adaptive operator extrapolation methods for smooth deterministic problems and related methods for nonsmooth and stochastic settings. Zhang et al.~\cite{zhang2024primal} proposed a constrained gradient method for variational inequality problems with function constraints, while Yang, Jordan, and Chavdarova \cite{YangJordanChavdarova2022} studied constrained VIs through an interior-point approach. Constraint-extrapolation ideas also appear in function-constrained optimization, where Boob, Deng, and Lan \cite{BoobDengLan2019} developed first-order methods based on linearization and extrapolation of function constraints.

Our study is closely related to this line of work, but the feasible set considered here has an additional semi-infinite structure. In our setting, each constraint must hold for all values of an index variable $y_i\in Y_i$. Therefore, feasibility is determined by infinitely many inequalities rather than by finitely many function constraints. This makes direct projection even more challenging than in finite function-constrained VI problems. 

From an optimization point of view, \cite{YaoLinYang2024} studied first-order methods for convex semi-infinite optimization. Though their work is closely related to ours, it is developed for convex optimization problems. In contrast, our paper studies monotone VI problems, which are substantially more general than optimization, effectively modeling generalized robust optimization, distributionally robust optimization with distributionally robust function constraints, saddle point problems with coupled constraints, and equilibrium problems including generalized Nash equilibrium. This also allows us to explore further structural information of saddle point problems to accelerate convergence rate guarantees in the presence of semi-infinite constraints.

The use of semi-infinite constraints arises naturally in robust and distributionally robust optimization. In robust optimization, feasibility must hold for all realizations of uncertain parameters in an uncertainty set \cite{BenTalNemirovski1998, BenTalElGhaouiNemirovski2009}. First-order robust optimization methods have been developed to avoid expensive robust counterparts or repeated adversarial subproblems \cite{HoNguyenKilincKarzan2018, PostekShtern2024}. Despite the recent focus on this topic, the developed algorithms require knowledge of the optimal Lagrange multiplier and do not analyze differentiated convergence for smooth and nonsmooth cases.  In distributionally robust optimization, the objective function considers the worst-case of expectation function over the {\em ambiguity set} of probability distributions. The foundation of this distributionally robust model was laid in the classical work of Scarf \cite{Scarf1958}. This research area became very active in recent years, with a focus on developing approximation schemes to reduce the computationally challenging distributionally robust objective functions with tractable (mostly convex) function, or considering different ambiguity sets based on application or computation requirements. Interested reader can see the recent survey \cite{RahimianMehrotra2019}. Despite many significant developments, the majority of the literature cannot handle distributional robustness in data-driven function constraints.
.Randomized and stochastic first-order methods for VIs have also been studied to reduce the cost of evaluating the operator, especially in large-scale settings. For example, stochastic operator extrapolation and related randomized first-order schemes use noisy or sampled operator evaluations instead of full deterministic evaluations \cite{KotsalisLanLi2022, BoobDengKhalafi2024}. These methods are important for scalability, but they typically assume that the feasible set is simple enough for efficient projection or that the constraints are finite-dimensional and can be handled directly. Hence, they do not directly address the projection difficulty created by semi-infinite function constraints.

Moreover, such randomized VI algorithms are not designed to handle generalized Nash equilibrium problems (GNEPs) with coupled constraints. In GNEPs, each player's feasible decision may depend on the decisions of other players through common constraints, and in our setting these constraints can also be semi-infinite. Therefore, existing randomized VI methods do not provide a direct framework for GNEPs with semi-infinite coupling constraints, which further motivates the primal--dual VI--SIP approach developed in this paper.

Clearly, the applications discussed are major fundamental research areas in their own right. However, there are no systematic algorithmic studies discussing faster, possibly differentiated convergence guarantees for the nonsmooth VI problem, ensuring scalability through simple single-loop design with only first-order information, and not imposing restrictive assumptions such as projection onto set $\wt{X}$ or knowledge of the optimal Lagrange multiplier. On the contrary, despite many important developments, existing first-order VI methods mainly rely on projection onto a simple feasible set, function-constrained VI methods primarily address finitely many function constraints, and convex semi-infinite optimization methods do not directly apply to general monotone VIs. This leaves open the development of first-order algorithms for monotone VI problems with semi-infinite function constraints, which have simple single loop design with feasibility ensured without a direct projection onto the set $\wt{X}$ and differentiated complexity results for smooth/nonsmooth components along with impacts of strong monotonicity or convexity. Below, we briefly discuss our contributions in this direction which, to the best of our knowledge, appear to be the first in the literature.

\subsection{Contributions}
\label{sec:contributions}

In this paper, we develop first-order methods for monotone variational inequality problems with semi-infinite functional constraints. Our main contributions are summarized as follows.

First, we provide a primal--dual reformulation of the VI--SIP problem. Since feasibility is determined by infinitely many inequalities indexed by \(y_i\in Y_i\), direct projection onto the feasible region is generally difficult. Under an appropriate constraint qualification, a solution \(x^\ast\) admits Lagrange multipliers \(\lambda^\ast\in\Rbb_+^m\) and worst-case index variables
\[
y_i^\ast\in\argmax_{y_i\in Y_i}g_i(x^\ast,y_i),\qquad i=1,\ldots,m,
\]
satisfying KKT-type optimality conditions. This characterization yields a structured primal--dual formulation in the lifted variable \(z=(x,\lambda,y)\). We further introduce a primal--dual gap function and decompose it into components associated with the primal, multiplier, and index-variable updates, which forms the main analytical tool for the convergence analysis.

Second, we propose a primal--dual constraint extrapolation method, referred to as PD-ConEx. The method updates the worst-case constraint indices, Lagrange multipliers, and primal variable within a single loop. The \(y\)-update performs blockwise proximal ascent, the \(\lambda\)-update uses a projected proximal step driven by extrapolated constraint values, and the \(x\)-update applies a projected first-order step involving the VI operator and linearized constraint terms. Each iteration requires only first-order evaluations and projections onto the simple sets \(X\), \(\Rbb_+^m\), and \(Y_i\), thereby avoiding direct projection onto the semi-infinite feasible region \(\wt{X}\).

Third, we establish differentiated componentwise convergence guarantees under several curvature and smoothness regimes. Let \(\mu_F\geq0\) denote the strong monotonicity modulus of the VI operator, let \(\mu_y\geq0\) denote the strong concavity modulus of \(g_i(x,\cdot)\), and let \(L_g^{yy}\) measure the smoothness of the constraints with respect to the index variable. When \(\mu_F>0\) and \(\mu_y>0\), the smooth operator, constraint, and initialization components decay at the accelerated rate \(O(1/K^2)\). When \(\mu_F>0\) and \(\mu_y=0\), the strongly monotone \(x\)-related smooth terms retain the \(O(1/K^2)\) rate, while the smooth \(y\)-related component, including the term involving \(L_g^{yy}\), decays as \(O(1/K)\). Thus, removing strong concavity in \(y\) affects only the index-variable component and does not eliminate acceleration of the remaining smooth terms. When \(\mu_F=0\) and \(\mu_y>0\), strong concavity improves the control of the \(y\)-block, but the absence of strong monotonicity leaves the dominant smooth rate at \(O(1/K)\). When \(\mu_F=\mu_y=0\), all smooth components attain the standard first-order rate \(O(1/K)\).

In the nonsmooth setting, the smooth components retain the rates described above, while nonsmoothness introduces additional componentwise errors. For \(\mu_F>0\) and \(\mu_y>0\), the smooth contribution remains \(O(1/K^2)\), whereas the nonsmooth \(x\)- and \(y\)-components both decay as \(O(1/K)\). For \(\mu_F>0\) and \(\mu_y=0\), the smooth \(x\)-related terms remain \(O(1/K^2)\), the smooth \(y\)-related term is \(O(1/K)\), while the nonsmooth \(x\)- and \(y\)-components decay as \(O(1/K)\) and \(O(1/\sqrt{K})\), respectively. For \(\mu_F=0\) and \(\mu_y>0\), the smooth contribution is \(O(1/K)\), while the nonsmooth \(x\)- and \(y\)-components contribute \(O(1/\sqrt{K})\) and \(O(\log(K)/K)\), respectively. Finally, when \(\mu_F=\mu_y=0\), the smooth terms remain \(O(1/K)\), while both nonsmooth components contribute \(O(1/\sqrt{K})\). Therefore, our results identify which components retain acceleration and which components determine the overall rate when individual curvature assumptions are removed. The constants may depend on \(\|\lambda^\ast\|_1\), but this dependence does not change the order in \(K\).

Fourth, we extend the VI--SIP framework to convex--concave saddle-point problems with coupled semi-infinite constraints:
\[
\min_{u\in U}\max_{v\in V}f(u,v)
\quad\text{subject to}\quad
g_i(u,v,y_i)\leq0,\qquad \forall y_i\in Y_i,\quad i=1,\ldots,m.
\]
The resulting method employs separate descent and ascent updates for the \(u\)- and \(v\)-blocks, together with multiplier and worst-case index updates. Let \(\mu_F^u\) and \(\mu_F^v\) denote the strong convexity and strong concavity moduli of the \(u\)- and \(v\)-blocks, respectively, and let \(L_F^u\) and \(L_F^v\) denote their block smoothness constants. Strong convexity in \(u\) allows the \(u\)-related smooth terms to attain an accelerated \(O(1/K^2)\) rate, whereas the absence of strong concavity in \(v\) generally leaves the \(v\)-related terms at \(O(1/K)\). In particular, when
\[
\mu_F^u>0,\qquad \mu_F^v=0,\qquad \mu_y=0,
\]
the \(u\)-related smooth terms decay as \(O(1/K^2)\), while the \(v\)- and \(y\)-related smooth terms decay as \(O(1/K)\). Consequently, the overall smooth rate is \(O(1/K)\), and the dominant nonsmooth rate is \(O(1/\sqrt{K})\). If \(f\) is linear in \(v\), then the Lipschitz constant associated with the variation of \(F_v\) with respect to \(v\) is zero. Consequently, the corresponding \(v\)-smoothness contribution vanishes, although coupling terms involving \(u\) and the semi-infinite constraints may remain.

Finally, we show that the proposed framework captures robust and distributionally robust optimization models. In particular, moment-based and Wasserstein ambiguity sets can be represented within the VI--SIP formulation. These applications naturally contain minimization, maximization, and worst-case index blocks, illustrating the relevance of the differentiated curvature regimes and the saddle-point extension. Numerical experiments on a robust quadratic semi-infinite problem and a moment-based distributionally robust learning problem using the Mushroom dataset demonstrate the behavior of PD-ConEx in reducing primal--dual gaps and feasibility violations.
\subsection{Notation and Preliminaries}
We use Euclidean norm throughout the paper, denoted as \(\norm{\cdot}\). and inner product as \(\inner{\cdot}{\cdot}\). 
We denote the set $\{1, \dots, m\}$ as $[m]$. The primal domain is a nonempty, convex, compact set \(X\subset\Rbb^{n}\).
There are \(m\) semi-infinite constraints. For each \(i \in [m]\), the index set
\(Y_i\subset\Rbb^{d_i}\) is nonempty, convex, and compact. Define the set
\(Y:=Y_1\times\cdots\times Y_m\) and the joint semi-infinite index variable $y:=
[ (y^{1})^\top,\dots,(y^{m})^\top ]^\top \in Y \subset \Rbb^{N}$. 
\(\Pi_X\) denotes Euclidean projection onto \(X\). We use similar notation for projections onto each \(Y_i\), $Y$ and the nonnegative orthant \(\Rbb_+^{m}\). 
The normal cone to \(X\) at \(x\in X\) is
\[
N_X(x):=\Bcbra{ v\in\Rbb^{n}:\ \inner{v}{x'-x}\le 0,\ \forall x'\in X }.
\]
The diameter of \(X\) is \(D_X:=\max_{x_1,x_2\in X}\norm{x_1-x_2}\). We use similar notation for set $Y$.\\
\textbf{Constraint mapping.}
The vector-valued constraint map is $g(x,y):=(g_1(x,y_1),\dots,g_m(x,y_m))\in\Rbb^{m}.$
Then, \(g_i(x,y_i)\le 0\) for all \(y_i\in Y_i\) and all \(i\) is equivalent to $g(x, y) \le 0$ for all $y\in Y$.\\
\textbf{Norms and weighted-norms.}
For multipliers \(\lambda=(\lambda_1,\dots,\lambda_m)^\top\in\Rbb_+^{m}\), we use the \(\ell_1\)-norm \(\norm{\lambda}_1=\tsum_{i=1}^m|\lambda_i|\).
Given \(y,\tilde y\in Y\), we define the \(\lambda\)-weighted-norm
\[
\norm{y-\tilde y}_{\lambda}^{2}:=\tsum_{i=1}^m \lambda_i \norm{y_i-\tilde y_i}^{2},
\]
which is used in proximal steps. Throughout the paper, all norms are assumed to be Euclidean.
We define block-diagonal matrix $\Lambda$ 
\begin{align}
\Lambda
:=
\text{Block-Diag}(
\lambda_1 I_{n_1},
\lambda_2 I_{n_2}, 
\dots, \lambda_m I_{n_m})
\in \Rbb^{N\times N}.\label{eq:Lambda_def}
\end{align}
\begin{assumption}[Convexity and Compactness] \label{ass:convxity_n_compactness}
We make the following assumptions throughout our paper.
\begin{enumerate}[label=(\alph*)]
    \item \(X\subset\Rbb^{n}\) and \(Y_i\subset\Rbb^{d_i}\) for \(i=1,\dots,m\) are convex and compact. \label{ass:compactness}
    \item There exists \(D_y\in\Rbb\) such that \(\norm{y_i-y_i'}\le D_y\) for any \(y_i,y_i'\in Y_i\) and all \(i=1,\dots,m\).
    \item The operator \(F:X\to\Rbb^{n}\) is \(\mu_F\)- monotone on \(X\) with \(\mu_F\ge 0\); i.e.,
    \[
    \inner{F(x_1)-F(x_2)}{x_1-x_2}\ \ge\ \mu_F\norm{x_1-x_2}^{2},\qquad \forall\,x_1,x_2\in X.
    \]
    If $\mu_F = 0$, then $F$ is monotone operator. Otherwise, it is a strongly monotone operator.
    \item For each \(i=1,\dots,m\) and every \(y_i\in Y_i\), the mapping \(x\mapsto g_i(x,y_i)\) is convex on \(X\).
    \item For each \(i=1,\dots,m\) and every \(x\in X\), the mapping \(y_i\mapsto g_i(x,y_i)\) is \(\mu_y\)-strongly concave on \(Y_i\), with \(\mu_y\ge 0\). If $\mu_y = 0$, then $g_i$ is concave on $Y_i$.
\end{enumerate}    
\end{assumption} 

\textbf{Subdifferentials.}
To incorporate nonsmoothness, we use subdifferentials with respect to both \(x\) and \(y_i\).
In view of Assumption \ref{ass:convxity_n_compactness}, we note that $g_i(\cdot, y_i)$ is a convex function for all $y_i \in Y_i$. Hence, we denote its subdifferential with respect to $x$ at point $(\wt{x}, \wt{y}_i) \in X \times Y_i$ as $\partial_xg_i(x,y_i)$.  We denote the subgradient of $g_i(x, y_i)$, with respect to $x$ 
as \( g'_{i,x}(\wt{x},\wt{y}_i)\). It is clear that $g'_{i,x}(\wt{x}, \wt{y}_i) \in \partial_{x}g_i(\wt{x}, \wt{y}_i)$. Note that $-g_i(x, \cdot)$ is a convex function of $y_i$. Hence, we define subdifferential of $g_i$ with respect to $y_i$ as $\partial_y g_i(\wt{x}, \wt{y}_i) := - \partial_y (-g_i(\wt{x}, \wt{y_i}))$. It is clear from the definition that for any $v \in \partial_y g_i(\wt{x}, \wt{y}_i)$, we have \[g_i(\wt{x}, y_i) \le g_i(\wt{x}, \wt{y}_i) + \inprod{v}{y_i - \wt{y}_i}.\]
This definition of subdifferential is a natural extension to the concave function $g_i(x, \cdot)$ (see \cite{rockafellar} for more details). We denote the subgradient of $g_i$ with respect to $y_i$ as $g'_{i,y}(\wt{x}, \wt{y}_i)$. Clearly, $g'_{i,y}(\wt{x}, \wt{y}_i) \in \partial_y g_i(x, y_i)$.

In this section, all assumptions and quantities are treated as deterministic. Specifically, the monotonicity of \(F\), the convexity/concavity properties of \(g\),and all Lipschitz constants introduced earlier are assumed to be fixed values.
\begin{assumption}[Existence of optimal solutions and multipliers] \label{ass:opt_sol}
For any solution \(x^{*}\) of the variational inequality \(\mathrm{VI}(F,\wt{X})\), there exist Lagrange multipliers \(\lambda^{*}\in\Rbb_{+}^{m}\) and the points
\(y_{i}^{*}\in\argmax_{y\in Y_{i}}  g_{i}\br{x^{*},y}\) for \(i=1,\dots,m\),
such that the following KKT-type conditions hold:
\begin{itemize}
\item \emph{Primal feasibility (semi-infinite):} \(g_{i}\br{x^{*},y}\le 0\) for all \(y\in Y_{i}\) and all \(i\), equivalently define the worst-case 
\emph{\(\wb{g}_{i}(x) := \sup_{y \in Y_{i}} g_{i}(x,y)\)},then primal feasibility means \emph{\(\wb{g}_{i}(x^{*}) \leq 0\) for every \(i\)}
\item \emph{Dual feasibility:} \(\lambda^{*}\ge 0\).
\item \emph{Complementarity:} \(\lambda_{i}^{*}  g_{i}^{*}\paran{x^{*}}=0\) for each \(i=1,\dots,m\).
\item \emph{Stationarity:}
\(0\in F\br{x^{*}}+\tsum_{i=1}^{m}\lambda_{i}^{*}\partial_{x} g_{i}\br{x^{*},y_{i}^{*}}+N_{X}\br{x^{*}}\).
\end{itemize}
\end{assumption}

The KKT-type stationarity condition for the VI--SIP states that, at the solution \(x^{*}\), the operator \(F\) is exactly balanced by a conic combination of subgradients of the semi-infinite constraints together with a normal vector that encodes the geometry of the feasible set \(X\). If \(x^{*}\in \inte X\), then \(N_X\br{x^{*}}=\cbra{0}\) and the condition simplifies to \(F\br{x^{*}}+\tsum_{i=1}^{m}\lambda_{i}^{*} s_{i}=0\).
In smooth cases we may take \(s_i = g'_{i,x}(x^{*}, y_i^{*})\); in nonsmooth cases any \(s_{i}\in\partial_{x} g_{i}\br{x^{*},y_{i}^{*}}\) is valid.
By complementarity, only indices with \(\lambda_{i}^{*}>0\) active constraints.

\begin{assumption}[Smoothness and Lipschitz Continuity] \label{ass:Smoothness_Lipcts}
There exist nonnegative constants \\
\(L_{F},L_{g}^{xx},L_{g}^{yx},L_{g}^{yy},H,H_{g}^{x},H_{g}^{y},M_{g}^{x},M_{g}^{y}\) such that:
\begin{itemize}
\item[(a)] 
\(\norm{F(x)-F(x')}\ \le\ L_{F}\norm{x-x'}\ +\ H,
 \forall x,x'\in X.\)
\item[(b)] 
\(g_i(x,y_i)-g_i(x',y_i)-\inner{g'_{i,x}(x',y_i)}{x-x'}\ 
\le\ \tfrac{L_{g}^{xx}}{2} \norm{x-x'}^{2} + H_{g}^{x} \norm{x-x'},\) for all \(x,x'\in X\), \(y_i\in Y_i\), and \(i=1,\dots,m\).
\item[(c)] 
\(\norm{ g'_{i,y}(x,y_i)-g'_{i,y}(x',y_i') }\ 
\le\ L_{g}^{yx} \norm{x-x'}\ +\ L_{g}^{yy} \norm{y_i-y_i'}\ +\ H_{g}^{y},\) for all \(x,x'\in X\), \(y_i,y_i'\in Y_i\), and \(i=1,\dots,m\).
\item[(d)] 
\(\abs{ g_i(x,y_i)-g_i(x',y_i) }\ 
\le\ M_{g} \norm{x-x'}\ 
,\) for all \(x,x'\in X\), 
and \(i=1,\dots,m\).
\end{itemize}
\end{assumption}
\noindent
This assumption covers both smooth and nonsmooth deterministic cases. In the \emph{smooth} regime we set \(H=H_{g}^{x}=H_{g}^{y}=0\); then items (a)–(c) assert Lipschitz properties for \(F\), \(g'_{i,x}(x,y_i)\), and \(g'_{i,y}(x,y_i)\) with \(L_{F}, L_{g}^{xx}, L_{g}^{yx}, L_{g}^{yy}\). In particular, when \(H=0\), item (a) reduces to \(F\) being Lipschitz continuous—often called a \emph{smooth operator} in the VI literature. In the \emph{nonsmooth} regime we take \(L_{F}=L_{g}^{xx}=L_{g}^{yx}=L_{g}^{yy}=0\); then (a)–(c) yield bounded-variation conditions and, together with (d), imply the uniform bounds
\begin{equation*}
\norm{g'_{i,x}(x,y_i)}\ \le\ M_{g},\qquad
\forall x\in X,\ 
\forall i \in [m],
\end{equation*}
where, in the nonsmooth case, \(\grad\) denotes a selected subgradient and \(\norm{\cdot}\) is the Euclidean norm. 
The convergence rates under the different settings are summarized in Table~\ref{tab:rate-summary}.
Based on Assumption~\ref{ass:Smoothness_Lipcts}, we say that the VI--SIP problem is smooth if $H=H_g^x=H_g^y=0.$
In this case, the operator $F$ is Lipschitz continuous and the relevant gradients of the constraint functions are Lipschitz continuous in the primal and index variables. If at least one of these smoothness conditions fails, we call the problem nonsmooth. In the nonsmooth setting, the functions may not be differentiable everywhere, and the analysis is performed using bounded subgradients.
\begin{table}[htbp]
\centering
\small
\renewcommand{\arraystretch}{1.15}
\setlength{\tabcolsep}{6pt}
\begin{tabular}{c c c c p{0.53\textwidth}}
\hline
\textbf{Method}
&
\textbf{Problem}
&
\textbf{Curvature in \(x\)}
&
\(\boldsymbol{\mu_y}\)
&
\textbf{Componentwise bound}
\\
\hline

\multirow{4}{*}{Alg.~\ref{alg:vi-sip-agsip}}
&
\multirow{4}{*}{VI--SIP}
&
\(\mu_F>0\)
&
\(>0\)
&
\(\displaystyle
O\!\left(
\frac{1+M_g^2/\mu_F}{K^2}
+
\frac{(H_g^x)^2}{\mu_F K}
+
\frac{(H_g^y)^2}{\mu_y K}
\right)\)
\\[2pt]

&
&
\(\mu_F>0\)
&
\(=0\)
&
\(\displaystyle
O\!\left(
\frac{L_g^{yy}}{K}
+
\frac{\mathcal L_{\mathrm{oth}}+M_g^2/\mu_F}{K^2}
+
\frac{(H_g^x)^2}{\mu_F K}
+
\frac{H_g^y}{\sqrt K}
\right)\)
\\[2pt]

&
&
\(\mu_F=0\)
&
\(>0\)
&
\(\displaystyle
O\!\left(
\frac{\mathcal L_{\mathrm{sm}}+M_g}{K}
+
\frac{H_g^x}{\sqrt K}
+
\frac{(H_g^y)^2\log K}{\mu_y K}
\right)\)
\\[2pt]

&
&
\(\mu_F=0\)
&
\(=0\)
&
\(\displaystyle
O\!\left(
\frac{\mathcal L_{\mathrm{sm}}+M_g}{K}
+
\frac{H_g^x+H_g^y}{\sqrt K}
\right)\)
\\
\hline

Alg.~\ref{alg:prob2}
&
Saddle-point VI--SIP
&
\(\mu_F^u>0,\ \mu_F^v=0\)
&
\(=0\)
&
\(\displaystyle
O\!\left(
\frac{L_F^v+M_g^2/L_F^v+L_g^{yy}}{K}
+
\frac{1+M_g^2/\mu_F^u}{K^2}
+
\frac{(H_g^u)^2}{\mu_F^u K}
+
\frac{H_g^v+H_g^y}{\sqrt K}
\right)\)
\\
\hline
\end{tabular}
\caption{Componentwise bounds for the proposed methods. The smooth cases
are obtained by setting the nonsmooth constants \(H_g^x,H_g^u,H_g^v,H_g^y\)
to zero. The quantities \(\mathcal{L}_{\mathrm{oth}}\) collects the smooth Lipschitz
constants appearing in the \(K^{-2}\) terms, while
\(\mathcal{L}_{\mathrm{sm}}\) collects the smooth Lipschitz and initial-distance constants appearing in the \(K^{-1}\) terms.}
\label{tab:rate-summary}
\end{table}

\section{The Primal-Dual Constraint Extrapolation Method}

Considering the Lagrangian $\Lcal(x,\lambda,y)$ and iterates $(x_k, y_k, \lambda_k)$, we recall that many saddle–point schemes update $(y_k,\lambda_k)$ simultaneously, while our analysis proceeds with  block updates and extrapolation made to order the VI with semi–infinite constraints\\
$x^*\in X:\ \inner{F(x^*)}{x-x^*}\ge 0,\ \forall x\in X,\ g_i(x,y_i)\le 0,\ \forall y_i\in Y_i$.\\
We linearize the constraints in $x$ using $\ell g(x; x', y) := g(x',y) + g'_{x}(x',y) (x-x')$, thereby admitting a proximal linear update for the dual variable and using the bounds by $L_g^{xx}$ and $H_g^x$.

The weighted norm $\norm{v}_\lambda^{2} = \tsum_{i=1}^m \lambda_i \norm{v_i}^{2}$ then appears naturally in the three–point lemma for the $y$–block, producing extra terms in the inequality. Accordingly, for each block $i$ we form the extrapolated gradient and take a proximal ascent step
\begin{align*}
u_i^k \gets (1+\theta_k) g'_{i,y}(x^k, y_i^k) - \theta_k g'_{i,y}(x^{k-1}, y_i^{k-1}),\\
y_i^{k+1} \gets \argmin_{y_i\in Y_i}\{-\inner{u_i^k}{y_i} + \tfrac{\sigma_k}{2}\norm{y_i - y_i^k}^{2}\},
\end{align*}
which implements momentum in the $y$–direction while controlling movement by the proximal term.\\
For the multiplier, we update by a projected proximal ascent by the linearized  $s_k$,
\begin{align*}
\lambda^{k+1} \gets \argmin_{\lambda\ge 0}\{- \inner{s_k}{\lambda}+\tfrac{\gamma_k}{2}\norm{\lambda-\lambda^k}^2\},
\end{align*}
Once we have updated $(y^{k+1}, \lambda^{k+1})$, we can generate $x^{k+1}$.
Under the standing assumptions, $F$ is $\mu_F$ strongly monotone and $L_F$ Lipschitz; each $g_i$ is smooth with Lipschitz constants $L_g^{xx}, L_g^{yx}, L_g^{yy}$ and remainder bounds $H_g^{x}, H_g^{y}$; and $g_i(x,\cdot)$ is $\mu_y$–strongly concave.
\begin{algorithm}[H]
\caption{Primal--Dual Constraint Extrapolated Method for VI--SIP (\pdconex)}
\label{alg:vi-sip-agsip}
\begin{algorithmic}[1]
\State \textbf{Inputs:} $x^0\in X$, $y^0\in Y$, $\lambda^0\in\Rbb_+^m$, $K\in\Nbb_+$; step sizes $\sigma_k, \gamma_k, \theta_k, \tau_k$.
\State \textbf{Init:} $x^{-2}=x^{-1}=x^0$, $y^{-1}=y^0$.
\For{$k=0,1,\dots,K-1$}
  \For{$i=1,\dots,m$}
    \State $u_i^k \gets (1+\theta_k) g'_{i,y}({x^k, y_i^k}-\theta_k g'_{i,y}({x^{k-1},y_i^{k-1}}$
    \State $y_i^{k+1} \gets \arg\min_{y_i\in Y_i}\{- \inner{u_i^k}{y_i}+\tfrac{\sigma_k}{2} \norm{y_i-y_i^k}^2\}$
  \EndFor
  \State $s_k \gets \ell_{g}({x^k; x^{k-1},y^{k+1}} \;+\; \theta_k[\ell_{g}({x^k; x^{k-1},y^{k}}-\ell_{g}({x^{k-1}; x^{k-2},y^{k}}]$
  \State $\lambda^{k+1} \gets \arg\min_{\lambda\ge 0}\{- \inner{s_k}{\lambda}+\tfrac{\gamma_k}{2} \norm{\lambda-\lambda^k}^2\}$
  \State $w_k \gets F({x^k} \;+\; \theta_k(F({x^k}-F({x^{k-1}}) \;+\; g'_{x}({x^k, y^{k+1}}^{\top}\lambda^{k+1}$
  \State $x^{k+1} \gets \arg\min_{x\in X}\{\inner{w_k}{x}+\tfrac{\tau_k}{2} \norm{x-x^k}^{2}\}$
\EndFor
\State \textbf{Output:}  $\displaystyle \wb{x}_K \gets \tfrac{\tsum_{k=0}^{K-1} t_k x^{k+1}}{\tsum_{k=0}^{K-1} t_k}$.
\end{algorithmic}
\end{algorithm}

\subsection{Convergence Analysis of \pdconex~method}
Let $z=({x,\lambda,y}$ and $z_k=({x_k,\lambda_k,y_k}$ for $k\ge 0$. We analyze convergence using the gap
\begin{equation} \label{eq:4(a)}
\begin{aligned}
Q({z_{k+1},z} := \inner{F({x}}{x_{k+1}-x} + \lambda^\top g({x_{k+1},y} - \lambda_{k+1}^\top g({x,y_{k+1}} .
\end{aligned}
\end{equation}
\vspace{-0.5em}
We decompose it as
\vspace{-0.5em}
\begin{align*}
Q({z_{k+1},z} = Q_0({z_{k+1},z} + Q_1({z_{k+1},z} + Q_2({z_{k+1},z} ,
\end{align*}
\vspace{-0.5em}
where
\begin{align*}
Q_0({z_{k+1},z} &:= \lambda^\top({g({x_{k+1},y}-g({x_{k+1},y_{k+1}}},\\
Q_1({z_{k+1},z} &:= ({\lambda-\lambda_{k+1}}^\top g({x_{k+1},y_{k+1}},\\
Q_2({z_{k+1},z} &:= \inner{F({x}}{x_{k+1}-x} + \lambda_{k+1}^\top({g({x_{k+1},y_{k+1}}-g({x,y_{k+1}}}.
\end{align*}
The convergence properties of the algorithm are derived by bounding $Q({z_{k+1},z}$ from above by summing up some upper bounds of $Q_0({z_{k+1},z}$, $Q_1({z_{k+1},z}$, and $Q_2({z_{k+1},z}$ and choosing parameters  $t_k,\theta_k,\tau_k,\gamma_k,\sigma_k \ge 0$ for $k=0,1,..$  appropriately. The following lemma is needed to obtain the upper bounds.
\begin{lemma}\label{lemma1}
Assume that $f : S \to \Rbb$ satisfies
\[
f(y) \ge f(x) + \inner{f'(x)}{y - x} + \tfrac{\mu}{2}\norm{y - x}^2,
\quad \forall  x,y \in S,
\]
for some $\mu \ge 0$, where $S$ is a closed convex set. If
\[
\bar{x} = \argmin_{x \in S} \{ f(x) + \tfrac{1}{2}\norm{x - x_{1}}^2 \},
\]
then for all $x \in S$,
\[
f(\bar{x}) + \tfrac{1}{2}\norm{\bar{x} - x_{1}}^2
+ \tfrac{(\mu+1)}{2}\norm{x - \bar{x}}^2
\le f(x) + \tfrac{1}{2}\norm{x - x_{1}}^2 .
\]
\end{lemma}
\begin{lemma}
Define $q_{\lambda}^k := \ell_g({x_k; x_{k-1}, y_k} - \ell_g({x_{k-1}; x_{k-2}, y_k}$. For any $z=(x,\lambda,y)$ with $\lambda\ge 0$, the upper bound for the combined 
gap term $Q_1(z_{k+1},z)+Q_2(z_{k+1},z)$.
\begin{align}
&Q_1({z_{k+1},z}+Q_2({z_{k+1},z}
\le \tfrac{\gamma_k}{2}({\norm{\lambda-\lambda_k}^{2}-\norm{\lambda-\lambda_{k+1}}^{2}})
- \inner{q_{\lambda}^{k+1}}{\lambda_{k+1}-\lambda}
+ \theta_k \inner{q_{\lambda}^{k}}{\lambda_k-\lambda} \nonumber\\
&\quad + \theta_k \inner{q_{\lambda}^{k}}{\lambda_{k+1}-\lambda_k}
- \tfrac{\gamma_k}{2}\norm{\lambda_{k+1}-\lambda_k}^{2}
+ \tfrac{\tau_k}{2}\norm{x-x_k}^{2} - ({\tfrac{\tau_k}{2}+\mu_F}\norm{x-x_{k+1}}^{2} \nonumber\\
&\quad + \inner{\Delta F_{k+1}}{x_{k+1}-x}
- \theta_k \inner{\Delta F_k}{x_k-x}
- \theta_k \inner{\Delta F_k}{x_{k+1}-x_k}
- \tfrac{\tau_k}{2}\norm{x_k-x_{k+1}}^{2} \nonumber\\
&\quad + \inner{\lambda}{ g({x_{k+1},y_{k+1}}-\ell_g({x_{k+1}; x_k, y_{k+1}} }.\label{eq:4.24}
\end{align}
\end{lemma}
\begin{sharedproof}
 Using the optimality condition of $x_{k+1}$ and Lemma \ref{lemma1},
\begin{equation} \label{eq:4.10}
\inner{w_k}{x_{k+1}-x}
\le \tfrac{\tau_k}{2}({\norm{x-x_k}^{2}-\norm{x_{k+1}-x_k}^{2}-\norm{x-x_{k+1}}^{2}}).
\end{equation}
By the $\mu_F$-strong monotonicity of $F$, we have
\begin{equation}\label{eq:4.11}
\inner{F({x_{k+1}}}{x_{k+1}-x}
\ge \inner{F({x}}{x_{k+1}-x} + \mu_F \norm{x_{k+1}-x}^{2}.
\end{equation}
Let \(\Delta F_k:=F(x_k)-F(x_{k-1})\). Then, $w_k = F({x_k} + \theta_k\sbra{F({x_k}-F({x_{k-1}}} +g'_{x}({x_k, y_{k+1}} \lambda_{k+1} = F(x_{k+1}) - \Delta F_{k+1} + \theta_k \Delta F_k + g'_{x}({x_k, y_{k+1}} \lambda_{k+1}$.
Using this relation in \eqref{eq:4.10}, then adding \eqref{eq:4.11}, and using the identity
\(x_{k+1}-x=(x_k-x)+(x_{k+1}-x_k)\), 
we obtain
\begin{align}
&\inner{F(x)}{x_{k+1}-x}
 +\inner{g'_x(x_k,y_{k+1})\lambda_{k+1}}{x_{k+1}-x}
 \nonumber\le \tfrac{\tau_k}{2}\norm{x-x_k}^{2}
 -(\tfrac{\tau_k}{2}+\mu_F)\norm{x-x_{k+1}}^{2}\nonumber\\
& +\inner{\Delta F_{k+1}}{x_{k+1}-x}
 -\theta_k\inner{\Delta F_k}{x_k-x}
 -\theta_k\inner{\Delta F_k}{x_{k+1}-x_k}
 -\tfrac{\tau_k}{2}\norm{x_{k+1}-x_k}^{2}.
\label{eq:4.14}
\end{align}
Applying Lemma \ref{lemma1} to the multiplier $\lambda_{k+1}$-update, noting the definition of $q^k_\lambda$ and using
\(s_k=\ell_g(x_k;x_{k-1},y_{k+1})+\theta_k q_\lambda^k\)
, we obtain
\begin{align}
&-\inner{\ell_g(x_{k+1};x_k,y_{k+1})}{\lambda_{k+1}-\lambda}
 +\inner{q_\lambda^{k+1}}{\lambda_{k+1}-\lambda}
\le \tfrac{\gamma_k}{2}
 (\norm{\lambda-\lambda_k}^{2}
  -\norm{\lambda-\lambda_{k+1}}^{2}
  -\norm{\lambda_{k+1}-\lambda_k}^{2})
 \nonumber\\
&+\theta_k\inner{q_\lambda^k}{\lambda_k-\lambda}
 +\theta_k\inner{q_\lambda^k}{\lambda_{k+1}-\lambda_k}.
\label{eq:4.18}
\end{align}
By convexity of \(g_i(\cdot,y_{k+1,i})\) for each \(i\) and
\(\lambda_{k+1}\geq0\),
\(\inner{g'_x(x_k,y_{k+1})\lambda_{k+1}}{x_{k+1}-x}
\geq
\inner{\lambda_{k+1}}
{\ell_g(x_{k+1};x_k,y_{k+1})-g(x,y_{k+1})}\).
Combining this inequality with \eqref{eq:4.14} and
\eqref{eq:4.18} gives
\begin{align}
&\inner{F(x)}{x_{k+1}-x}
 -\inner{g(x,y_{k+1})}{\lambda_{k+1}}
 +\inner{\ell_g(x_{k+1};x_k,y_{k+1})}{\lambda}
 \nonumber\\
&\le \tfrac{\gamma_k}{2}
 (\norm{\lambda-\lambda_k}^{2}
  -\norm{\lambda-\lambda_{k+1}}^{2}
  -\norm{\lambda_{k+1}-\lambda_k}^{2})
 -\inner{q_\lambda^{k+1}}{\lambda_{k+1}-\lambda}
 \nonumber\\
&+\theta_k\inner{q_\lambda^k}{\lambda_k-\lambda}
 +\theta_k\inner{q_\lambda^k}{\lambda_{k+1}-\lambda_k}
 +\tfrac{\tau_k}{2}\norm{x-x_k}^{2}
 -(\tfrac{\tau_k}{2}+\mu_F)\norm{x-x_{k+1}}^{2}
 \nonumber\\
&+\inner{\Delta F_{k+1}}{x_{k+1}-x}
 -\theta_k\inner{\Delta F_k}{x_k-x}
 -\theta_k\inner{\Delta F_k}{x_{k+1}-x_k}
 -\tfrac{\tau_k}{2}\norm{x_{k+1}-x_k}^{2}.
\label{eq:4.22}
\end{align}
 Recalling $Q_1$ and $Q_2$,
\begin{equation} \label{eq:4.23}
Q_1({z_{k+1},z}+Q_2({z_{k+1},z}
= \inner{F({x}}{x_{k+1}-x}
- \inner{g({x,y_{k+1}}}{\lambda_{k+1}}
+ \inner{g({x_{k+1},y_{k+1}}}{\lambda}.
\end{equation}
Substituting \eqref{eq:4.23} into \eqref{eq:4.22} yields
\eqref{eq:4.24} and conclude the proof.
\end{sharedproof}
Before we proceed to bound on $Q_0$, we define the $y$-gradient of mapping $g$ as
\begin{align*}
g'_{y}({x, y}
&:=
[ g'_{1,y}(x, y^{1})^\top,\ \dots,\ g'_{m,y}(x, y^{m})^\top ]^\top
\in \mbb{R}^{N},
\qquad N := \tsum_{i=1}^{m} n_i .
\end{align*}
Then, $y$-\emph{momentum  vector} is $q_{y}^{k}:=
g'_{y}({x_k, y_k}
-
g'_{y}({x_{k-1}, y_{k-1}}
=
[
(q_{y}^{k,1})^\top,\dots,(q_{y}^{k,m})^\top 
]^\top$,
where 
$q_{y}^{k,i} := g'_{i,y}(x_k, y_k^{i}) - g'_{i,y}(x_{k-1}, y_{k-1}^{i})$.
\begin{lemma}\label{lemma3}
Recall $\Lambda$ from \eqref{eq:Lambda_def}. Then, for any $z=(x,\lambda,y) \in X \times \Rbb^m_+ \times Y$, we have 
\begin{equation} \label{eq:4.30}
\begin{aligned}
Q_0({z_{k+1},z}
&\le -\inner{\Lambda q_y^{k+1}}{y_{k+1}-y}
+ \theta_k \inner{\Lambda q_y^k}{y_k-y}
+ \theta_k \inner{\Lambda q_y^k}{y_{k+1}-y_k} \\
&\quad - \tfrac{\sigma_k}{2}\norm{y_{k+1}-y_k}_{\lambda}^{2}
+ \tfrac{\sigma_k}{2}\norm{y-y_k}_{\lambda}^{2}
- \tfrac{\sigma_k}{2}\norm{y-y_{k+1}}_{\lambda}^{2}
- \tfrac{\mu_y}{2}\norm{y_{k+1}-y}_{\lambda}^{2}.
\end{aligned}
\end{equation}
\end{lemma}
\begin{sharedproof}
 Using the optimality condition of $y_{k+1}$ and the three–point lemma,
\begin{equation}\label{eq:4.25}
\inner{-u_k^i}{y^i_{k+1}-y^i}
\le \tfrac{\sigma_k}{2}\norm{y^i-y_k^i}^{2}
- \tfrac{\sigma_k}{2}\norm{y^i_{k+1}-y_k^i}^{2}
- \tfrac{\sigma_k}{2}\norm{y^i-y^i_{k+1}}^{2}.
\end{equation}
 Using the definition of $q_y^{k,i}$ inside \eqref{eq:4.25}, we have 
\begin{equation} \label{eq:4.26}
\begin{aligned}
\inner{- g'_{i,y}(x_{k+1},y_{k+1}^i)}{y_{k+1}^i-y^i}
&\le -\inner{q_y^{k+1,i}}{y_{k+1}^i-y^i}
+ \theta_k \inner{q_y^{k,i}}{y_k^i-y^i}
+ \theta_k \inner{ q_y^{k,i}}{y_{k+1}^i-y_k^i} \\
&\quad - \tfrac{\sigma_k}{2}\norm{y_{k+1}^i-y_k^i}^{2}
+ \tfrac{\sigma_k}{2}\norm{y^i-y_k^i}^{2}
- \tfrac{\sigma_k}{2}\norm{y^i-y_{k+1}^i}^{2}.
\end{aligned}
\end{equation}
 By $\mu_y$–strong concavity of $g_i$ in $y^i$,
\begin{equation}\label{eq:4.28}
g_i({x_{k+1},y^i}-g_i({x_{k+1},y_{k+1}^i}
\le -\inner{g'_{i,y}({x_{k+1},y_{k+1}^i}}{y_{k+1}^i-y^i}
- \tfrac{\mu_y}{2}\norm{y_{k+1}^i-y^i}^{2}.
\end{equation}
Summing \eqref{eq:4.26} and \eqref{eq:4.28}, multiplying resulting relation by by $\lambda_i$, summing over $i \in [m]$, using the definition of $q_y^k$, and recalling $\Lambda$ and $\norm{\cdot}_\lambda$, we get
relation \eqref{eq:4.30} and conclude the proof.
\end{sharedproof}
\begin{lemma}
Summing \eqref{eq:4.24} and \eqref{eq:4.30}, and recalling
\eqref{eq:4(a)}, we obtain
\begin{align}
Q(z_{k+1},z)
&\le \tfrac{\gamma_k}{2}
 (\norm{\lambda-\lambda_k}^{2}
 -\norm{\lambda-\lambda_{k+1}}^{2})
 -\inner{q_\lambda^{k+1}}{\lambda_{k+1}-\lambda}
 +\theta_k\inner{q_\lambda^k}{\lambda_k-\lambda}
 \nonumber\\
&+\theta_k\inner{q_\lambda^k}{\lambda_{k+1}-\lambda_k}
 -\tfrac{\gamma_k}{2}\norm{\lambda_k-\lambda_{k+1}}^{2}
 +\tfrac{\tau_k}{2}\norm{x-x_k}^{2}
 -(\tfrac{\tau_k}{2}+\mu_F)\norm{x-x_{k+1}}^{2}
 \nonumber\\
&+\inner{\Delta F_{k+1}}{x_{k+1}-x}
 -\theta_k\inner{\Delta F_k}{x_k-x}
 -\theta_k\inner{\Delta F_k}{x_{k+1}-x_k}
 -\tfrac{\tau_k}{2}\norm{x_k-x_{k+1}}^{2}
 \nonumber\\
&+\inner{\lambda}
 {g(x_{k+1},y_{k+1})-\ell_g(x_{k+1};x_k,y_{k+1})}
 -\inner{\Lambda q_y^{k+1}}{y_{k+1}-y}
 \nonumber\\
&+\theta_k\inner{\Lambda q_y^k}{y_k-y}
 +\theta_k\inner{\Lambda q_y^k}{y_{k+1}-y_k}
 -\tfrac{\sigma_k}{2}\norm{y_{k+1}-y_k}_\lambda^{2}
 +\tfrac{\sigma_k}{2}\norm{y-y_k}_\lambda^{2}
 \nonumber\\
&-\tfrac{\sigma_k}{2}\norm{y-y_{k+1}}_\lambda^{2}
 -\tfrac{\mu_y}{2}\norm{y_{k+1}-y}_\lambda^{2}.
\label{eq:4.31}
\end{align}
\end{lemma}
Considering the upper bound on the optimality gap of $\bar{x}_K, K=i+1$ ,where $\bar{x}_K$ is averaged iterate $ \wb{x}_K = \tfrac{\tsum_{k=0}^{K-1} t_k  x_{k+1}}{\tsum_{k=0}^{K-1} t_k}$.Then the gap function we derived before in \eqref{eq:4(a)} yields the averaged inequality
\begin{align}\label{eq:4.32}
\langle F(x), \bar{x}_{K}-x\rangle
+ \lambda^\top g(\bar{x}_{K},y)
- \tfrac{1}{\tsum_{k=0}^{K-1} t_k}\tsum_{k=0}^{K-1} t_k \ \lambda_{k+1}^\top g(x,y_{k+1})
\le
\tfrac{1}{\tsum_{k=0}^{K-1} t_k}\tsum_{k=0}^{K-1} t_k \ Q(z_{k+1},(x,\lambda,y)),
\end{align}
$\forall x \in X,\ \lambda \ge 0,\ y \in Y$ ,
and choosing $\lambda=\lambda_0=0$ and $y=y_0 \in Y$, for any $x \in X$ we obtain the upper bound of the optimality gap
\begin{gather} \label{eq:4.33}
\forall x \in \tilde{X}; 
\langle F(x), \bar{x}_K - x \rangle
\leq \tfrac{1}{\tsum_{k=0}^{K-1} t_k} 
\tsum_{k=0}^{K-1} t_k Q(z_{k+1},\hat{z})), \quad 
\hat{z}= \begin{bmatrix}
x &
0 &
y_0
\end{bmatrix}^T
\end{gather}
where $\tilde{X} := X \cap \left\{ x : g_i(x,y) \leq 0,\ i = 1, \dots, m \right\}$\\
By Assumption \ref{ass:opt_sol} (existence of multipliers for VI-SIP), there exists 
$x^{*} \in \tilde{X}$ and $\lambda^{*} \geq 0$ with complementarity $\lambda_i^{*} g_i^{*}(x^{*}) = 0$  and feasibility $g^{*}(x^{*}) \leq 0$. Then, we define
\begin{equation} \label{eq:4.34}
\hat{\lambda}_K 
:= \textstyle{\arg\max_\lambda}  \{ (\lambda - \lambda^*)^\top g^*(\bar{x}_K) | \lambda \geq \lambda^*,  \|\lambda - \lambda^*\| \leq 1 \}
\end{equation}
and $\hat{y}_K := (\hat{y}_{K,1}^\top, \ldots, \hat{y}_{K,m}^\top)^\top \in Y$ where 
$\hat{y}_K \in \arg\max_{y \in Y} g(\bar{x}_K, y)$.
Then, we have
\begin{align}
(\hat{\lambda}_K)^\top g^*(\bar{x}_K) 
&= \|[g^*(\bar{x}_K)]^+\|_\infty + (\lambda^*)^\top g^*(\bar{x}_K)\label{eq:4.35},\\
g(\bar{x}_K, \hat{y}_K) &= g^*(\bar{x}_K) \label{eq:4.36}.
\end{align}
Using \eqref{eq:4.35},\eqref{eq:4.36}, feasibility of $x^*$, and applying the averaged gap inequality \eqref{eq:4.32} with the choice $(x,\lambda,y)=(x^*,\hat{\lambda}_K,\hat{y}_K)$, we obtain the following upper bound on the infeasibility of $\bar{x}_K$:
\begin{align}\label{eq:4.37}
&\|[g^*(\bar{x}_K)]^+\|_\infty 
= (\hat{\lambda}_K)^\top g^*(\bar{x}_K) - (\lambda^*)^\top g^*(\bar{x}_K) \nonumber\\
&\leq \langle F(x^*), \bar{x}_K - x^* \rangle 
+ \hat{\lambda}_K^\top g(\bar{x}_K, \hat{y}_K) 
- \tfrac{\tsum_{k=0}^{K-1} t_k \lambda_{k+1}^\top g(x^*)}{\tsum_{k=0}^{K-1} t_k} \leq \tfrac{1}{\tsum_{k=0}^{K-1} t_k} 
\tsum_{k=0}^{K-1} t_k Q(z_{k+1}, (x^*, \hat{\lambda}_K, \hat{y}_K))
\end{align}

\subsubsection{Step-size policies for various settings}
In this section, we present our main theoretical results on the convergence rate of the algorithm. First, we prove the upper bound of the gap function over $k=1,...,{K-1}$ and derive the convergence rates  of the optimality gap and infeasibility of $\bar{x}_K$.
\begin{proposition}
 Let $\{t_k\}_{k \ge 0}$ be a nondecreasing sequence. Suppose $\{t_k, \theta_k, \tau_k, \sigma_k,\gamma_k\}_{k \ge 1}$ in Algorithm \ref{alg:vi-sip-agsip} satisfies the following conditions:
\begin{gather}
t_k \tau_k \leq t_{k-1} [\tau_{k-1} + 2 \mu_F], \quad 
t_k \gamma_k \leq t_{k-1} \gamma_{k-1}, \quad
t_k \sigma_k \leq t_{k-1} (\sigma_{k-1}+\mu_y), \quad
t_k \theta_k = t_{k-1}, \label{eq:4.38}\\
\tfrac{6 t_k M_g^2}{\gamma_{k+1}} + t_k L_F + \tfrac{t_{k-1} L_F}{2} 
+ \tfrac{24 t_{k+1} M_g^2}{\gamma_{k+2}} 
- \tfrac{7 t_k \tau_k}{32} \leq 0, \label{eq:4.39}\\
\tfrac{8 t_k (L_g^{yy})^2}{\sigma_{k+1}} 
+ \tfrac{3 t_{k-1} \sigma_k}{16} 
- \tfrac{t_k \sigma_k}{2} \leq 0. \label{eq:4.41}
\end{gather}
Define $\Gamma_{i+1} := \sum_{k=0}^i t_k$. Then, for any $x\in X$, $y\in Y$, and $\lambda\in\mathbb{R}_{+}^{m}$ satisfying
$\|\lambda\|_{1}\leq \Lambda$, we have
\label{eq:4.62}
\begin{align}
&\Gamma_K\inner{F(x)}{\bar{x}_K-x}
-\Gamma_K\inner{g(\bar{x}_K,y)}{\lambda}
-\tsum_{k=0}^{K-1}t_k\inner{g(x,y_{k+1})}{\lambda_{k+1}}
\leq
\tfrac{t_0\tau_0}{2}\norm{x-x_0}^{2}
+
\tfrac{t_0\gamma_0}{2}\norm{\lambda-\lambda_0}^{2}
\nonumber\\&+
\tfrac{t_0\sigma_0}{2}\norm{y-y_0}_{\lambda}^{2}
+
\Psi(\Lambda,K)
+
\tsum_{k=0}^{K-1}
\tfrac{4t_k(H_g^x\Lambda)^2}{\tau_k}
+
4(
\tsum_{k=1}^{K-1}\tfrac{t_k}{\sigma_{k+1}}
+
\tsum_{k=1}^{K-1}\tfrac{t_{k-1}}{\sigma_k})
(H_g^y)^2\Lambda .
\label{eq:T-Lambda-K}\nonumber\\
&\text{where }
\Psi(\Lambda,K)
:=
\tsum_{k=1}^{K-1}
[
t_kL_{xx}^{g}\Lambda
+
\tfrac{8t_k\Lambda(L_g^{yx})^2}{\sigma_{k+1}}
-
\tfrac{7t_k\tau_k}{32}
]
\norm{x_{k+1}-x_k}^{2}.
\end{align}
\end{proposition}
\begin{sharedproof}
Recall that, for each \(i\in[m]\),
\begin{equation}\label{eq:4.42}
g_i(x,y_i)-g_i(x',y_i)
-\inner{g'_{i,x}(x',y_i)}{x-x'}
\le L_g^{xx}\norm{x-x'}^{2}
   +H_g^x\norm{x-x'}.
\end{equation}
Setting \(x'=x_k\), \(x=x_{k+1}\), and \(y_i=y_{i,k+1}\)
in \eqref{eq:4.42}, then multiplying by \(\lambda_i\geq0\)
and summing over \(i\), gives
\begin{equation}\label{eq:4.43}
\inner{\lambda}
{g(x_{k+1},y_{k+1})
 -\ell_g(x_{k+1};x_k,y_{k+1})}
\le L_g^{xx}\norm{\lambda}_1\norm{x_{k+1}-x_k}^{2}
   +H_g^x\norm{\lambda}_1\norm{x_{k+1}-x_k}.
\end{equation}
Substituting \eqref{eq:4.43} into \eqref{eq:4.31} yields
\begin{align}
Q(z_{k+1},z)
&\le \tfrac{\gamma_k}{2}
 (\norm{\lambda-\lambda_k}^{2}
 -\norm{\lambda-\lambda_{k+1}}^{2})
 -\inner{q_\lambda^{k+1}}{\lambda_{k+1}-\lambda}
 +\theta_k\inner{q_\lambda^k}{\lambda_k-\lambda}
 \nonumber\\
&+\theta_k\inner{q_\lambda^k}{\lambda_{k+1}-\lambda_k}
 -\tfrac{\gamma_k}{2}\norm{\lambda_k-\lambda_{k+1}}^{2}
 +\tfrac{\tau_k}{2}\norm{x-x_k}^{2}
 -(\tfrac{\tau_k}{2}+\mu_F)\norm{x-x_{k+1}}^{2}
 \nonumber\\
&+\inner{\Delta F_{k+1}}{x_{k+1}-x}
 -\theta_k\inner{\Delta F_k}{x_k-x}
 -\theta_k\inner{\Delta F_k}{x_{k+1}-x_k}
 -\tfrac{\tau_k}{2}\norm{x_k-x_{k+1}}^{2}
 \nonumber\\
&+L_g^{xx}\norm{\lambda}_1\norm{x_{k+1}-x_k}^{2}
 +H_g^x\norm{\lambda}_1\norm{x_{k+1}-x_k}
 -\inner{\Lambda q_y^{k+1}}{y_{k+1}-y}
 \nonumber\\
&+\theta_k\inner{\Lambda q_y^k}{y_k-y}
 +\theta_k\inner{\Lambda q_y^k}{y_{k+1}-y_k}
 -\tfrac{\sigma_k}{2}\norm{y_{k+1}-y_k}_\lambda^{2}
 \nonumber\\
&+\tfrac{\sigma_k}{2}\norm{y-y_k}_\lambda^{2}
 -\tfrac{\sigma_k+\mu_y}{2}\norm{y-y_{k+1}}_\lambda^{2}.
\label{eq:4.44}
\end{align}
 Multiplying \eqref{eq:4.44} by $t_k$, 
summing the resulting relation from $k=0$ to $k=i$, then using the stepsize condition \eqref{eq:4.38}, applying Jensen’s inequality for the convex functions $\inner{g({\cdot,y}}{\lambda}$, and the fact that $q_\lambda^{0}=0$, $\Delta F_{0}=0$, $q_y^{0}=0$, we obtain
\begin{align}\label{eq:4.46}
&\Gamma_{i+1}\inner{F({x}}{\bar{x}_{i+1}-x}
- \Gamma_{i+1}\inner{g({\bar{x}_{i+1},y}}{\lambda}
- \tsum_{k=0}^{i} t_k \inner{g({x,y_{k+1}}}{\lambda_{k+1}} \nonumber\\
&\le \tfrac{t_0\tau_0}{2}\norm{x-x_0}^{2}
+ \tfrac{t_0\gamma_0}{2}\norm{\lambda-\lambda_0}^{2}
+ \tfrac{t_0\sigma_0}{2}\norm{y-y_0}_{\lambda}^{2}
+ t_i \inner{\Delta F_{i+1}}{x_{i+1}-x}
- t_i \inner{\Lambda q_y^{i+1}}{y_{i+1}-y} \nonumber\\
&-t_i\inner{q_\lambda^{i+1}}{\lambda_{i+1}-\lambda}
-t_i(\tfrac{\tau_i}{2}+\mu_F)\norm{x-x_{i+1}}^{2}
-\tfrac{t_i\gamma_i}{2}\norm{\lambda-\lambda_{i+1}}^{2}
-\tfrac{t_i(\sigma_i+\mu_y)}{2}\norm{y-y_{i+1}}_{\lambda}^{2}
\nonumber\\
& + \tsum_{k=1}^{i}
[t_{k-1}\inner{q_\lambda^{k}}{\lambda_{k+1}-\lambda_k}
-t_{k-1}\inner{\Delta F_k}{x_{k+1}-x_k}
+t_{k-1}\inner{\Lambda q_y^{k}}{y_{k+1}-y_k}]\nonumber \\
&+
\tsum_{k=0}^{i}
[-\tfrac{t_k\tau_k}{2}\norm{x_k-x_{k+1}}^{2}
-\tfrac{t_k\gamma_k}{2}\norm{\lambda_{k+1}-\lambda_k}^{2}
-\tfrac{t_k\sigma_k}{2}\norm{y_{k+1}-y_k}_{\lambda}^{2}]
\nonumber\\
& + 
\tsum_{k=0}^{i}
[t_k L_{xx}^{g}\norm{\lambda}_{1}\norm{x_{k+1}-x_k}^{2}
+t_k H_g^{x}\norm{\lambda}_{1}\norm{x_{k+1}-x_k}]
,\end{align}
where $\bar{x}_{i+1}=\tfrac{1}{\Gamma_{i+1}}\tsum_{k=0}^{i} t_k x_{k+1}$ and recall $\Gamma_{i+1}=\tsum_{k=0}^{i} t_k$.\\
Recall that $q_\lambda^{k}=\ell_g({x_k;x_{k-1},y_k}-\ell_g({x_{k-1};x_{k-2},y_k}$. Using the definition $\ell_g({x;x',y}=g({x',y}+g'_{x}({x',y}({x-x'}$,
\begin{align}
\norm{q_\lambda^k}
&=\norm{
 g(x_{k-1},y_k)
 +g'_x(x_{k-1},y_k)^\top(x_k-x_{k-1})
 -g(x_{k-2},y_k)
 -g'_x(x_{k-2},y_k)^\top(x_{k-1}-x_{k-2}) }
 \nonumber\\
&\le \norm{g(x_{k-1},y_k)-g(x_{k-2},y_k)}
 +\norm{g'_x(x_{k-1},y_k)^\top(x_k-x_{k-1})} +
 \norm{g'_x(x_{k-2},y_k)^\top(x_{k-1}-x_{k-2})}
 \nonumber\\
&\le 2M_g\norm{x_{k-1}-x_{k-2}}
 +M_g\norm{x_k-x_{k-1}}.
\label{eq:4.47}
\end{align}
Using Young’s inequality, \eqref{eq:4.47} and the fact that $(a+b)^2 \le 2a^2+2b^2$ for all $a,b \in \Rbb$, we have
\begin{align}\label{eq:4.48}
t_{k-1}\inner{q_\lambda^{k}}{\lambda_{k+1}-\lambda_k}
&\le \tfrac{t_{k-1}\gamma_k}{8}\norm{\lambda_{k+1}-\lambda_k}^{2}
+ \tfrac{2t_{k-1}}{\gamma_k}\norm{q_\lambda^{k}}^{2} \nonumber\\
&\le \tfrac{t_{k-1}\gamma_k}{8}\norm{\lambda_{k+1}-\lambda_k}^{2}
+ \tfrac{4 t_{k-1} M_g^{2}}{\gamma_k}\norm{x_k-x_{k-1}}^{2}
+ \tfrac{16 t_{k-1} M_g^{2}}{\gamma_k}\norm{x_{k-1}-x_{k-2}}^{2}.
\end{align}
Using the definition of $\Delta F_k$ and $L_F$–Lipschitz continuity of operator $F$, we have
\begin{equation}\label{eq:4.49}
\begin{aligned}
- t_{k-1}\inner{\Delta F_k}{x_{k+1}-x_k}
&\le t_{k-1}\norm{\Delta F_k} \norm{x_{k+1}-x_k}
\le t_{k-1} L_F \norm{x_k-x_{k-1}} \norm{x_{k+1}-x_k} \\
&\le \tfrac{t_{k-1}L_F}{2}\norm{x_k-x_{k-1}}^{2}
+ \tfrac{t_{k-1}L_F}{2}\norm{x_{k+1}-x_k}^{2}.
\end{aligned}
\end{equation}
Using the Lipschitz properties, we have
\begin{equation}\label{eq:4.50}
\begin{aligned}
\norm{g'_{i,y}({x_k,y_k^{i}}-g'_{i,y}({x_{k-1},y_{k-1}^{i}}}
&\le L_g^{yx}\norm{x_k-x_{k-1}}
+ L_g^{yy}\norm{y_k^{i}-y_{k-1}^{i}} + H_g^{y}
\end{aligned}
\end{equation}
for any $x_k,x_{k-1} \in X$ and $y_k^i, y_{k-1}^i \in Y^i$.
Recalling $q_{y}^k = g'_{y}({x_k,y_k} - g'_{y}({x_{k-1},y_{k-1}}$
\begin{align}\label{eq:4.51}
&\inner{\Lambda q_y^{k}}{y_{k+1}-y_k}
= \tsum_{i=1}^{m}\lambda_i \inner{g'_{i,y}({x_k,y_k^{i}}-g'_{i,y}({x_{k-1},y_{k-1}^{i}}}{y_{k+1}^{i}-y_k^{i}} \nonumber\\
&\le \tsum_{i=1}^{m}\lambda_i [
L_g^{yx}\norm{x_k-x_{k-1}}\norm{y_{k+1}^{i}-y_k^{i}}
+ L_g^{yy}\norm{y_k^{i}-y_{k-1}^{i}}\norm{y_{k+1}^{i}-y_k^{i}}
+ H_g^{y}\norm{y_{k+1}^{i}-y_k^{i}}] \nonumber\\
&\le \tsum_{i=1}^{m}\lambda_i [
\tfrac{4}{\sigma_k}({L_g^{yx}\norm{x_k-x_{k-1}}}^{2}
+ \tfrac{\sigma_k}{8}\norm{y_{k+1}^{i}-y_k^{i}}^{2}
+ \tfrac{4}{\sigma_k}({L_g^{yy}\norm{y_k^{i}-y_{k-1}^{i}}}^{2}
+ H_g^{y}\norm{y_{k+1}^{i}-y_k^{i}}].
\end{align}
By Cauchy–Schwarz Inequality, 
$ \tsum_{i=1}^{m}\lambda_i \norm{y_{k+1}^{i}-y_k^{i}}
\le \norm{\lambda}_{1}^{1/2} \norm{y_{k+1}-y_k}_{\lambda}.$
Therefore combining \eqref{eq:4.50}–\eqref{eq:4.51} we have
\begin{equation}\label{eq:4.53}
\begin{aligned}
t_{k-1}\inner{\Lambda q_y^{k}}{y_{k+1}-y_k}
&\le \tfrac{4t_{k-1}}{\sigma_k}\norm{\lambda}_{1}({L_g^{yx}}^{2}\norm{x_k-x_{k-1}}^{2}
+ \tfrac{4t_{k-1}}{\sigma_k}({L_g^{yy}}^{2}\norm{y_k-y_{k-1}}_{\lambda}^{2} \\
&\quad + \tfrac{t_{k-1}\sigma_k}{8}\norm{y_{k+1}-y_k}_{\lambda}^{2}
+ t_{k-1} H_g^{y} \norm{\lambda}_{1}^{1/2} \norm{y_{k+1}-y_k}_{\lambda}.
\end{aligned}
\end{equation}
Substituting \eqref{eq:4.48}, \eqref{eq:4.49}, and
\eqref{eq:4.53} into \eqref{eq:4.46} and combining the
terms, we obtain
\begin{align}
&\tsum_{k=1}^{i}
 [t_{k-1}\inner{q_\lambda^k}{\lambda_{k+1}-\lambda_k}
 -t_{k-1}\inner{\Delta F_k}{x_{k+1}-x_k}
 +t_{k-1}\inner{\Lambda q_y^k}{y_{k+1}-y_k}]
 \nonumber\\
&\le \tsum_{k=1}^{i}
 [\tfrac{t_{k-1}\gamma_k}{8}\norm{\lambda_{k+1}-\lambda_k}^2
 +\tfrac{4t_{k-1}M_g^2}{\gamma_k}\norm{x_k-x_{k-1}}^2
 +\tfrac{16t_{k-1}M_g^2}{\gamma_k}\norm{x_{k-1}-x_{k-2}}^2
 \nonumber\\
&+\tfrac{t_{k-1}L_F}{2}\norm{x_k-x_{k-1}}^2
 +\tfrac{t_{k-1}L_F}{2}\norm{x_{k+1}-x_k}^2
 +\tfrac{4t_{k-1}(L_g^{yx})^2}{\sigma_k}
   \norm{\lambda}_1\norm{x_k-x_{k-1}}^2
 \nonumber\\
&+\tfrac{4t_{k-1}(L_g^{yy})^2}{\sigma_k}
   \norm{y_k-y_{k-1}}_\lambda^2
 +\tfrac{t_{k-1}\sigma_k}{8}
   \norm{y_{k+1}-y_k}_\lambda^2
 +t_{k-1}H_g^y\norm{\lambda}_1^{1/2}
   \norm{y_{k+1}-y_k}_\lambda]
\label{eq:4.54}
\end{align}
By Young’s inequality, for each $k$, and summing over the corresponding indices, we obtain
\begin{align}\label{eq:4.56}
&\tsum_{k=0}^{i} t_k H_g^{x}\norm{\lambda}_{1}\norm{x_{k+1}-x_k}
+ \tsum_{k=1}^{i} t_{k-1} H_g^{y}\norm{\lambda}_{1}^{1/2}\norm{y_{k+1}-y_k}_{\lambda} \nonumber\\ 
&\le \tsum_{k=0}^{i}
[\tfrac{t_k(4H_g^{x}\norm{\lambda}_{1})^{2}}{\tau_k}
+\tfrac{t_k\tau_k}{16}\norm{x_{k+1}-x_k}^{2}]
+\tsum_{k=1}^{i}
[\tfrac{4t_{k-1}(H_g^{y})^{2}\norm{\lambda}_{1}}{\sigma_k}
+\tfrac{t_{k-1}\sigma_k}{16}\norm{y_{k+1}-y_k}_{\lambda}^{2}].
\end{align}
 Using the same logic we used before in \eqref{eq:4.48}–\eqref{eq:4.53} we obtain,
\begin{gather}
\label{eq:4.58}
\begin{aligned}
t_i \inner{\Delta F_{i+1}}{x_{i+1}-x}
&\le
\tfrac{t_i L_F}{2}\norm{x_{i+1}-x_i}^{2}
+
\tfrac{t_i L_F}{2}\norm{x_{i+1}-x}^{2}.
\end{aligned}
\\
\label{eq:4.59}
\begin{aligned}
- t_i \inner{q_\lambda^{i+1}}{\lambda_{i+1}-\lambda}
&\le
\tfrac{t_i \gamma_{i+1}}{4}\norm{\lambda_{i+1}-\lambda}^{2}
+
\tfrac{8 t_i M_g^{2}}{\gamma_{i+1}}\norm{x_i-x_{i-1}}^{2}
+
\tfrac{2 t_i M_g^{2}}{\gamma_{i+1}}\norm{x_{i+1}-x_i}^{2}.
\end{aligned}
\\
\label{eq:4.60}
\begin{aligned}
- t_i \inner{\Lambda q_y^{i+1}}{y_{i+1}-y}
&\le
\tfrac{4 t_i}{\sigma_{i+1}}\norm{\lambda}_{1}(L_g^{yx})^{2}
\norm{x_{i+1}-x_i}^{2}
+
\tfrac{4 t_i}{\sigma_{i+1}}(L_g^{yy})^{2}
\norm{y_{i+1}-y_i}_{\lambda}^{2}
\\
&
+
\tfrac{t_i\sigma_{i+1}}{8}\norm{y_{i+1}-y}_{\lambda}^{2}
+
t_i H_g^{y}\norm{y_{i+1}-y}_{\lambda}
\\
&\le
\tfrac{4 t_i}{\sigma_{i+1}}\norm{\lambda}_{1}(L_g^{yx})^{2}
\norm{x_{i+1}-x_i}^{2}
+
\tfrac{4 t_i}{\sigma_{i+1}}(L_g^{yy})^{2}
\norm{y_{i+1}-y_i}_{\lambda}^{2}
\\
&
+
\tfrac{3 t_i\sigma_{i+1}}{16}\norm{y_{i+1}-y}_{\lambda}^{2}
+
\tfrac{4 t_i (H_g^{y})^{2}\norm{\lambda}_{1}}{\sigma_{i+1}}.
\end{aligned}
\end{gather}
Substituting \eqref{eq:4.54}-\eqref{eq:4.60} into \eqref{eq:4.46}, and rearranging, the final bound can be obtained as:
\begin{align}
& \Gamma_{i+1}\inner{F(x)}{\bar{x}_{i+1}-x}
- \Gamma_{i+1}\inner{g(\bar{x}_{i+1},y)}{\lambda}
- \tsum_{k=0}^{i} t_k \inner{g(x,y_{k+1})}{\lambda_{k+1}}
\nonumber\\
& \leq \tfrac{t_0 \tau_0}{2}\norm{x-x_0}^{2} 
+ \tfrac{t_0 \gamma_0}{2}\norm{\lambda-\lambda_0}^{2} 
+ \tfrac{t_0 \sigma_0}{2}\norm{y-y_0}_{\lambda}^{2}
\nonumber\\
&\quad+(\tfrac{t_i L_F}{2} - t_i(\tfrac{\tau_i}{2}+\mu_F))\norm{x_{i+1}-x}^{2}
+ (\tfrac{t_i \gamma_{i+1}}{4}-\tfrac{t_i \gamma_i}{2})\norm{\lambda_{i+1}-\lambda}^{2}
+(\tfrac{3 t_i \sigma_{i+1}}{16}-\tfrac{t_i(\sigma_i+\mu_y)}{2})\norm{y_{i+1}-y}_{\lambda}^{2}
\nonumber\\
&\quad
+ \tsum_{k=1}^{i} (\tfrac{t_{k-1} \gamma_k}{8}-\tfrac{t_k \gamma_k}{2})\norm{\lambda_{k+1}-\lambda_k}^{2}
+ \tsum_{k=1}^i [\tfrac{8 t_k (L_g^{yy})^2}{\sigma_{k+1}}+\tfrac{3 t_{k-1}\sigma_k}{16}-\tfrac{t_k \sigma_k}{2}]\norm{y_{k+1}-y_k}_{\lambda}^{2}
\nonumber\\
&\quad + \tsum_{k=1}^{i} (\tfrac{6 t_k M_g^2}{\gamma_{k+1}}+t_k L_F+\tfrac{t_{k-1}L_F}{2}+\tfrac{24 t_{k+1}M_g^2}{\gamma_{k+2}}-\tfrac{7 t_k \tau_k}{32})\norm{x_{k+1}-x_k}^{2}
\nonumber\\
&\quad + \tsum_{k=1}^{i} (t_k L_{xx}^g\norm{\lambda}_1+\tfrac{8 t_k}{\sigma_{k+1}}\norm{\lambda}_1(L_g^{yx})^2-\tfrac{7 t_k \tau_k}{32})\norm{x_{k+1}-x_k}^{2}
\nonumber\\
&\quad + \tsum_{k=0}^{i}\tfrac{4t_k (H_g^x \norm{\lambda}_1)^2}{\tau_k}
+4(\tsum_{k=1}^i \tfrac{t_k}{\sigma_{k+1}}+\tsum_{k=1}^i \tfrac{t_{k-1}}{\sigma_k})(H_g^y)^2\norm{\lambda}_1.
\label{eq:4.63}
\end{align}
Applying conditions \eqref{eq:4.38}, \eqref{eq:4.39}, and \eqref{eq:4.41} along with the fact that $\{t_k\}$ is a nondecreasing sequence to \eqref{eq:4.63}, and dropping the resulting non-positive terms on the right-hand side, we obtain \eqref{eq:4.62}. This concludes the proof.
\end{sharedproof}
\begin{corollary}\label{corollary1}
Suppose the parameters $t_k,\tau_k,\sigma_k,\gamma_k$ satisfy the conditions
\eqref{eq:4.38}, \eqref{eq:4.39}, \eqref{eq:4.41}, for $k=0,\ldots,K-1$.  
It holds that
\begin{align} 
\langle F(x), \bar{x}_K - x \rangle 
&\leq \tfrac{t_0 \tau_0}{2\tsum_{k=0}^{K-1} t_k}  
\| x - x_0 \|^2,
\label{eq:4.64}
\\
\| [g^*(\bar{x}_K)]_+ \|_\infty 
&\leq \tfrac{1}{\tsum_{k=0}^{K-1} t_k} [
\tfrac{t_0 \tau_0}{2} \| x - x_0 \|^2 
+ \tfrac{t_0 \gamma_0}{2} ( \| \lambda^* \|_1 + 1 )^2
+ \tfrac{t_0 \sigma_0}{2} D_y^2 ( \| \lambda^* \|_1 + 1 )
+ \Psi(\| \lambda^* \|_1 + 1,K)
] \nonumber\\
&\quad+ \tfrac{1}{\tsum_{k=0}^{K-1} t_k}
[
4 (H_g^x)^2 ( \| \lambda^* \|_1 + 1 )^2  
\tsum_{k=0}^{K-1} \tfrac{t_k}{\tau_k}
+ 8 (H_g^y)^2 ( \| \lambda^* \|_1 + 1 )
\tsum_{k=0}^{K} \tfrac{t_k}{\sigma_k}
],
\label{eq:4.65}
\end{align}
where $\Psi(\| \lambda^* \|_1 + 1,K)$ defined as $\eqref{eq:T-Lambda-K}$ with $\Lambda=\| \lambda^* \|_1 +1$.
\end{corollary}
\begin{sharedproof} Using the conditions \eqref{eq:4.38}, \eqref{eq:4.39}, and \eqref{eq:4.41}, for any 
$x \in X, y \in Y, \lambda \in \mathbb{R}_+^m$ with $\|\lambda\|_1 \leq \Lambda$ we have
\begin{equation} \label{eq:4.66}
\begin{aligned}
\tsum_{k=0}^{K-1} t_k Q(z_{k+1},z)  
&\leq  \tfrac{t_0 \tau_0}{2} \| x - x_0 \|^2 
+ \tfrac{t_0 \gamma_0}{2} \| \lambda - \lambda_0 \|^2 
+ \tfrac{t_0 \sigma_0}{2} \| y - y_0 \|^2_\Lambda+\Psi( \Lambda,K) \\
&\quad+ ( 4 (H_g^x)^2 \Lambda^2  ) \tsum_{k=0}^{K-1} \tfrac{t_k}{\tau_k} 
+ 8 (H_g^y)^2 \Lambda \tsum_{k=0}^{K} \tfrac{t_k }{\sigma_k}.
\end{aligned}
\end{equation}
Recall that $\lambda_0=0$. Applying \eqref{eq:4.66} to the right-hand side of \eqref{eq:4.37} leads to \eqref{eq:4.65} with $\Lambda= \|\lambda^*\|_1 +1$, and using the facts  
$\|\hat{\lambda}_{K}-\lambda_{0}\| = \|\hat{\lambda}_{K}\|\le \|\lambda^{*}\|_{1}+1$,  
and $\|\hat{y}_{K}-y_{0}\|_{\hat{\lambda}_{K}}^{2}\le \|\hat{\lambda}_{K}\|_{1} D_{y}^{2}$.
So, here we consider the convergence properties for different scenarios of the problem using following theorems by choosing the parameters $t_k, \theta_k, \tau_k, \sigma_k$ and $\gamma_k$.
\end{sharedproof}
\begin{theorem}
Suppose \(\mu_F > 0\), \(\mu_y=0\) and the algorithm uses parameter sequences \(\{t_k, \theta_k, \tau_k, \gamma_k, \sigma_k\}\)
\begin{equation}\label{eq:4.76}
 t_k = k + k_0 +1, \theta_k = \tfrac{k +k_0}{k +k_0 +1},  \gamma_k = \tfrac{55 M_g^2}{\mu_F t_k}, \tau_k = t_k \mu_F,  \sigma_k = \tfrac{\max \left\{\sqrt{128/5}(K+k_0) L_g^{yy},(H_g^y/D_y)K^{1/2}(K+k_0+1)
\right\}}{k+k_0+1}
\end{equation}
where $k_0 \geq \tfrac{192L_F}{\mu_F}$. Then, the following holds.
\begin{gather}
\begin{aligned}
\inner{F(x)}{\bar{x}_K-x}
&\le
\tfrac{t_0\tau_0}{2K(K+1+2k_0)}
\norm{x-x_0}^{2},
\end{aligned}
\nonumber
\\
\label{eq:4.77}
\begin{aligned}
\norm{[g^*(\bar{x}_K)]_+}_\infty
&\leq
\tfrac{(k_0+1)^2\mu_F}{2K(K+2k_0+1)}
\norm{x-x_0}^{2}
+
\tfrac{55M_g^2(\norm{\lambda^*}_1+1)^2}
{K(K+2k_0+1)\mu_F}
+
\tfrac{\sqrt{32/5}(K+k_0)L_g^{yy}D_y^2(\norm{\lambda^*}_1+1)}
{K(K+2k_0+1)}
\\
&\quad
+
\tfrac{8(H_g^x)^2(\norm{\lambda^*}_1+1)^2}
{\mu_F(K+2k_0+1)}
+
\tfrac{34H_g^yD_y(\norm{\lambda^*}_1+1)(K+k_0+1)}{\sqrt{K}{(K+2k_0+1)}}
+
\tfrac{2\Psi(\norm{\lambda^*}_1+1,K)}
{K(K+2k_0+1)} ,
\end{aligned}
\end{gather}
where $\Psi(\|\lambda^*\|_1+1,K)$ is defined in \eqref{eq:T-Lambda-K} and this term is bounded by a constant independent of $K$,i.e., $\Psi(\|\lambda^*\|_1+1,K)=O(1)$.
In addition  $\Psi(\|\lambda^*\|_1+1,K)\leq 0$ when
\begin{equation}\label{eq:4.78}
k_0 \geq \max \left\{
\tfrac{128\sqrt{5} (L_g^{yx})^2 (\|\lambda^*\|_1 + 1)}{7\sqrt{2}L_g^{yy} \mu_F} ,
\quad
\tfrac{16 L_g^{xx} (\|\lambda^*\|_1 + 1)}{\mu_F}
\right\} - 1.
\end{equation}
\end{theorem}
\begin{sharedproof}
It is easy to verify that \eqref{eq:4.38} holds. In addition,
\begin{equation*}
\begin{aligned}
\tfrac{8t_k (L_g^{yy})^2}{\sigma_{k+1}} 
&\leq \tfrac{8 (k + k_0 + 1)^2 (L_g^{yy})^2}{\sqrt{128/5}(K+k_0)L_g^{yy}}\leq \tfrac{5(\sqrt{128/5})(K+k_0)L_g^{yy}}{16}
\leq \tfrac{5 t_k \sigma_k}{16},
\end{aligned}
\end{equation*}
above equation holds \eqref{eq:4.41} where $t_{k-1}\leq t_k$, and $k\leq K-1$. Then, by definitions of $t_k, \sigma_k, \tau_k$, and $\gamma_k$, we can write
\begin{equation}\label{eq:4.80}
\begin{aligned}
\tsum_{k=0}^{K} \tfrac{t_k}{\sigma_k}
= \tsum_{k=0}^{K} \tfrac{(k+k_0+1)^2}{(H_g^y/Dy)K^{1/2}(K+k_0+1)}
\leq \tfrac{2D_yK^{1/2}(K+k_0+1)}{H_g^y},
\end{aligned}
\end{equation}

where the inequality follows because $\tsum_{k=0}^{K} (k+k_0+1)^2 \leq (K+1)(K+k_0+1)^2 \leq 2K(K+k_0+1)^2$.  
\begin{gather}
\label{eq:4.81}
\tsum_{k=0}^{K-1} t_k
=
\tfrac{K (K + 1 + 2 k_0)}{2},
\\
\label{eq:4.82}
\tsum_{k=0}^{K-1} \tfrac{t_k}{\tau_k}
=
\tfrac{K}{\mu_F},
\end{gather}
\begin{equation}\label{eq:4.83}
\begin{aligned}
\tfrac{6 t_k M_g^2}{\gamma_{k+1}} + t_k L_F + \tfrac{t_{k-1} L_F}{2} + \tfrac{24 t_{k+1} M_g^2}{\gamma_{k+2}} 
&\leq \tfrac{3t_kL_F}{2}+\tfrac{6(k+k_0+1)^2M_g^2}{55M_g^2/\mu_F} + \tfrac{24(k+k_0+2)^2M_g^2}{55M_g^2/\mu_F} \\
&\leq \tfrac{t_k\tau_k}{128}+\tfrac{30.5t_k\tau_k}{55} 
= \tfrac{7t_k\tau_k}{32},
\end{aligned}
\end{equation}
where the inequality follows because $t_{k-1}\leq t_k$, $\tau_k \geq k_0 
\mu_F$, where $k_0=\tfrac{192L_F}{\mu_F}$,  
and $k+k_0+2 \leq (1.01)(k+k_0+1)$.  
Define using \eqref{eq:T-Lambda-K} 
\begin{gather}
\label{eq:4.84}
\begin{aligned}
\zeta
&=
\tsum_{k=0}^{K-1}
g[
\tfrac{8 t_k (L_g^{yx})^2 (\norm{\lambda^*}_1 + 1)}{\sigma_{k+1}}
+
t_k L_g^{xx} (\norm{\lambda^*}_1 + 1)
-
\tfrac{7}{32} t_k \tau_k
g]
\\
&\leq
\tsum_{k=0}^{\hat{K}-1}
g[
\tfrac{8 t_k (L_g^{yx})^2 (\norm{\lambda^*}_1 + 1)}{\sigma_{k+1}}
+
t_k L_g^{xx} (\norm{\lambda^*}_1 + 1)
-
\tfrac{7}{32} t_k \tau_k
g]
\end{aligned}
\\
\label{eq:4.85}
\begin{aligned}
\text{where}\quad
\hat{K}
&\geq
[
\max
\{
\tfrac{128\sqrt{5} (L_g^{yx})^2 (\norm{\lambda^*}_1 + 1)}
{7\sqrt{2}L_g^{yy}\mu_F},
\;
\tfrac{16 L_g^{xx}(\norm{\lambda^*}_1 + 1)}
{\mu_F}
\}
-k_0-1
].
\end{aligned}
\end{gather} 
In fact when $k\geq \hat{K}$ we have,

\begin{gather}
\label{eq:4.86}
\begin{aligned}
&t_k L_g^{xx}\norm{\lambda}_1
+
\tfrac{8t_k\norm{\lambda}_1}{\sigma_{k+1}}
(L_g^{yx})^2
-
\tfrac{7t_k\tau_k}{32}
\\
&\leq
(k+k_0+1)L_g^{xx}(\norm{\lambda^*}_1+1)
+
\tfrac{
8(k+k_0+1)^2(\norm{\lambda^*}_1+1)(L_g^{yx})^2
}
{
\sqrt{128/5}(K+k_0)L_g^{yy}
}
-
\tfrac{7(k+k_0+1)^2\mu_F}{32}
\\
&\leq
(k+k_0+1)L_g^{xx}(\norm{\lambda^*}_1+1)
+
\tfrac{
8(k+k_0+1)(\norm{\lambda^*}_1+1)(L_g^{yx})^2
}
{
\sqrt{128/5}L_g^{yy}
}
-
\tfrac{7(k+k_0+1)^2\mu_F}{32}\leq 0 .
\end{aligned}
\end{gather}
where the first inequality is from the definitions of $t_k,\theta_k,\tau_k,\sigma_k$, and the second inequality is because $k\leq K-1$, and the last inequality is because $k\geq\hat{K}$, when $\hat{K}$ is defined as \eqref{eq:4.85}.  
In addition we have $\Psi(\|\lambda^*\|_1+1,K)\leq 0$ and thus \eqref{eq:4.41} $\leq 0$ when \eqref{eq:4.78} holds.  We obtain \eqref{eq:4.77} by applying \eqref{eq:4.80}, \eqref{eq:4.81}, \eqref{eq:4.82} to the right hand side of \eqref{eq:4.65}.
\end{sharedproof}

\begin{theorem}
Suppose \(\mu_F>0\), \(\mu_y > 0\), and the algorithm uses parameter sequences \(\{t_k, \theta_k, \tau_k, \gamma_k, \sigma_k\}\)
\begin{equation}\label{eq:4.95}
\begin{aligned}
    t_k &= k + k_0 + 1, \quad 
    \theta_k = \tfrac{k + k_0}{k + k_0 + 1},  \tau_k = t_k \mu_F, \quad 
    \sigma_k = \tfrac{(k + k_0) \mu_y}{2}, \quad 
    \gamma_k = \tfrac{55  M_g^2}{\mu_F t_k},
\end{aligned}
\end{equation}
where $k_0 \geq max\{\sqrt{512/5}L_g^{yy}/\mu_y, 32L_F/\mu_F\}$.\\
Then,
 \begin{equation}
\begin{aligned}
    \langle F(x), \bar{x}_K - x \rangle
    \le \tfrac{t_o\tau_0 \| x - x_0 \|^{2}}{2K(K+1+2k_0)},
\end{aligned} \nonumber
\end{equation}
\begin{equation}\label{eq:4.96}
\begin{aligned}
    \| [g^*(\bar{x}_K)]_+ \|_\infty &\leq \tfrac{(k_0+1)^2\mu_F}{2K(K+2k_0+1)}\|x-x_0\|^2+\tfrac{55M_g^2(\|\lambda^*\|_1+1)^2}{K(K+2k_0+1)\mu_F}+ \tfrac{k_0(k_0+1)\mu_yD_y^2(\|\lambda^*\|_1+1)}{K(K+2k_0+1)}\\
    &+\tfrac{8(H_g^x)^2(\|\lambda^*\|_1+1)^2}{\mu_F(K+2k_0+1)}+\tfrac{16(H_g^y)^2(\|\lambda^*\|_1+1)}{\mu_y(K+2k_0+1)}+\tfrac{2\Psi(\|\lambda^*\|_1+1,K)}{K(K+2k_0+1)}.
\end{aligned}
\end{equation}
where $\Psi(\|\lambda^*\|_1+1,K)$ is defined in \eqref{eq:T-Lambda-K} and since this term is bounded by a constant independent of $K$,then $\Psi(\|\lambda^*\|_1+1,K)=O(1)$.
In addition $\Psi(\|\lambda^*\|_1+1,K) \leq 0$ when
\begin{equation}\label{eq:4.97}
\begin{aligned}
    k_0 \geq \max \{ 
\sqrt{\tfrac{(2048) (L_g^{yx})^2 (\|\lambda^*\|_1 + 1)}{7\mu_y \mu_F}}, \quad
\tfrac{32 L_g^{xx} (\|\lambda^*\|_1 + 1)}{\mu_F} 
\} - 1.
\end{aligned}
\end{equation}
 \end{theorem}
\begin{sharedproof}
It is easy to verify that \eqref{eq:4.38} holds. In addition, 
\begin{equation*}\label{eq:4.98}
\begin{aligned}
     \tfrac{8t_k (L_g^{yy})^2}{\sigma_{k+1}} = \tfrac{8 (k + k_0 + 1)(L_g^{yy})^2}{(k+k_0+1)\mu_y} \leq\tfrac{5k_0^2\mu_y}{32}\leq \tfrac{5 t_k \sigma_k}{16},
\end{aligned}
\end{equation*}
where the first inequality is because $k_0 \geq \sqrt{512/5}L_g^{yy}/\mu_y$. This verifies \eqref{eq:4.41}.\\
Furthermore, using \eqref{eq:4.83} we can verify \eqref{eq:4.39}.
By the definitions of \( t_k, \theta_k, \tau_k, \sigma_k \) and \( \gamma_k \), we have
\begin{equation}\label{eq:4.99}
\begin{aligned}
    \tsum_{k=0}^{K-1} t_k = \tfrac{K(K+1+2k_0)}{2},
\end{aligned}
\end{equation}
\begin{equation}\label{eq:4.100}
\begin{aligned}
    \tsum_{k=0}^{K-1} \tfrac{t_k}{\tau_k} = \tfrac{K}{\mu_F}
\end{aligned}
\end{equation}
\begin{equation}\label{eq:4.101}
\begin{aligned}
    \tsum_{k=0}^{K} \tfrac{t_k }{\sigma_k} =\tsum_{k=0}^{K} \tfrac{2(k+k_0+1)}{(k+k_0)\mu_y} \leq \tfrac{2K}{\mu_y},
\end{aligned}
\end{equation}
To show the second part of the conclusion 
\begin{equation*}
\begin{aligned}
\hat{K} = \lceil 
\max \{ 
\sqrt{\tfrac{(2048) (L_g^{yx})^2 (\|\lambda^*\|_1 + 1)}{7\mu_y \mu_F}}, \quad
\tfrac{32 L_g^{xx} (\|\lambda^*\|_1 + 1)}{\mu_F} 
\} 
- k_0 - 1 .
\rceil
\end{aligned}
\end{equation*} 
where $\zeta$ is defined as \eqref{eq:4.84}
In fact , when $k\geq \hat{K}$, it can be verified that 
\begin{equation*}
\begin{aligned}
        t_k L_g^{xx} \| \lambda \|_1  
+ \tfrac{8 t_{k}}{\sigma_{k+1}} \|\lambda\|_1 \left( L_g^{yx} \right)^2 \leq 
\tfrac{7 t_k \tau_k}{64}+\tfrac{7 t_k \tau_k}{64}=\tfrac{7 t_k \tau_k}{32},
\end{aligned}
\end{equation*}
The sum in the definition of $\zeta$ with indexes $k\geq \hat{K}$ is less than or equal to zero. In addition, we have $\hat{K}\leq 0$ and thus $\zeta\leq 0$ when \eqref{eq:4.97} holds.
Using the above \eqref{eq:4.99},\eqref{eq:4.100},\eqref{eq:4.101} we obtain \eqref{eq:4.96}
\end{sharedproof}
\begin{theorem}
Suppose \(\mu_F=0\), \(\mu_y > 0\) , and the algorithm uses parameter sequences \(\{t_k, \theta_k, \tau_k, \gamma_k, \sigma_k\}\)
\begin{equation}\label{eq:4.87}
\begin{aligned}
t_k=1 , \theta_k = 1,\tau_k = \max\{ \tau',   \tau \sqrt{K} \}, \sigma_k = (k + k_0 + 1)   \mu_y, \gamma_k = 6M_g
\end{aligned}
\end{equation}
where $\tau' \ge \max\left\{
\tfrac{16}{7}(3L_F+10M_g),
\tfrac{16}{7}(L_g^{yx}+2L_g^{xx})(\|\lambda^{*}\|_{1}+1)
\right\}$, and $\tau\geq H_g^x (\|\lambda^{*}\|_{1}+1)$,  and \\
$k_0\ge max\{\tfrac{\sqrt{128/5}   L_g^{yy}}{\mu_y},16L_g^{yx}/\mu_y\}$. Then

\begin{align}
    \langle F(x), \bar{x}_K - x \rangle
    \le \tfrac{\tau' \| x - x_0 \|^{2}}{2K}
    + \tfrac{\tau \| x - x_0 \|^{2}}{2\sqrt{K}},
\end{align}
\begin{equation}\label{eq:4.88}
\begin{aligned}
    \| [g^{*}(\bar{x}_K)]_{+} \|_{\infty}
\le&\tfrac{\tau' \| x - x_0 \|^{2}}{2K}
+\tfrac{\tau \| x - x_0 \|^{2}}{2\sqrt{K}}
+\tfrac{ 3M_g(\| \lambda^{*} \|_{1} + 1)^{2}}{K}+\tfrac{(k_0 + 1)\mu_y D_y^2(\|\lambda^*\|_1 + 1)}{2K}\\
&+\tfrac{8(H_g^y)^2(\|\lambda^*\|_1 + 1)\ln({K + k_0 + 1}/{k_0})}{K \mu_y}
+\tfrac{4 H_g^x (\| \lambda^{*} \|_{1} + 1)}{\sqrt{K}}.
\end{aligned}
\end{equation}
\end{theorem}
\begin{sharedproof}
It is easy to verify that \eqref{eq:4.38} holds. In addition,
\begin{equation*}
\begin{aligned}
    \tfrac{8t_k(L_g^{yy})^2}{\sigma_{k+1}} 
    = \tfrac{8(L_g^{yy})^2}{(k + k_0 + 2)\mu_y}
    \leq \tfrac{5 t_k \sigma_k}{16},
\end{aligned}
\end{equation*}
where the inequality holds because \( k_0 \geq \tfrac{\sqrt{128/5}   L_g^{yy}}{\mu_y} \). This verifies \eqref{eq:4.41}.\\
By the definitions of \(t_k,\theta_k,\tau_k,\sigma_k\), and \(\gamma_k\), we have
\vspace{-0.5em}
\begin{gather}
\label{eq:4.90}
\tsum_{k=0}^{K-1} t_k = K,
\\
\label{eq:4.91}
\tsum_{k=0}^{K-1} \tfrac{t_k}{\tau_k}
\leq
\tfrac{K}{\tau\sqrt{K}}
\leq
\tfrac{K^{1/2}}{H_g^x(\norm{\lambda^*}_1+1)},
\\
\label{eq:4.92}
\tsum_{k=0}^{K} \tfrac{t_k}{\sigma_k}
=
\tsum_{k=0}^{K}
\tfrac{1}{(k+k_0+1)\mu_y}
\leq
\tfrac{\ln\!(\tfrac{K+k_0+1}{k_0})}{\mu_y}.
\end{gather}

Furthermore,
\begin{equation*}
\begin{aligned}
    \tfrac{L_F}{2}+L_F
+ \tfrac{6 t_k M_g^2}{\gamma_{k+1}}
+ \tfrac{24 t_{k+1} M_g^2}{\gamma_{k+2}}= \tfrac{3L_F}{2}
+ \tfrac{10 M_g}{2}
\leq \tfrac{7t_k \tau_k}{32},
\end{aligned}
\end{equation*}
where \(\tau_k \geq \tau' \geq \tfrac{16}{7}(3L_F + 10M_g)\). This verifies \eqref{eq:4.39}.
In addition, $\Psi(\|\lambda^*\|_1+1,K)\leq0$ in \eqref{eq:4.65} . We have for any $k\ge0$,
\begin{equation*}
\begin{aligned}
\tfrac{8t_k(L_g^{yx})^{2}(\|\lambda^{*}\|_{1}+1)}{\sigma_{k+1}}+t_kL_g^{xx}(\|\lambda^{*}\|_{1}+1)&=\tfrac{8t_k(L_g^{yx})^{2}(\|\lambda^{*}\|_{1}+1)}{(k+k_0+2)\mu_y}+t_kL_g^{xx}(\|\lambda^{*}\|_{1}+1)\\
&\le
\tfrac{L_g^{yx}(\|\lambda^{*}\|_{1}+1)}{2}+L_g^{xx}(\|\lambda^{*}\|_{1}+1)\le
\tfrac{7t_k\tau_k}{32} ,
\end{aligned}
\end{equation*}
where the first inequality because $k_0 \geq 16 L_g^{yx}/\mu_y$, and the second inequality because $t_k\tau_k\geq \tau'\geq \tfrac{16}{7}(L_g^{yx}+2L_g^{xx})(\|\lambda^{*}\|_{1}+1)$
\\
Applying \eqref{eq:4.90}, \eqref{eq:4.91}, \eqref{eq:4.92} to the right hand side of  \eqref{eq:4.65} we can obtain \eqref{eq:4.88}
\end{sharedproof}
\noindent\textbf{Remark.}
The choices of \(\tau\) and \(\tau'\) above require an upper bound
on \(\|\lambda^*\|_1\). Section~\ref{subsec:adaptive-pdconex} presents a new algorithm that removes this dependence without affecting the convergence rate when \(\mu_y>0\).
\begin{theorem}\label{thm:monotone_conv}
Suppose \(\mu_F = 0\) ( \(F\) is monotone ), \(\mu_y=0\) and the algorithm uses parameter sequences \(\{t_k, \theta_k, \tau_k, \gamma_k, \sigma_k\}\).  
\begin{align}\label{eq:4.67}
& t_k = 1, \quad \theta_k = 1, \gamma_k = 12 M_g,  \tau_k=  \max\{ \tau', \tau\sqrt{K} \}, \sigma_k = \max\{\sqrt{128/5}  L_g^{yy},  16L_g^{yx},  \tfrac{H_g^y}{D_y} \sqrt{K}\}
\end{align}
where $\tau' \ge \max\left\{ \tfrac{16}{7}(3L_F+5M_g), 
\tfrac{16}{7}(L_g^{yx}+2L_g^{xx})(\|\lambda^{*}\|_{1}+1)\right\}, 
\quad 
\tau \geq H_g^x (\|\lambda^{*}\|_{1}+1).$
Then the following holds:
\vspace{-0.5em}
\begin{gather}
\label{eq:4.68}
\inner{F(x)}{\bar{x}_K-x}
\le
\tfrac{\tau'}{2K}\norm{x-x_0}^{2}
+
\tfrac{\tau}{2\sqrt{K}}\norm{x-x_0}^{2},
\\
\label{eq:4.69}
\begin{aligned}
\norm{[g^{*}(\bar{x}_K)]_{+}}_{\infty}
\leq&
\tfrac{\tau'}{2K}\norm{x-x_0}^{2}
+
\tfrac{\tau}{2\sqrt{K}}\norm{x-x_0}^{2}
+
\tfrac{6M_g(\norm{\lambda^{*}}_{1}+1)^{2}}{K}
+
\tfrac{
\max\{\sqrt{32/5} L_g^{yy}, 8L_g^{yx}\}
D_y^{2}(\norm{\lambda^{*}}_{1}+1)
}{K}\\
&+
\tfrac{17H_g^{y}D_y(\norm{\lambda^{*}}_{1}+1)}{2\sqrt{K}}
+
\tfrac{4H_g^x(\norm{\lambda^{*}}_{1}+1)}{\sqrt{K}} .
\end{aligned}
\end{gather}

\end{theorem}
\begin{sharedproof}
It is easy to verify that \eqref{eq:4.38} holds.  
In addition,
\begin{equation*}
\begin{aligned}
\tfrac{8 t_k(L_g^{yy})^{2}}{\sigma_{k+1}}
\le \tfrac{8(L_g^{yy})^{2}}{\sqrt{128/5} L_g^{yy}}
\le \tfrac{5t_k \sigma_k}{16},
\end{aligned}
\end{equation*}
where the second inequality follows from \(\sigma_k \ge \sqrt{128/5} L_g^{yy}\).  
This verifies \eqref{eq:4.41}.  
Furthermore,

\begin{equation*}
\begin{aligned}
\tfrac{L_F}{2}+L_F
+ \tfrac{6 t_k M_g^2}{\gamma_{k+1}}
+ \tfrac{24 t_{k+1} M_g^2}{\gamma_{k+2}}
= \tfrac{3L_F}{2}
+ \tfrac{M_g}{2}
+ 2M_g
\leq \tfrac{7t_k \tau_k}{32},
\end{aligned}
\end{equation*}
since \(\tau_k \ge \tau' \ge \tfrac{16}{7}(3L_F+5M_g)\). This verifies \eqref{eq:4.39}.  
By the definitions of \(t_k,\tau_k,\sigma_k,\gamma_k\), we have
\vspace{-0.5em}
\begin{gather}
\label{eq:5.35}
\tsum_{k=0}^{K-1} t_k = K,
\\
\label{eq:4.73}
\tsum_{k=0}^{K-1} \tfrac{t_k}{\tau_k}
\le
\tfrac{K}{\tau\sqrt{K}}
\le
\tfrac{\sqrt{K}}{H_g^x(\norm{\lambda^*}_1+1)},
\\
\label{eq:4.74}
\tsum_{k=0}^{K} \tfrac{t_k}{\sigma_k}
\le
\tfrac{K}{(H_g^{y}/D_y)\sqrt{K}}
=
\tfrac{D_y\sqrt{K}}{H_g^{y}} .
\end{gather}
Furthermore, $\Psi(\|\lambda^*\|_1+1,K)\leq0$ in \eqref{eq:4.65} . We have for any $k\ge0$,
\begin{equation*}
\begin{aligned}
\tfrac{8t_k(L_g^{yx})^{2}(\|\lambda^{*}\|_{1}+1)}{\sigma_{k+1}}
+ t_k L_g^{xx}(\|\lambda^{*}\|_{1}+1)
&\le \tfrac{L_g^{yx}(\|\lambda^{*}\|_{1}+1)}{2}
+ L_g^{xx}(\|\lambda^{*}\|_{1}+1) \le \tfrac{7t_k \tau_k}{32}.
\end{aligned}
\end{equation*}

Finally, we obtain \eqref{eq:4.68} and \eqref{eq:4.69} by applying \eqref{eq:5.35}, \eqref{eq:4.73}, \eqref{eq:4.74} to the right-hand side of \eqref{eq:4.64}, and \eqref{eq:4.65}.
\end{sharedproof} 

\subsection{Adaptive PD-Lagrangian Extrapolation}
\label{subsec:adaptive-pdconex}

We finally present an adaptive variant of PD-ConEx for the smooth case with
$\mu_y>0$. The motivation is that the primal KKT operator is not uniformly
Lipschitz with respect to the multiplier variable. Define
\begin{equation}
P(x,y,\lambda):=F(x)+g_x'(x,y)^\top\lambda .
\label{eq:adap-P}
\end{equation}
For fixed $y$, the standing smoothness assumptions give
\[
\|P(x,y,\lambda)-P(x',y,\lambda')\|
\le
\big(L_F+L_g^{xx}\|\lambda'\|_2\big)\|x-x'\|
+M_g\|\lambda-\lambda'\|_2.
\]
Thus, the effective Lipschitz constant depends on the generated multiplier
sequence. Instead of imposing an a priori bound on $\lambda$, we adapt the
step sizes to $\|\lambda_k\|_2$. Among the smooth adaptive cases considered
here, the strongly concave case $\mu_y>0$ is the one for which we establish
convergence without requiring an a priori bound on $\lambda^*$.\\
Let
$\bar L:=\max\{L_F,M_g,1\}$, $c_1:=32$, $c_2:=16$,
$\beta:=2$, and $\rho_x:=1/8$. with $\lambda_0=0$,define $
a_k:=c_1\bar L+c_2L_g^{xx}\|\lambda_k\|_2,\qquad
a_0:=c_1\bar L,
$\\
Choose $
\tau_k:=a_k,\qquad
\gamma_k:=\beta a_k,\qquad
t_k:=\tfrac{a_0}{a_k},\qquad
\theta_0:=0,\qquad
\theta_k:=\tfrac{t_{k-1}}{t_k}\quad(k\ge1).
$\\
For the $y$-block, 
$
b_0\ge4L_g^{yy},\qquad
b_{k+1}:=b_k+\mu_yt_k,\qquad
\sigma_k:=\tfrac{b_k}{t_k}.
$\\
Hence
$
t_k\tau_k=a_0,\qquad
t_k\gamma_k=\beta a_0,\qquad
t_k\theta_k=t_{k-1},\qquad
t_k\sigma_k=b_k.
$
\begin{algorithm}[H]
\caption{Adaptive PD-Lagrangian Extrapolation}
\label{alg:adaptive-pdconex}
\small
\begin{algorithmic}[1]
\State Choose $x^0\in X$, $y^0\in Y$, $\lambda^0=0$; set
$x^{-1}=x^0$, $y^{-1}=y^0$, $\lambda^{-1}=\lambda^0$.
\For{$k=0,\ldots,K-1$}
\For{$i=1,\ldots,m$}
\State $u_i^k=(1+\theta_k)g'_{i,y}(x^k,y_i^k)
-\theta_kg'_{i,y}(x^{k-1},y_i^{k-1})$.
\State $y_i^{k+1}=\arg\min_{y_i\in Y_i}
\left\{-\langle u_i^k,y_i\rangle
+\tfrac{\sigma_k}{2}\|y_i-y_i^k\|^2\right\}$.
\EndFor
\State $s_k=g(x^k,y^{k+1})
+\theta_k[g(x^k,y^k)-g(x^{k-1},y^k)]$.
\State $\lambda^{k+1}=\arg\min_{\lambda\ge0}
\left\{-\langle s_k,\lambda\rangle
+\tfrac{\gamma_k}{2}\|\lambda-\lambda^k\|_2^2\right\}$.
\State $w_k=P(x^k,y^{k+1},\lambda^k)
+\theta_k[P(x^k,y^k,\lambda^k)
-P(x^{k-1},y^k,\lambda^{k-1})]$.
\State $x^{k+1}=\arg\min_{x\in X}
\left\{\langle w_k,x\rangle
+\tfrac{\tau_k}{2}\|x-x^k\|^2\right\}$.
\EndFor
\State $\Gamma_K:=\sum_{k=0}^{K-1}t_k$,\qquad
$\bar x_K:=\Gamma_K^{-1}\sum_{k=0}^{K-1}t_kx^{k+1}$.
\end{algorithmic}
\end{algorithm}
The adaptive parameter choice satisfies
\begin{equation}
t_{k-1}L_F\le\tfrac{a_0}{c_1},\qquad
t_{k-1}M_g\le\tfrac{a_0}{c_1},\qquad
t_{k-1}L_g^{xx}\|\lambda_{k-1}\|_2
\le\tfrac{a_0}{c_2}.
\label{eq:adap-step-control}
\end{equation}
These bounds allow the primal cross terms to be absorbed by the proximal terms, while the increasing sequence \(b_k\) controls the \(x\)--\(y\) coupling through the strong
concavity of $g_i(x,\cdot)$.
For the convergence analysis, we use the gap function $
Q(z_{k+1},z):=
\langle F(x),x_{k+1}-x\rangle
+\lambda^\top g(x_{k+1},y)
-(\lambda_{k+1})^\top g(x,y_{k+1}),
$
where \(z=(x,\lambda,y)\). The proof is given in
Appendix~\ref{app:adaptive-proof}.

\begin{theorem}[Adaptive PD-Lagrangian Extrapolation]
\label{thm:adaptive-pdconex}
Suppose the standing smoothness assumptions hold and $\mu_y>0$. Let
$\{x_k,\lambda_k,y_k\}$ be generated by
Algorithm~\ref{alg:adaptive-pdconex}. Then there exists a constant
$R_{\rm ad}<\infty$, independent of $K$, such that
$
\|\lambda_k\|_2\le R_{\rm ad},
 k=0,\ldots,K.
$\\
Consequently, there exists a constant $c_{\Gamma}>0$, independent of $K$,
such that $\Gamma_K\ge c_{\Gamma}K.$
Moreover, there exist constants $C_{\rm gap},C_{\rm feas}>0$, independent
of $K$, such that
\begin{align}
\sup_{x\in\widetilde X}
\langle F(x),\bar x_K-x\rangle
&\le\tfrac{C_{\rm gap}}{K},
\label{eq:adaptive-gap}\\
\|[\phi(\bar x_K)]_+\|_\infty
&\le\tfrac{C_{\rm feas}}{K}.
\label{eq:adaptive-feas}
\end{align}
Thus, Adaptive PD-Lagrangian Extrapolation achieves an $O(K^{-1})$ rate for both the VI gap
and the worst-case semi-infinite feasibility violation without requiring
an a priori bound on the generated multipliers.
\end{theorem}
See Appendix~\ref{app:adaptive-proof} for the proof.

\section{Saddle Point Problems with Coupling Constraints}
In this section, we consider the semi-infinite saddle-point problem
\begin{align}
\min_{u\in U}\;\max_{v\in V}\; f(u,v)
\quad\text{s.t.}\quad
g(u,v,y)\le 0 \quad \forall y\in Y.
\label{5.1}
\end{align}
The objective is convex in \(u\) and concave in \(v\), while
\(g(\cdot,\cdot,y)\) is jointly convex in \((u,v)\). Feasibility requires
the constraints to hold for every \(y\in Y\), where \(Y\) may be
infinite. Since constraints involve both \(u\) and \(v\), the
feasible choices of either variable depend on the other, as in
generalized Nash equilibrium problems
\cite{NabetaniTsengFukushima2011}. The semi-infinite constraints
also make direct projection onto the feasible set potentially
impractical. We therefore use a first-order primal--dual method with
extrapolated steps for \(u\) and \(v\) and proximal updates for the
multipliers and index variables.

Problem \eqref{5.1} extends the earlier VI--SIP formulation from a
single primal block to \(x=(u,v)\). The associated operator is
\begin{align}
F(x):=\big(\nabla_u f(u,v),-\nabla_v f(u,v)\big).
\label{eq:5.3}
\end{align}
Defining the feasible set by
$
\widetilde X
:=\{(u,v)\in U\times V:g(u,v,y)\le 0\ \forall y\in Y\},
$
we study the associated VI: find \(x^\star\in\widetilde X\) such that
\[
\langle F(x^\star),x-x^\star\rangle\ge 0
\qquad \forall x\in\widetilde X.
\]
The two-block formulation allows us to exploit strong convexity in
\(u\), when available, even if \(f\) is not strongly concave in \(v\).
This difference in blockwise curvature motivates the separate
convergence analysis below.
\begin{algorithm}[h]
\caption{PD-ConEx}
\label{alg:prob2}
\begin{algorithmic}[1]

\State \textbf{Inputs:} $u^0 \in \mathcal{U}$, $v^0 \in \mathcal{V}$, $\lambda^0 \in \mathbb{R}^m_+$, $K\in\mathbb{N}_+$; step sizes $\{ \sigma_k, \gamma_k, \theta_k, \tau_k^u, \tau_k^v \}$.
\State $u^{-2} = u^{-1} = u^0$, \; $v^{-2} = v^{-1} = v^0$, \; $y^{-1} = y^0$.

\For{$k=0,1,\dots,K-1$}
\For{$i=1,\dots,m$}
    \State $p_k^i \gets (1+\theta_k)\nabla_y g_i(u_k, v_k, y_k^i) 
           \;-\; \theta_k\nabla_y g_i(u_{k-1}, v_{k-1}, y_{k-1}^i)$
    \State $y_{k+1}^i \gets 
    \arg\min_{y^i \in \mathcal{Y}^i}
        -\langle p_k^i, y^i\rangle 
        + \tfrac{\sigma_k}{2} \|y^i - y_k^i\|^2$
\EndFor

\State $s_k \gets 
l_g(x_k; x_{k-1}, y_{k+1})
\;+\;
\theta_k [
    l_g(x_k; x_{k-1}, y_k)
    - l_g(x_{k-1}; x_{k-2}, y_k)
]$

\State $\lambda_{k+1} \gets 
\arg\min_{\lambda \ge 0}
    -\langle s_k, \lambda\rangle
    + \tfrac{\gamma_k}{2}\|\lambda - \lambda_k\|^2
$

\State {\small $v_{k+1} \gets 
\arg\min_{v\in\mathcal{V}}
-\big\langle
    (1+\theta_k)\nabla_v f(u_k, v_k)
    - \theta_k \nabla_v f(u_{k-1}, v_{k-1})
    + \lambda_{k+1}\nabla_v g(u_k, v_k, y_{k+1}),
    \; v
\big\rangle
+ \tfrac{\tau_k^v}{2}\|v - v_k\|^2
$}

\State {\small $u_{k+1} \gets 
\arg\min_{u\in\mathcal{U}}
\big\langle
    (1+\theta_k)\nabla_u f(u_k, v_k)
    - \theta_k \nabla_u f(u_{k-1}, v_{k-1})
    + \lambda_{k+1}\nabla_u g(u_k, v_k, y_{k+1}),
    \; u
\big\rangle
+ \tfrac{\tau_k^u}{2}\|u - u_k\|^2
$
}

\EndFor

\State \textbf{Output:} 
$\bar{x}_K 
= (\tsum_{k=0}^{K-1} t_k x_{k+1})/
        (\tsum_{k=0}^{K-1} t_k)$.

\end{algorithmic}
\end{algorithm}
\subsection{Convergence Guarantees for Algorithm~\ref{alg:prob2}}
Let \(f: U\times V\to\mathbb R\) be differentiable,
convex in \(u\), and concave in \(v\), and let
\(x=(u,v)\in\tilde X:= U\times V\).
We define the saddle-point operator \(F:\mathcal X\to\tilde X\) by
\begin{equation*}
F(x)
:=
\begin{pmatrix}
\nabla_u f(u,v)\\[2mm]
-\nabla_v f(u,v)
\end{pmatrix}
=
\begin{pmatrix}
F_u(x)\\[1mm]
F_v(x)
\end{pmatrix},
\end{equation*}
where $F_u(x):=\nabla_u f(u,v), 
F_v(x):=-\nabla_v f(u,v).$
Then, we know that the point \(x^*=(u^*,v^*)\) is a saddle point of \eqref{5.1} if and only if it solves the corresponding variational inequality
\begin{equation*}
\left\langle F(x^*),x-x^*\right\rangle \geq 0,
\qquad
\forall x\in\tilde X.
\end{equation*}
We define the following primal--dual gap function to measure the convergence for problem \eqref{5.1}
\begin{align*}
Q(z_{k+1}, z)
:=
\langle F(x), x_{k+1}-x \rangle
+
\langle \lambda, g(x_{k+1},y) \rangle
-
\langle \lambda_{k+1}, g(x,y_{k+1}) \rangle .
\end{align*}
Define the blockwise operator differences by
\[
\Delta F_k^u
:=
F_u(x_k)-F_u(x_{k-1}),
\qquad
\Delta F_k^v
:=
F_v(x_k)-F_v(x_{k-1}).
\]
Applying the optimality conditions of the \(u\)- and \(v\)-updates,
together with the three-point inequality and the block strong
monotonicity of \(F\), yields the corresponding blockwise
estimates. Combining these estimates with the convexity of
\(g(\cdot,\cdot,y)\), the convexity--concavity of \(f\), and the
multiplier and index-variable update bounds gives the following
weighted primal--dual gap estimate. 
\subsection{Step Size Policies}
\begin{proposition} \label{prop:summed-gap-alg2}
Suppose the parameters
${t_k, \theta_k,\tau^u_k,\tau^v_k ,\gamma_k,\sigma_k}$
of the algorithm satisfy the following conditions:
\begin{gather}
\label{eq:stepsize_conditon1}
t_k \tau^u_k
\leq
t_{k-1}(\tau^u_{k-1}+2\mu_F^u),
t_k \tau^v_k
\leq
t_{k-1}\tau^v_{k-1},
t_k\gamma_k
\leq
t_{k-1}\gamma_{k-1}, 
t_k\sigma_k
\leq
t_{k-1}(\sigma_{k-1}+\mu_y),
t_k\theta_k=t_{k-1},\\
\label{eq:stepsize_condition2}
\tfrac{6 t_k M_g^2}{\gamma_{k+1}}
+
t_k L_F^u
+
\tfrac{t_{k-1}L_F^u}{2}
+
\tfrac{24t_{k+1}M_g^2}{\gamma_{k+2}}
-
\tfrac{7t_k\tau_k^u}{32}
\leq 0 ,
\\
\label{eq:stepsize_condition2(2)}
\tfrac{6 t_k M_g^2}{\gamma_{k+1}}
+
t_k L_F^v
+
\tfrac{t_{k-1}L_F^v}{2}
+
\tfrac{24t_{k+1}M_g^2}{\gamma_{k+2}}
-
\tfrac{7t_k\tau_k^v}{32}
\leq 0 ,
\\
\label{eq:stepsize_cond4}
\tfrac{8t_k(L_g^{yy})^2}{\sigma_{k+1}}
+
\tfrac{3t_{k-1}\sigma_k}{16}
-
\tfrac{t_k\sigma_k}{2}
\leq 0 .
\end{gather}
for any \((u,v)\in X\), \(y\in Y\), and
\(\lambda\in\mathbb{R}_{+}^{m}\) satisfying
\(\|\lambda\|_1\leq \Lambda\), we have
\begin{equation}\label{eq:gap_prop3}
\begin{aligned}
&\Gamma_K \left[ f(\bar{u}_K,v)-f(u,\bar{v}_K) \right]
+\Gamma_K \left\langle g(\bar{u}_K,\bar{v}_K,y),\lambda \right\rangle
-\tsum_{k=0}^{K-1} t_k
\left\langle g(u,v,y_{k+1}),\lambda_{k+1}\right\rangle\\
&\leq
\tfrac{t_0\tau_0^u}{2}\|u-u_0\|^2
+
\tfrac{t_0\tau_0^v}{2}\|v-v_0\|^2
+
\tfrac{t_0\gamma_0}{2}\|\lambda-\lambda_0\|^2
+
\tfrac{t_0\sigma_0}{2}\|y-y_0\|_\lambda^2\\
&
+
\Psi_u(\Lambda,K)
+
\Psi_v(\Lambda,K)
+
\tsum_{k=0}^{K-1}
\left[
\tfrac{4t_k(H_g^u\Lambda)^2}{\tau_k^u}
+
\tfrac{4t_k(H_g^v\Lambda)^2}{\tau_k^v}
\right]
+
4\left(
\tsum_{k=1}^{K-1}\tfrac{t_k}{\sigma_{k+1}}
+
\tsum_{k=1}^{K-1}\tfrac{t_{k-1}}{\sigma_k}
\right)
(H_g^y)^2\Lambda .
\end{aligned}
\end{equation}
\text{where}
\begin{gather}
\label{eq:T-Lambda-K_u}
\Psi_u(\Lambda,K)
:=
\tsum_{k=1}^{K-1}
\left[
t_kL_g^{uu}\Lambda
+
\tfrac{8t_k\Lambda(L_g^{yu})^2}{\sigma_{k+1}}
+
\tfrac{t_kL_g^{uv}\Lambda}{2}
-
\tfrac{7t_k\tau_k^u}{32}
\right]
\|u_{k+1}-u_k\|^2,
\\[-0.3em]
\label{eq:T-Lambda-K_v}
\Psi_v(\Lambda,K)
:=
\tsum_{k=1}^{K-1}
\left[
t_kL_g^{vv}\Lambda
+
\tfrac{8t_k\Lambda(L_g^{yv})^2}{\sigma_{k+1}}
+
\tfrac{t_kL_g^{uv}\Lambda}{2}
-
\tfrac{7t_k\tau_k^v}{32}
\right]
\|v_{k+1}-v_k\|^2.
\end{gather}
\end{proposition}
 The detailed derivation of
Proposition~\ref{prop:summed-gap-alg2} is
provided in Appendix~\ref{app:proof-gap-alg2}\\
Let \(x=(u,v)\in\widetilde X\), and choose
\(\lambda=\lambda_0=0\) and \(y=y_0\),
the definition of \(Q\), followed by convexity--concavity of \(f\)
and Jensen's inequality, gives
\begin{align*}
f(\bar u_K,v)-f(u,\bar v_K)
\leq \tfrac{1}{\Gamma_K}
   \tsum_{k=0}^{K-1}t_k
   [f(u_{k+1},v)-f(u,v_{k+1})]\leq \tfrac{1}{\Gamma_K}
   \tsum_{k=0}^{K-1}t_k Q(z_{k+1},(u,v,0,y_0)).
\end{align*}
Hence, applying the bound from Proposition \ref{prop:summed-gap-alg2} and since we choose $\lambda=0$, we have $\Lambda=\|\lambda\|_1=0$. Hence all
positive $\lambda$-dependent terms in $\Psi_u(\Lambda,K)$ and
$\Psi_v(\Lambda,K)$ vanish.
Therefore, these terms can be dropped from the upper bound.
\begin{align*}
f(\bar u_K,v)-f(u,\bar v_K)
\leq
\tfrac{1}{\Gamma_K}
[
\tfrac{t_0\tau_0^u}{2}\|u-u_0\|^2
+
\tfrac{t_0\tau_0^v}{2}\|v-v_0\|^2
].
\end{align*}
For the infeasibility bound, we choose $ z=(u^*,v^*,\widehat{\lambda}_K,\widehat{y}_K), $ where $(u^*,v^*)$ is a feasible saddle-point solution and $\widehat{\lambda}_K$ is chosen and $\|\widehat{\lambda}_K\|_1\leq \|\lambda^*\|_1+1$,
Then
\begin{align*}
Q( z_{k+1}, \; z )
\ge
f( u_{k+1}, v^* )
-
f( u^*, v_{k+1} )
+
\left\langle 
\hat{\lambda}_K,\;
g( u_{k+1}, v_{k+1}, \hat{y}_K )
\right\rangle
-
\left\langle
\lambda_{k+1},\;
g( u^*, v^*, y_{k+1} )
\right\rangle.
\end{align*}
since$
f( \bar{u}_K, v^* )
-
f( u^*, \bar{v}_K )
\ge 0.$\\
Suppose 
$
\left\|
g^*[
( \bar{u}_K, \bar{v}_K )_+
]
\right\|_\infty
\;\le\;
\tfrac{1}{\Gamma_K}
\tsum_{k=0}^{K-1}
t_k  Q( z_{k+1}, \; z ).
$
In particular, we have:
\begin{align*}
\|
g^*( \bar{x}_K )_+
\|_\infty &\le\tfrac{1}{\tsum_{k=0}^{K-1} t_k}
[
 \tfrac{t_0\tau_0^u}{2}  \| u - u_0 \|^2
+
 \tfrac{t_0\tau_0^v}{2}  \| v - v_0 \|^2
+
 \tfrac{t_0\gamma_0}{2}  ( \| \lambda^* \|_1 + 1 )^2
+
 \tfrac{t_0\sigma_0}{2}  D_y^2  ( \| \lambda^* \|_1 + 1 )
]\\&\nonumber+\tfrac{1}{\tsum_{k=0}^{K-1} t_k}
[\Psi_u(\|\lambda^*\|_1+1,K)+\Psi_v(\|\lambda^*\|_1+1,K)]\\&\nonumber+ [\tfrac{4 ( (H_g^u)^2 ( \| \lambda^* \|_1 + 1 )^2  )}{\tsum_{k=0}^{K-1} t_k} \tsum_{k=0}^{K-1} \tfrac{t_k}{\tau^u_k}+\tfrac{4 ( (H_g^v)^2 ( \| \lambda^* \|_1 + 1 )^2  )}{\tsum_{k=0}^{K-1} t_k} \tsum_{k=0}^{K-1} \tfrac{t_k}{\tau^v_k}
+ \tfrac{8 (H_g^y)^2 ( \| \lambda^* \|_1 + 1 )}{\tsum_{k=0}^{K-1} t_k} \tsum_{k=0}^{K} \tfrac{t_k }{\sigma_k}
].
\end{align*}
\begin{theorem}
Suppose \(\mu^u_F > 0\), \(\mu^v_F = 0\),  \(\mu_y=0\) and the algorithm uses parameter sequences {$t_k,\theta_k, \tau^u_k,\tau^v_k,\gamma_k,\sigma_k$}
\begin{align*}
&t_k = k + k_0 +1, \theta_k = \tfrac{k +k_0}{k +k_0 +1}, \tau^u_k = \tfrac{t_k \mu^u_F}{2}, \tau^v_k = \tfrac{1}{t_k}\max\{\tfrac{96(K+k_0+1)L^{v}_F}{7}, H_g^v \sqrt{K}(K+k_0+1)\},
\end{align*}
\begin{align}
\gamma_k = \tfrac{1}{t_k}\max\{ \tfrac{260 M_g^2}{\mu^u_F }, \tfrac{65M_g^2K}{L^{v}_F}\},\sigma_k = \max \{\tfrac{\sqrt{128/5}(K+k_0) L_g^{yy},(H_g^y/D_y)K^{1/2}(K+k_0+1)}{t_k}
\},\label{eq:168}
\end{align}
Then the following holds;
\begin{align}
f(\bar u_K,v)-f(u,\bar v_K)
\leq
\tfrac{1}{2(K+1+2k_0)}
[
t_0\tau_0^u\|u-u_0\|^2
+
t_0\tau_0^v\|v-v_0\|^2
].
\label{eq:optimality_bound_6}
\end{align}
\begin{align}
    \| [g^*(\bar{x}_K)]_+ \|_\infty &\leq \tfrac{(k_0+1)^2\mu^u_F}{2K(K+2k_0+1)}\|u-u_0\|^2+\tfrac{96(K+k_0+1)L^{v}_F\|v-v_0\|^2}{7K(K+1+2k_0)}+\tfrac{H_g^v(K+k_0+1)}{\sqrt{K}(K+1+2k_0)}+\tfrac{260M_g^2(\|\lambda^*\|_1+1)^2}{K(K+2k_0+1)\mu^u_F}\nonumber\\ &+ \tfrac{65M_g^2K(k_0+1)(\|\lambda^*\|_1+1)^2}{L^{v}_FK(K+1+2k_0)}+ \tfrac{(K+k_0)\sqrt{128/5}D_y^2L_g^{yy}(\|\lambda^*\|_1+1)+H_g^yD_y\sqrt{K}(K+k_0+
    1)(\|\lambda^*\|_1+1)}{2K(K+2k_0+1)}\nonumber\\
    &+\tfrac{8(H_g^u)^2(\|\lambda^*\|_1+1)^2}{\mu^u_F(K+2k_0+1)}+\tfrac{16H^v_g(\|\lambda^*\|_1+1)^2(K+k_0+1)}{\sqrt{K}(K+2k_0+1)}+\tfrac{32H_g^yD_y(\|\lambda^*\|_1+1)(K+k_0+1)}{\sqrt{K}(K+1+2k_0 )}+\tfrac{2\Psi_u(\|\lambda^*\|_1+1,K)}{K(K+2k_0+1)} \nonumber\\ &+\tfrac{2\Psi_v(\|\lambda^*\|_1+1,K)}{K(K+2k_0+1)}\label{eq:infgap(2)},
   \end{align} 
where $\Psi_u(\|\lambda^*\|_1+1,K)$ is defined in \eqref{eq:T-Lambda-K_u} and it is $\Psi_u(\|\lambda^*\|_1+1,K)\leq 0$ when
\begin{align}
k_0 \geq \max\{
\tfrac{128\sqrt{5} (L_g^{yu})^2 (\|\lambda^*\|_1 + 1)}{7\sqrt{2}L_g^{yy} \mu^u_F} ,
\quad
\tfrac{16 L_g^{uu} (\|\lambda^*\|_1 + 1)}{\mu^u_F}
\} - 1\label{eq:eq_k0(2)}.
\end{align}
\end{theorem}  
\begin{sharedproof}
It is easy to verify that \eqref{eq:stepsize_conditon1} holds.
In addition, 
\begin{equation*}
     \tfrac{8t_k (L_g^{yy})^2}{\sigma_{k+1}}  \leq \tfrac{8 (k+k_0+1)^2 (L_g^{yy})^2}{\sqrt{128/5}(K+k_0)L_g^{yy}}\leq \tfrac{5(\sqrt{128/5})(K+k_0)L_g^{yy}}{16}\leq \tfrac{5 t_k \sigma_k}{16},
\end{equation*}
where $t_{k-1} \leq t_k$, $k\leq K-1$ and \eqref{eq:stepsize_cond4} holds.\\
Then by definitons of $t_k, \sigma_k, \tau_k$, and $\gamma_k$ we can write that;
\begin{equation*}
    \tsum_{k=0}^{K} \tfrac{t_k}{\sigma_k}= \tsum_{k=0}^{K} \tfrac{(k+k_0+1)^2}{(H_g^y/Dy)K^{1/2}(K+k_0+1)}\leq \tfrac{2D_yK^{1/2}(K+k_0+1)}{H_g^y} ,\label{eq:81}
\end{equation*}
where the inequality because $\tsum_{k=0}^{K} (k+k_0+1)^2 \leq (K+1)(k+k_0+1)^2 \leq 2K(K+k_0+1)^2$
Furthermore 
\begin{equation*}
    \tsum_{k=0}^{K-1} t_k = \tfrac{K (K + 1 + 2 k_0)}{2} ,
\end{equation*}
\begin{equation*}
    \tsum_{k=0}^{K-1} \tfrac{t_k}{\tau^u_k} = \tfrac{2K}{\mu^u_F} ,
\end{equation*}
\begin{equation*}
    \tsum_{k=0}^{K-1} \tfrac{t_k}{\tau_k^v}
= \tsum_{k=0}^{K-1} \tfrac{ (k+k_0+1)^2}{H_g^v (K+k_0+1)\sqrt{K}} \leq\tsum_{k=0}^{K-1} \tfrac{2K^{1/2}(K+k_0+1)}{H_g^v}.
\end{equation*}
In addition,
\begin{align*}
\tfrac{6 t_k M_g^2}{\gamma_{k+1}} + t_k L^u_F + \tfrac{t_{k-1} L^{u}_F}{2} + \tfrac{24 t_{k+1} M_g^2}{\gamma_{k+2}} &\leq \tfrac{3t_kL^{u}_F}{2}+\tfrac{6(k+k_0+1)^2M_g^2}{260Mg^2/\mu^u_F}+ \tfrac{24(k+k_0+2)^2M_g^2}{260M_g^2/\mu^u_F}\nonumber \\&\leq \tfrac{3t_k\tau^u_k}{32}+\tfrac{32.5t_k\tau^u_k}{260} = \tfrac{3t_k\tau^u_k}{32}+\tfrac{4t_k\tau^u_k}{32}=\tfrac{7t_k\tau^u_k}{32},
\end{align*}
where the inequality because$t_{k-1} \leq t_k$, $\tau^u_k\geq \tfrac{k_0\mu^u_F}{2} \geq \tfrac{128L^{u}_F}{7}\geq 16L^{u}_F$ where $k_0=\tfrac{128L^u_F}{7}$  and $k+k_0+2 \leq (1.05)(k+k_0+1)$ and \eqref{eq:stepsize_condition2} satisfied.
\begin{align*}
\tfrac{6 t_k M_g^2}{\gamma_{k+1}} + t_k L^v_F + \tfrac{t_{k-1} L^v_F}{2} + \tfrac{24 t_{k+1} M_g^2}{\gamma_{k+2}} &\leq \tfrac{3t_kL^{v}_F}{2}+\tfrac{6(k+k_0+1)^2M_g^2}{64M_g^2K/L^{v}_F}+ \tfrac{24(k+k_0+2)^2M_g^2}{64M_g^2K/L^{v}_F}\nonumber \\
&\leq \tfrac{3t_kL^{v}_F}{2}+\tfrac{32.5t_k^2L^{v}_F}{65K} \leq3t_kL^{v}_F(K+k_0+1)=\tfrac{7t_k\tau^v_k}{32},
\end{align*}
where the inequality because $t_{k-1} \leq t_k$ , $k+k_0+2 \leq (1.05)(k+k_0+1)$  and \eqref{eq:stepsize_condition2(2)} satisfied. \\
Define using \eqref{eq:stepsize_cond4} 
\begin{align*}
    \zeta' = \tsum_{k=0}^{K-1}
[
& \tfrac{8 t_k (L_g^{yu})^2 (\|\lambda^*\|_1 + 1)}{\sigma_{k+1}}
+ t_k L_g^{uu} (\|\lambda^*\|_1 + 1)
- \tfrac{7}{32} t_k \tau^u_k
] \nonumber \\
\leq 
\tsum_{k=0}^{\hat{K}-1}
[
& \tfrac{8 t_k (L_g^{yu})^2 (\|\lambda^*\|_1 + 1)}{\sigma_{k+1}}
+ t_k L_g^{uu} (\|\lambda^*\|_1 + 1)
- \tfrac{7}{32} t_k \tau^u_k
],
\end{align*}
where 
\begin{align} \label{eq:k_hat(2)} 
\hat{K}\geq [ \max \{
\tfrac{128\sqrt{5} (L_g^{yu})^2 (\|\lambda^*\|_1 + 1)}{7\sqrt{2}L_g^{yy} \mu^u_F} ,
\tfrac{16 L_g^{uu} (\|\lambda^*\|_1 + 1)}{\mu^u_F}
\} -k_0- 1] .
\end{align}
In fact, when $k\geq \hat{K}$ we have
\begin{align}
&t_k L_g^{uu} \| \lambda \|_1  
+ \tfrac{8 t_{k} \|\lambda\|_1}{\sigma_{k+1}} ( L_g^{yu} )^2 
-\tfrac{7 t_k \tau^u_k}{32}\nonumber\\
&\leq(k+k_0+1)L_g^{uu}(\|\lambda^*\|_1+1)+\tfrac{8(k+k_0+1)^2(\|\lambda^*\|_1+1)( L_g^{yu})^2}{\sqrt{128/5}(K+k_0)L_g^{yy}}-\tfrac{7 (k+k_0+1)^2\mu^u_F}{64}\nonumber\\
&\leq(k+k_0+1)L_g^{uu}(\|\lambda^*\|_1+1)+\tfrac{8(k+k_0+1)(\|\lambda^*\|_1+1)( L_g^{yu})^2}{\sqrt{128/5}L_g^{yy}}-\tfrac{7 (k+k_0+1)^2\mu^u_F}{64}\nonumber\\
&\leq 0.
\label{eq:k_hat(3)}
\end{align}
where the first inequality is from the definitions of $t_k,\theta_k,\tau^u_k,\sigma_k$ and the second inequality is because $k\leq K-1$, and the last inequality is because $k\geq\hat{K}$, when $\hat{K}$ is defined as \eqref{eq:k_hat(2)}. In addition,the condition \eqref{eq:eq_k0(2)} implies that $\widehat K\le 0$. 
Since $\widehat K\le 0$, inequality \eqref{eq:k_hat(3)} holds for every 
$k=0,\ldots,K-1$. Therefore, every coefficient appearing in 
$\Psi_u(\|\lambda^*\|_1+1,K)$ is nonpositive, and hence
$\Psi_u(\|\lambda^*\|_1+1,K)\le 0$.
In more applications, including robust and distributionally
robust models, the semi-infinite constraint does not depend on the maximization variable \(v\), or its dependence on \(v\) is separable or linear so that the corresponding mixed curvature constants vanish. In
particular, if $g(u,v,y)=g(u,y),$ then$ L_g^{vv}=L_g^{yv}=L_g^{uv}=0.$
Therefore, under this structural condition 
$\Psi_v(\|\lambda^*\|_1+1\*,K) = -\tsum_{k=1}^{K-1} \tfrac{7t_k\tau_k^v}{32} \|v_{k+1}-v_k\|^2 \le 0.$
\end{sharedproof}
\subsection{Special Case: Robust Optimization}\label{sec:robust_opt}
We now consider a special case of the general semi-infinite saddle-point problem \eqref{5.1}, where the constraint function depends only on the minimization variable. Specifically, we study
\begin{equation}\label{6.1}
    \min_{u \in \mathcal{U}} \ \max_{v \in \mathcal{V}} f(u,v) 
    \quad \text{s.t.} \quad g(u,y) \leq 0, \ \forall y \in \mathcal{Y}.
\end{equation}
Here, $f(\cdot,v)$ is convex in $u$ for each $v\in\mathcal V$, and $f(u,\cdot)$ is concave in $v$ for each $u\in\mathcal U$. The constraint function is independent of the maximization variable $v$, that is,
\[
g(u,v,y)\equiv g(u,y).
\]
This formulation is a direct specialization of the general semi-infinite saddle-point model. The variable $v$ still interacts with $u$ through the saddle-point objective $f(u,v)$, but it does not appear in the semi-infinite feasibility constraint. Hence, the constraint-induced coupling appears only in the $u$-block. Consequently, all terms involving the derivative of the constraint with respect to $v$ vanish, including$\nabla_v g, L_g^{vv},L_g^{yv}, L_g^{uv}$.
Therefore, the convergence analysis follows from the general case with simplified constants and without the $v$-block constraint remainder.
\\
\textbf{Convergence Guarantees:}
The saddle-point operator remains
\vspace{-0.5em}\[
F(u,v)=(\nabla_u f(u,v),-\nabla_v f(u,v)),\]
while the primal-dual gap bound reduces to the same form as in the general problem, but without any contribution from $\nabla_v g$. In particular, the only constraint-induced linearization error is the $u$-block remainder
\[
\Psi_u(\Lambda,K)
:=
\tsum_{k=1}^{K-1}
[
t_kL_g^{uu}\Lambda
+
\tfrac{8t_k\Lambda(L_g^{yu})^2}{\sigma_{k+1}}
-
\tfrac{7t_k\tau_k^u}{32}
]
\|u_{k+1}-u_k\|^2.
\]
Thus, under the same step-size conditions as in the general analysis, with the $v$-constraint constants set to zero, the averaged iterates $(\bar u_K,\bar v_K)$ satisfy the corresponding optimality and feasibility guarantees. In particular, taking$
\Lambda=\|\lambda^*\|_1+1,$
the feasibility violation is bounded by the initial distance terms, the dual and $y$-block terms, and the remaining $u$-block remainder $\Psi_u(\Lambda,K)$.
Moreover, if $k_0$ is chosen such that\[
k_0
\ge
\max\{
\tfrac{128\sqrt{5}(L_g^{yu})^2\Lambda}
{7\sqrt{2}L_g^{yy}\mu_F^u},
\tfrac{16L_g^{uu}\Lambda}{\mu_F^u}
\}-1,
\]
then $\Psi_u(\Lambda,K)\le 0$. Hence, this remainder does not increase the upper bound. Therefore, the convergence rate obtained for the general semi-infinite saddle-point problem remains valid in this special case, with improved constants due to the absence of constraint-induced coupling in the $v$-direction.
\subsection{Generalized Nash Equilibrium Problem (GNEP)}

A Nash equilibrium is a fundamental solution concept in game theory. It describes a stable outcome of a game in which no player has an incentive to change its strategy unilaterally. In a two-player game, each player chooses a strategy from its own feasible strategy set, and the payoff of each player may depend on the strategies chosen by both players. A strategy pair is called a Nash equilibrium if each player's strategy is optimal given the strategy of the other player. Therefore, at equilibrium, no player can improve its objective by changing only its own decision while the other player's decision remains fixed.

For a standard saddle-point problem without additional constraints,
\[
\min_{u\in U}\max_{v\in V} f(u,v),
\]
the minimization player chooses \(u\in U\), while the maximization player chooses \(v\in V\). A saddle point \((u^\star,v^\star)\) satisfies
\[
f(u^\star,v)\le f(u^\star,v^\star)\le f(u,v^\star),
\qquad \forall u\in U,\ \forall v\in V.
\]
This condition means that \(u^\star\) is optimal for the minimization player when \(v^\star\) is fixed, and \(v^\star\) is optimal for the maximization player when \(u^\star\) is fixed. Hence, the saddle-point condition can be viewed as the Nash equilibrium condition of a two-player zero-sum game.

A generalized Nash equilibrium extends this idea to games with coupling constraints. Unlike a standard Nash equilibrium, where each player's feasible set is independent, a GNEP allows the feasible strategy set of one player to depend on the other player's decision. Therefore, unilateral deviations are allowed only if they preserve the common feasibility requirements. This makes the equilibrium concept stronger than the standard Nash equilibrium, since a solution must satisfy both optimality and shared feasibility.

This distinction is important for the semi-infinite saddle-point problem
\[
\min_{u\in U}\max_{v\in V} f(u,v)
\quad \text{s.t.}\quad
g(u,v,y)\le 0,\quad \forall y\in Y.
\]
Here, the constraint \(g(u,v,y)\le 0\) depends on both saddle-point variables \(u\) and \(v\). Hence, the feasible region is no longer a simple product set \(U\times V\). Instead, the two players are coupled through the shared semi-infinite constraint. This coupling is the key feature that connects the constrained saddle-point problem to a generalized Nash equilibrium problem.

Therefore, the semi-infinite saddle-point problem can be interpreted as a generalized Nash equilibrium problem, because the solution must satisfy both saddle-point optimality and the shared semi-infinite feasibility condition. This viewpoint also motivates the primal-dual method used in this study: since projection onto the full constrained feasible set is generally difficult, the algorithm updates the saddle-point variables together with multiplier and index variables to control both the saddle-point gap and the feasibility violation.
\section{Application to Distributionally Robust Optimization}
Distributionally robust optimization (DRO) is widely used as an alternative to robust optimization (RO) in Section \ref{sec:robust_opt}. An RO model considers the worst-case possibility for each choice, which may have almost zero probability. A DRO model is only moderately conservative; it strikes a balance between standard {\em stochastic optimization} (SO) and the absolute worst case of RO by considering the underlying probability model and using historical data information. Consider the following three variations:
\[\text{SO}: \min_{x\in X} \Ebb_{\xi \sim P}[f(x,\xi)],\qquad \text{RO}:\min_{x\in X} \max_{\xi \in \Xi} f(x, \xi) ,\qquad \text{DRO}: \min_{x \in X} \max_{Q \in \mathcal{P}} \Ebb_{\xi \sim Q}[f(x,\xi)].\]

As one can see, a standard SO assumes complete knowledge of the probability distribution $P$. Here, $\xi$ is a random variable on a probability space $(\Xi, \Fcal, P)$ where $\Xi$ is the sample space with associated $\sigma$-algebra $\Fcal$. The RO objective overlooks the probability distribution completely and assumes that, for each choice of $x \in X$, the outcome $\xi \in \Xi$ must be the absolute worst-case even though the probability of that event, according to measure $P$, can be zero. The DRO model uses a set of probability distributions $\Pcal$, defined over the measurable space $(\Xi, \Fcal)$. The set $\Pcal$ is also referred to as an {\em ambiguity set} in the DRO literature. The ambiguity set $\Pcal$ is chosen such that $P \in \Pcal$. The maximum over the ambiguity set aims to strike a balance. While it does not assume complete knowledge of $P$ as is the case of SO, it still associates a distribution $Q$ - possibly in the neighborhood of $P$ - to consider the historical distributional information. This information is usually ignored by RO. The concept of DRO dates back to Scarf \cite{Scarf1958} and interested readers can refer the review paper \cite{RahimianMehrotra2019} for a comprehensive survey of the topic.

Even though DRO strikes a balance from a modeling point of view, the uncertainty in the distribution yields the optimization objective $\wt{f}(x) := \max_{Q \in \mathcal{P}} \Ebb_{\xi \sim Q}[f(x,\xi)]$ non-tractable for many important choices of ambiguity set $\Pcal$ due to the fact that variations are modeled over distributional space compared to finite-dimensional vector spaces of RO. This is true even when $f(\cdot, \xi)$ is convex for every $\xi \in \Xi$. A common technique used in the literature is to approximate $\wt{f}$ with a convex lower bounding function $\underline{f}(x)$ and optimize $\underline{f}$ as a proxy to keep the problem algorithmically tractable. Various techniques are used to create such lower bounds, and their effectiveness is usually showcased near the optimal solution based on the objective value achieved. 

Here, we are interested in solving a distributionally robust function-constrained optimization (DRFCO) as follows:
\begin{equation}\label{eq:DRFCO}
    \min_{x \in X} \max_{Q \in \Qcal}\Ebb_{\xi \sim Q}[f(x, \xi)] \quad \text{s.t.}\quad \wt{g}(x):= \max_{P \in \Pcal} \Ebb_{\xi \sim P}[g(x, \xi)] \le 0.
\end{equation}
In \eqref{eq:DRFCO}, we have introduced a distributionally robust function constraint to a standard DRO problem. This can be seen as a generalized DRO problem, similar to the generalized RO problem in \eqref{6.1}, which introduces a robust function constraint. The usual technique of obtaining a convex lower bound $\underline{g}(x)$ on distributionally robust functions $\wt{g}(x)$ does not work well when dealing with function constraints $\wt{g}(x) \le 0$. This is due to the fact that these lower bounds are shown to be effective empirically at only optimal solutions. To our best knowledge, there are no  theoretical or empirical guarantees of their accuracy near the level set  $\{x: \wt{g}(x) = 0\}$. We need such guarantees to ensure that tractable representation $\{x: \underline{g}(x) \le 0\}$ is a close outer-approximation of the set $\{x: \wt{g}(x)\le 0\}$. Hence, we consider a data-driven approximation of such functions proposed by \cite{GangammanavarBansal2022}, where they approximate the non-tractable ambiguity set (e.g., $\Pcal$ or $\Qcal$ here) based on available empirical data and the type of ambiguity set. We consider two widely used ambiguity sets here: (1) moment-based ambiguity and (2) Wasserstein distance-based ambiguity. For each type of ambiguity set, we construct a tractable data-driven approximation of \eqref{eq:DRFCO} that falls into our set up of \eqref{eq:vi}. Moreover, the approximation is provably asymptotically tight as the number of samples increases to infinity. The empirical effectiveness of such approximations is demonstrated in~\cite{GangammanavarBansal2022}. Throughout the remainder of the discussion, we only consider the ambiguity set $\Pcal$. The same arguments apply to the ambiguity set $\Qcal$ of the objective function.
\subsection{Moment-based ambiguity set}  
We first consider a moment-based ambiguity set for the distributionally robust function-constrained problem.  
Let \((\Xi,\mathcal F)\) be the measurable space of the uncertain parameter \(\xi\), and let \(\mathcal P(\Xi)\) denote the set of probability measures on \((\Xi,\mathcal F)\). A moment-based ambiguity set is defined by 
\begin{equation}
\mathcal P_{\rm mom} = \{ P\in \mathcal P(\Xi): \textstyle{\int}_{\Xi}\psi_j(\xi)\,dP(\xi)=b_j,\; j=1,\ldots,J \}.
\label{eq:population-moment-ambiguity}
\end{equation}where \(\psi_j:\Xi\to\mathbb R\) are prescribed moment functions and \(b_j\) are the corresponding moment values. The first condition ensures that \(P\) is a probability distribution, while the remaining conditions impose the moment restrictions. Thus, \(\mathcal P_{\rm mom}\) contains all probability distributions that match the specified moment information. For a function \(h:X\times \Xi\to\mathbb R\), the corresponding distributionally robust function is
\begin{equation} g_r(x) =\sup_{P\in \mathcal P_{\rm mom}} \mathbb E_{\xi\sim P}[g(x,\xi)]. \label{eq:robust-function-moment} 
\end{equation} In particular, the robust functional constraint \(g_r(x)\leq 0\) is equivalent to \begin{equation} \mathbb E_{\xi\sim P}[g(x,\xi)]\leq 0, \quad \forall P\in \mathcal P_{\rm mom}. 
\label{eq:semi-infinite-dro-constraint} \end{equation} 
Therefore, the ambiguity set \(\mathcal P_{\rm mom}\) plays the role of an uncertainty index set. This shows the direct connection between the DRO constraint and the semi-infinite constraint structure. Although \eqref{eq:semi-infinite-dro-constraint} gives the desired semi-infinite interpretation, it is not directly suitable for the numerical method. The reason is that the ambiguity set \(\mathcal P_{\rm mom}\) is defined over probability measures on \((\Xi,\mathcal F)\), and hence the maximization is over an infinite-dimensional distributional space. Moreover, in data-driven problems, the support and the moment values are generally unknown. Therefore, instead of working directly with \(\mathcal P_{\rm mom}\), we follow the data-driven approximation approach of~\cite{GangammanavarBansal2022}. In this approach, the unknown distribution is represented by probability weights assigned to the observed samples. As a result, given a finite sample \(\{\xi_i\}_{i=1}^N\), the ambiguity set is converted into a finite-dimensional set over \(p\in\mathbb R^N\). The moment information $b_j$ is estimated using the empirical moment $\widehat b_{j,N} = \tfrac{1}{N}\tsum_{i=1}^N \psi_j(\xi_i)$ for $j=1,\ldots,J$. Then, 
the sample-based moment ambiguity set is given by
\begin{equation}
\widehat{\mathcal P}_{N}^{\rm mom}
=\{ p\in \mathbb R_+^N: \tsum_{i=1}^N p_i=1,\;
\tsum_{i=1}^N p_i\psi_j(\xi_i)=\widehat b_{j,N},\;
j=1,\ldots,J \}.
\label{eq:sample-moment-ambiguity}
\end{equation} Here, \(p_i\) is the probability weight assigned to the sample point \(\xi_i\). The simplex constraint ensures that \(p\) is a probability vector, while the linear moment constraints ensure that the weighted empirical distribution satisfies the required empirical moment conditions. Using \eqref{eq:sample-moment-ambiguity}, the finite-sample robust function is
\begin{equation} \widehat g_N(x) =\max_{p\in \widehat{\mathcal P}_{N}^{\rm mom}} \tsum_{i=1}^N p_i g(x,\xi_i).
\label{eq:finite-sample-robust-function} \end{equation}
For each fixed \(x\), the maximization problem in \eqref{eq:finite-sample-robust-function} is a linear program over the probability weights \(p\). Since the ambiguity set is described by linear equality constraints and nonnegativity constraints, it is a convex finite-dimensional set. Therefore, the DRO approximation produces a tractable finite-dimensional robust function. The constraint $\widehat g_N(x)\leq 0 $ is equivalent to
\begin{equation} \tsum_{i=1}^N p_i g(x,\xi_i)\leq 0, \quad \forall p\in \widehat{\mathcal P}_{N}^{\rm mom}.
\label{eq:sample-semi-infinite-constraint} \end{equation} 
As one can see, \eqref{eq:sample-semi-infinite-constraint} is a semi-infinite function constraint. 
Hence, the moment-based DRO approximation fits naturally into the semi-infinite framework here. Moreover, applying the same idea to the DRO objective $f_r(x) := \max_{Q \in \mathcal{P}_{\rm mom}} \Ebb_{\xi \sim Q}[f(x,\xi)]$, we have the function $\widehat{f}_n(x) := \max_{q \in \widehat{\mathcal P}_{n}^{\rm mom}}\tsum_{i=1}^n q_if(x, \bar{\xi}_i)$ as a tractable approximation of $f_r(x)$. Here, $\{\bar{\xi}_i\}_{i=1}^n$ is possibly a different sample for the objective and $\widehat{\mathcal P}_{n}^{\rm mom}$ represents this sample based moment-ambiguity set similar to \eqref{eq:sample-moment-ambiguity}. In particular, we have 
\begin{equation}
\widehat{\mathcal P}_{n}^{\rm mom}
=\{ q\in \mathbb R_+^n: \tsum_{i=1}^n q_i=1,\;
\tsum_{i=1}^n q_i\psi_j(\bar{\xi}_i)=\widehat b_{j,n},\;
j=1,\ldots,I \}.
\label{eq:sample-moment-ambiguity_obj}
\end{equation}
Hence, we can approximate the DRO problem \eqref{eq:DRFCO} as follows:
\begin{equation}\label{eq:DRFCO_appx}
    \min_{x \in X} \max_{q \in \widehat{\mathcal P}_{n}^{\rm mom}}\wt{f}(x, q):= \tsum_{i=1}^n q_if(x, \bar{\xi}_i) \quad \text{s.t. } \wt{g}(x, p):= \tsum_{i=1}^Np_ig(x, \xi_i) \le 0, \quad \forall p \in \widehat{\mathcal P}_{N}^{\rm mom},
\end{equation}
Clearly, \eqref{eq:DRFCO_appx} is a special case of \eqref{5.1}. Hence, we can apply Algorithm \ref{alg:prob2} to solve the DRO problem.\\[2mm]
{\bf The quality of DRO approximation.} \hspace{1em}According to \cite[Proposition 2.1]{GangammanavarBansal2022}, the Hausdorff distance between $\widehat{\mathcal P}_{N}^{\rm mom}$ and ${\mathcal P}_{\rm mom}$ goes to zero as $N \to \infty$, almost surely. Hence, when the sample size $N$ (and $n$)is large, we get a high-quality approximation of ${\mathcal P}_{\rm mom}$ and consequently of $h_r(x)$ using \eqref{eq:finite-sample-robust-function}. The setting is tailor-made for large-scale data-driven optimization problems that are widely used across various applications relevant to DRO. It is also important to note that the tight asymptotic bounds are more useful than convex lower bounds for the function constraint, as mentioned earlier. 
\\[2mm]
{\bf Computational issues with projection onto $\widehat{\mathcal P}_{N}^{\rm mom}$.}\hspace{1em} Note that to apply Algorithm \ref{alg:prob2} for \eqref{eq:DRFCO_appx}, we need to project onto feasible sets $\widehat{\mathcal P}_{N}^{\rm mom}$ and $ \widehat{\mathcal P}_{n}^{\rm mom}$. Even though these sets are polytopes, they may not have suitable projection oracles due to presence of moment-based constraints, e.g., linear equality constraint $\tsum_{i=1}^N p_i\psi_j(\xi_i)=\widehat b_{j,N}$ in $\widehat{\mathcal P}_{N}^{\rm mom}$. 
This issue can be easily managed by incorporating Lagrange multipliers for moment constraints while the simplex constraint can be handled through projection onto the probability simplex. In particular, under assumption that Lagrange multipliers exist for \eqref{eq:DRFCO_appx},  we have its equivalent form
\begin{align*}
    &\min_{x\in X} \max_{q\in \widehat{\mathcal P}_{n}^{\rm mom}} \wt{f}(x,q)+ \max_{\lambda \ge 0}\lambda^\top [\max_{p \in \widehat{\mathcal P}_{N}^{\rm mom}} \wt{g}(x, p)],\\
    \equiv&\min_{x\in X}\max_{q\in \Delta_n}\min_{\mu_q \in \Rbb^I} \wt{f}(x,q)+\mu_q^\top (A_nq-\wh{b}_n) + \max_{\lambda \ge 0} \max_{p \in \widehat{\mathcal P}_{N}^{\rm mom}} \lambda^\top\wt{g}(x, p),\\
    \equiv& \min_{x\in X}\min_{\mu_q \in \Rbb^I}\max_{q\in \Delta_n} \wt{f}(x,q)+\mu_q^\top (A_nq-\wh{b}_n) + \max_{\lambda \ge 0} \max_{p \in \widehat{\mathcal P}_{N}^{\rm mom}} \lambda^\top \wt{g}(x,p).
\end{align*}
where the first equivalence follows by introducing the Lagrange multiplier $\mu_q$ for the linear moment constraint $A_nq=\widehat{b}_n$ and by using the nonnegativity of $\lambda$. In particular, since $\lambda\geq 0$, multiplying the inner maximization over $p$ by $\lambda$ does not change the maximizer. Moreover, the constraint $A_nq=\widehat{b}_n$ represents the moment-based approximation constraints
$\tsum_{i=1}^n q_i\psi_j(\bar{\xi}_i)=\widehat{b}_{j,n}$ for $j=1,\dots,I$, while $\Delta_n$ denotes the standard simplex in $\Rbb^n$. The intersection of the simplex constraint and these linear moment constraints yields $\widehat\Pcal_n^{\rm mom}$ in \eqref{eq:sample-moment-ambiguity_obj}. The second equivalence follows by the bilinear coupling between $\mu_q$ and $q$, which allows the order of $\min_{\mu_q}$ and $\max_q$ to be interchanged. Applying the same argument to the moment-based linear equality constraints defining $\widehat\Pcal_N^{\rm mom}$, we obtain the following equivalent reformulation:
\begin{equation}
    \min_{x\in X}\min_{\mu_q\in \Rbb^I}\min_{\mu_p \in \Rbb^J} \max_{q \in \Delta_n} \wt{f}(x,q) + \mu_q^\top (A_nq-\widehat b_n) + \max_{\lambda \ge 0}\max_{p \in\Delta_N } \lambda^\top[\wt{g}(x,p)+ \mu_p^\top(A_Np-\widehat b_N)],
    \label{eq:DRFCO_appx_full_min_max_form}
\end{equation}
where we used an additional switch between the order of $\min_{\mu_p}$ and $\max_\lambda$ due to bilinear coupling between $\mu_p$ and $\lambda$. Now, we treat minimization variables $(x, \mu_p, \mu_q) \in X \times \Rbb^I \times \Rbb^J$ as the joint primal variable, $q \in \Delta^n$ as the dual for the objective, $p \in \Delta^N$ as the semi-infinite variable for the constraints and finally $\lambda \ge 0$ as the Lagrange multiplier. This formulation satisfies all convexity and Lipschitz assumptions despite unbounded primal variables $\mu_p$ and $\mu_q$ and projection onto their feasible domain is easy. Now, we can easily apply Algorithm \ref{alg:prob2} without projection related issues.

This form is used in the numerical DRO experiment, where we verify the linear equality constraints are indeed satisfied by the output of our algorithm for both empirical ambiguity sets in the objective and the constraint.

\subsection{Wasserstein Distance-Based Ambiguity Set} The Wasserstein distance-based ambiguity set is another widely used ambiguity set in DRO. Unlike moment-based ambiguity sets, which describe distributional uncertainty through prescribed moment information, Wasserstein ambiguity sets describe distributional uncertainty through a distance from a reference distribution. 
Let \((\Xi,\mathcal F)\) be the measurable space of the uncertain parameter \(\xi\), and let \(\mathcal P(\Xi)\) denote the set of probability distributions on \(\Xi\). For two probability distributions \(P_1,P_2\in\mathcal P(\Xi)\), the Wasserstein distance is defined by 
\begin{equation}
W(P_1,P_2) = \inf_{\Pi\in\mathcal J(P_1,P_2)}
\scalebox{0.7}{$\displaystyle\int$}_{\Xi\times\Xi}
c(\xi,\xi')\,d\Pi(\xi,\xi'),
\label{eq:wasserstein-distance}
\end{equation}
where \(c(\xi,\xi')\) is the transportation cost between \(\xi\) and \(\xi'\), and \(\mathcal J(P_1,P_2)\) is the set of all joint distributions with marginals \(P_1\) and \(P_2\). Each \(\Pi\in\mathcal J(P_1,P_2)\) is called a transport plan. It describes how probability mass from \(P_1\) is transported to form \(P_2\). Hence, the Wasserstein distance is the minimum total transportation cost over all admissible transport plans. Given a reference distribution \(P_0\) and a radius \(\rho>0\), the Wasserstein ambiguity set is defined as
\begin{equation} \mathcal P_W = \left\{ P\in\mathcal P(\Xi): W(P,P_0)\leq \rho \right\}. \label{eq:population-wasserstein-ambiguity} \end{equation} 
The radius \(\rho\) controls the size of the ambiguity set. For a measurable function \(h:X\times\Xi\to\mathbb R\), the corresponding Wasserstein distributionally robust function is\begin{equation}
g_W(x) = \sup_{P\in\mathcal P_W} \mathbb E_{\xi\sim P}[g(x,\xi)].
\label{eq:wasserstein-robust-function}
\end{equation}
Therefore, the robust functional constraint \(g_W(x)\leq 0\) is equivalent to
\begin{equation}
\mathbb E_{\xi\sim P}[g(x,\xi)]\leq 0, \quad \forall P\in\mathcal P_W.
\label{eq:wasserstein-semi-infinite-population}
\end{equation}
Thus, the probability distribution \(P\) plays the role of an uncertainty index, and the Wasserstein DRO constraint has a semi-infinite interpretation. As earlier, this form is not implementable due to the infinite-dimensional space in which $\mathcal{P}_W$ resides. 
The sample-based approximation of the Wasserstein DRO constraints also admit a semi-infinite structure compatible with our framework. Let \(\Xi_N=\{\xi_1,\ldots,\xi_N\}\) be the observed sample set, and let \(\widehat p_i\) denote the empirical probability assigned to sample point \(\xi_i\). For example, if all observations are distinct and equally weighted, then \(\widehat p_i=1/N\). The empirical reference distribution is represented by the probability vector $\widehat p=(\widehat p_1,\ldots,\widehat p_N).$ The adversarial distribution supported on the same sample set is represented by another probability vector \(p=(p_1,\ldots,p_N)\). In this discrete finite-sample setting, the transport plan is represented by a nonnegative matrix \(\eta=(\eta_{ij})\), where \(\eta_{ij}\) denotes the amount of probability mass transported from \(\xi_j\) to \(\xi_i\). Hence, the integral transportation cost in \eqref{eq:wasserstein-distance} is replaced by the finite summation \[ \tsum_{i=1}^N\tsum_{j=1}^N c_{ij}\eta_{ij}, \] where \(c_{ij}=c(\xi_i,\xi_j)\). The finite-dimensional Wasserstein ambiguity set is then given by 
\begin{align} \label{eq:finite-wasserstein-ambiguity}
\widehat{\mathcal P}_{N}^{\rm mom} = \{ p\in\mathbb R_+^N:\; \tsum_{i=1}^N p_i=1,\; \exists\,\eta\in\mathbb R_+^{N\times N} \ \text{s.t.}\ \tsum_{j=1}^N \eta_{ij}=p_i,\; \tsum_{i=1}^N \eta_{ij}=\widehat p_j,\; \tsum_{i=1}^N\tsum_{j=1}^N c_{ij}\eta_{ij}\leq \rho \}.  \end{align}   
The constraints \(\tsum_{j=1}^N\eta_{ij}=p_i\) ensure that the total mass arriving at \(\xi_i\) equals the adversarial probability \(p_i\), while the constraints \(\tsum_{i=1}^N\eta_{ij}=\widehat p_j\) ensure that the total mass leaving \(\xi_j\) equals its empirical probability. The transportation-budget constraint limits the total transport cost. Using \eqref{eq:finite-wasserstein-ambiguity}, the finite-sample Wasserstein robust function becomes 
\begin{equation} 
\widehat h_N^W(x) = \max_{p\in\widehat{\mathcal P}_{N}^{\rm mom}} \tsum_{i=1}^N p_i g(x,\xi_i). \label{eq:finite-wasserstein-robust-function} \end{equation} 
For each fixed \(x\), this is a finite-dimensional linear optimization problem over the adversarial probability vector \(p\), with the transport variables \(\eta\) used to enforce the Wasserstein-distance constraint~\cite{GangammanavarBansal2022,EsfahaniKuhn2018}. The finite-sample robust constraint \(\widehat h_N^W(x)\leq 0\) can equivalently be written as 
\begin{equation} \tsum_{i=1}^N p_i g(x,\xi_i)\leq 0, \quad \forall p\in\widehat{\mathcal P}_{N}^{\rm mom}. \label{eq:finite-wasserstein-semi-infinite} \end{equation} 
The above formulation is a tractable semi-infinite function constraints that fits in our framework. Moreover, from \cite[Proposition 2.3]{GangammanavarBansal2022} shows the Hausdorff distance between $\mathcal{P}_W$ and $\widehat{\mathcal P}_{N}^{\rm mom}$ goes to zero as $N \to \infty$. Then, following the steps from moment-based ambiguity discussion, we can dualize linear constraints involving variable $\eta$ in \eqref{eq:finite-wasserstein-ambiguity} to get a projection friendly simplex set for variable $p$. The reformulation here is even more consequential since the number of linear constraints increase very fast for Wasserstein ambiguity set.

\section{Numerical Experiments}
\subsection{Experiment 1: Robust Quadratic Semi-Infinite Optimization}
In this section, we evaluate the practical performance of the PD-ConEx Method, on robust quadratic optimization instances and compare it with the \emph{Subgradient Saddle-Point} (SGSP) method. The goal of the experiments is to assess the behavior of the PD-ConEx algorithm in terms of both objective and constraint satisfaction on semi-infinite robust optimization problems.

We consider the robust optimization problem
\begin{equation}\label{eq:num_prob}
\min_{u \in  U} \ \sup_{z \in  Z} \ g_0(u,z)
\quad \text{s.t.} \quad
g_i(u,y) \le 0, \ \forall y \in Y,\ i\in[m],
\end{equation}
where the objective and constraint functions are of robust quadratic type. In particular, the uncertain functions are constructed in the form
\begin{equation}\label{eq:robust_quad_constraint}
g_i(u,y) =\|
(
P_{i0} + \tsum_{k=1}^{K_z} P_{ik} y_k
)u
\|^2
+
b_i^\top u + c_i,
\qquad i=1,\dots,m,
\end{equation}
and similarly for the objective function \(g_0(u,z)\). The uncertainty sets are Euclidean balls,
\[
\ Y=\{y\in\mathbb R^{K_z}:\|y\|_2\le D_y\},\qquad
\ Z=\{z\in\mathbb R^{K_z}:\|z\|_2\le D_z\},
\]
and the primal variable is constrained to
\[
\ U=\{u\in\mathbb R^{n_u}:\|u\|_2\le D_u\}.
\]
In all experiments, we use \(D_u=D_y=D_z=1\).
This class of problems is consistent with the robust quadratic benchmark used in prior work on first-order methods for robust optimization, and is particularly suitable for testing algorithms designed for semi-infinite constraints.
We compare the following methods.
\begin{itemize}
    \item \textbf{PD-ConEx.}
This is the proposed primal-dual first-order method developed in this work. It uses structured updates for the primal variable, the dual multipliers, and the auxiliary uncertainty variables, together with extrapolation and the step-size policy prescribed by the Theorems.
 \item \textbf{SGSP}
As a baseline, we implement the Subgradient Saddle-Point method applied to the robust formulation. For SGSP, we use the theoretically motivated horizon-based step-size rule
\[
\alpha = \tfrac{\tilde\tau}{\sqrt{N}},\qquad
\theta = \tfrac{\tilde\theta}{\sqrt{N}},
\]
with two choices of parameters:
$
(\tilde\tau,\tilde\theta)=(5,5),  $and $(\tilde\tau,\tilde\theta)=(1,1).
$
These two settings are included to illustrate the sensitivity of SGSP to the practical choice of scaling constants.
\end{itemize}
\textbf{Experimental setup:}
We evaluate the performance of the PD-ConEx method against a standard baseline on the robust quadratic semi-infinite optimization problems described above. The algorithms are run for a fixed iteration budget \(N\), and performance is measured in terms of the level function and feasibility violation. The same problem instances and parameter settings are used across all methods to ensure a fair comparison.
Two representative problem instances are generated:
\begin{itemize}
    \item \textbf{Set 1:} \((n_u,n_v,m,K_z,L)=(10,10,3,10,10)\),
    \item \textbf{Set 2:} \((n_u,n_v,m,K_z,L)=(600,25,3,15,15)\).
\end{itemize}
Here, \(n_u\) is the primal dimension, \(n_v\) is the dimension of the auxiliary \(v\)-block in the implementation, \(m\) is the number of robust constraints, \(K_z\) is the uncertainty dimension, and \(L\) controls the size of the quadratic data matrices.

All problem data are randomly generated and normalized and for each instance, all methods are initialized from the same infeasible starting point \(u^0\), generated so that at least one robust constraint is violated initially. This allows us to directly evaluate the ability of each method to reduce infeasibility while improving the objective value.
The maximum iteration budget is set to $N = 10^7$
and performance is recorded throughout the run.
To compare the methods, we consider the following two quantities.

First, we evaluate the \emph{feasibility violation}
\[
\max\{0,g(u)\},
\]
where \(g(u)\) denotes the worst robust constraint value over all constraints and uncertainty realizations. This quantity measures how far the current iterate is from satisfying the semi-infinite constraints.

\begin{figure}[htbp]
    \centering
    \begin{subfigure}[b]{0.45\textwidth}
        \centering
        \includegraphics[width=\textwidth]{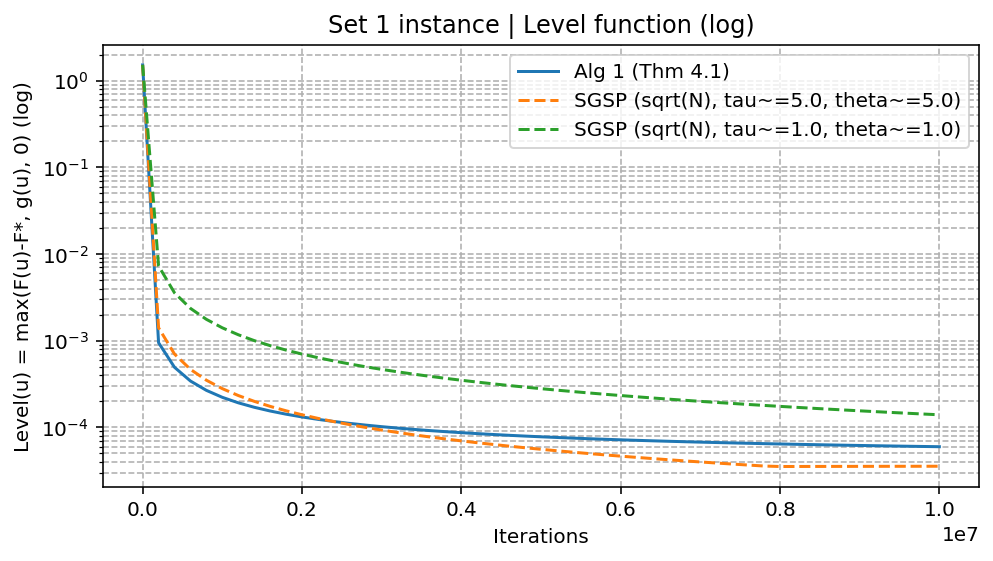}
        \caption{Level function for set 1}
        \label{set1_level}
    \end{subfigure}
    \hfill
    \begin{subfigure}[b]{0.45\textwidth}
        \centering
        \includegraphics[width=\textwidth]{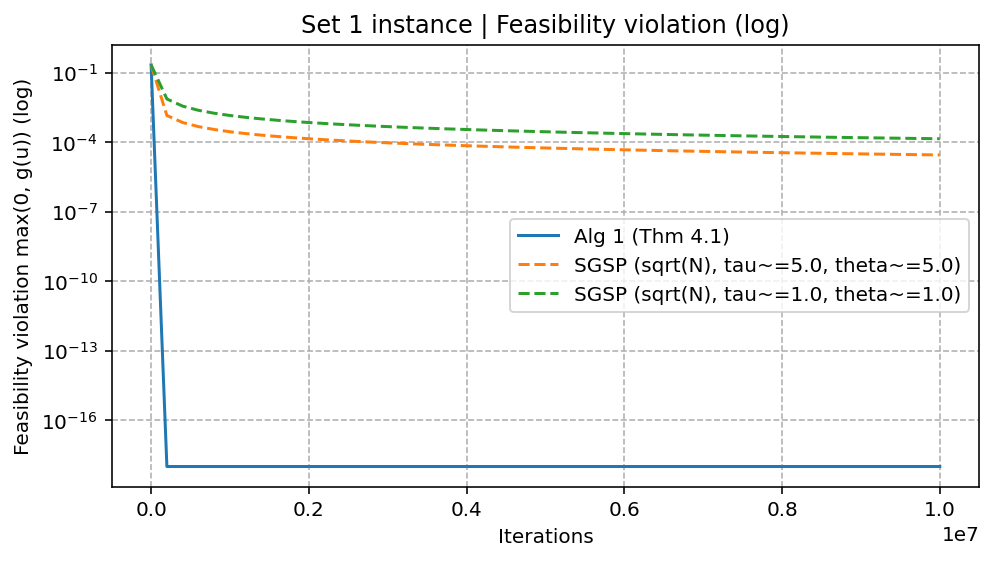}
        \caption{Feasibility violation for set 1}
        \label{set1_feas}
    \end{subfigure}
    \caption{Performance of the algorithms for set 1}
    \label{fig:combined}
\end{figure}

Second, we measure the \emph{level function}
\[
\mathrm{Level}(u)=\max\{f(u)-f^\star,\ g(u),\ 0\},
\]
where \(f(u)\) is the worst-case robust objective value and \(f^\star\) is estimated from the best feasible objective value observed across the runs. This criterion jointly captures progress in objective minimization and feasibility restoration. Smaller values indicate better overall performance. Here, all the curves are plotted on a logarithmic scale.\\
\textbf{Results:}
For the implementation, the step-size parameters used in the experiment were selected according to the theoretical policy in \ref{thm:monotone_conv}. When considering the \textbf{Results of set 1}, Figure~\ref{set1_level} reports the level function, while Figure~\ref{set1_feas} shows the feasibility violation.

The results indicate a clear advantage of the PD-ConEx algorithm over both SGSP variants. In terms of the level function, PD-ConEx decreases rapidly in the early stage and remains below the two SGSP curves over most of the run. While the SGSP method with \((\tilde\tau,\tilde\theta)=(5,5)\) eventually approaches a similar level of objective accuracy, the PD-ConEx method achieves lower values much earlier and with more stable progress.

PD ConEx algorithm drives  the constraint violation almost zero and does this very quickly, using only a small number of iterations. This behavior suggests that the PD-ConEx method is particularly effective at handling the semi-infinite constraint structure and restoring feasibility quickly from an infeasible initialization.

\begin{figure}[htbp]
    \centering
    \begin{subfigure}[b]{0.45\textwidth}
        \centering
        \includegraphics[width=\textwidth]{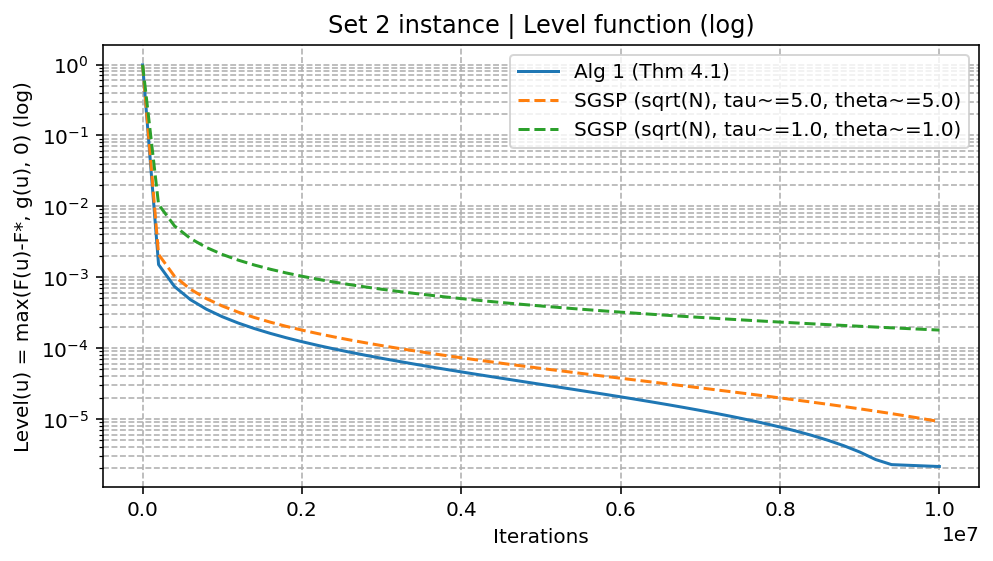}
        \caption{Level function for set 2}
        \label{set2_level}
    \end{subfigure}
    \hfill
    \begin{subfigure}[b]{0.45\textwidth}
        \centering
        \includegraphics[width=\textwidth]{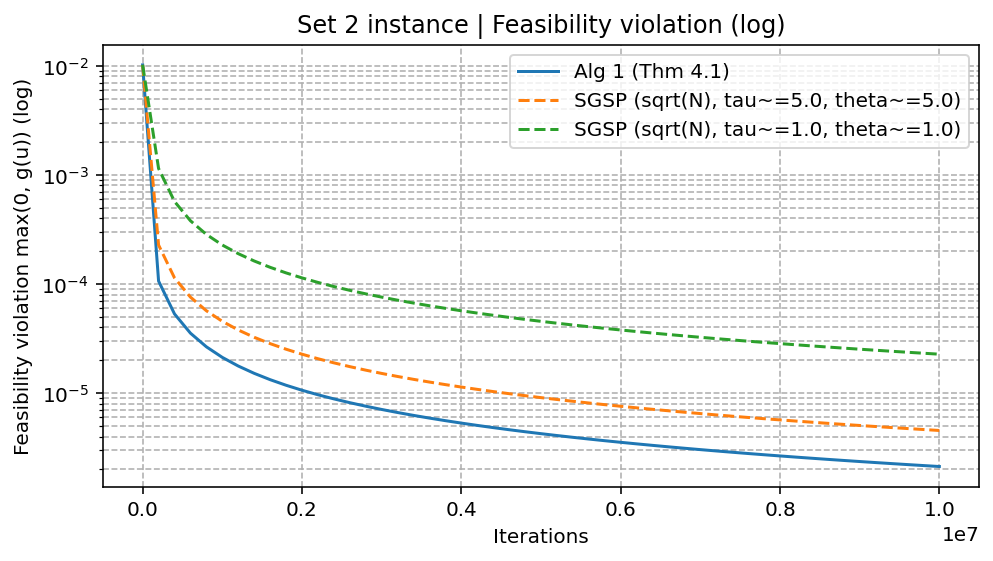}
        \caption{Feasibility violation for set 2}
        \label{set2_feas}
    \end{subfigure}
    \caption{Performance of the algorithms for set 2}
    \label{fig:combined}
\end{figure}

Overall, on Set 1, PD-ConEx algorithm exhibits both faster practical convergence and significantly stronger feasibility behavior than SGSP.
\\
Then we consider the \textbf{Results of set 2,}
Figures~\ref{set2_level} and \ref{set2_feas} report the corresponding results, which are substantially higher-dimensional and therefore more challenging.

The same general pattern is observed. PD-ConEx consistently outperforms both SGSP variants in the level function and achieves the smallest final value among all the methods tested. The improvement is more visible in the later stage of the run, where the proposed method continues to decrease steadily while the SGSP curves flatten earlier.

A similar advantage is seen in the feasibility violation. The PD-ConEx algorithm remains uniformly below both SGSP variants and reaches the lowest final constraint violation. Although the PD-ConEx method still shows a clear and consistent benefit over the baseline. This indicates that the practical advantage of PD-ConEx is preserved even as the dimension and complexity of the problem increase.
\subsection{Experiment 2: DRO Problem on the Mushroom Dataset}

In this experiment, we apply the moment-based distributionally robust optimization (DRO) formulation introduced above to the Mushroom dataset \cite{uci_mushroom}. The Mushroom dataset represents a binary classification problem in which each mushroom is classified as either edible or poisonous. We consider the projection-friendly reformulation in \eqref{eq:DRFCO_appx_full_min_max_form}, obtained from the moment-based DRO formulation in \eqref{eq:DRFCO_appx}, and apply the proposed primal-dual method to solve the resulting function-constrained DRO problem.

In this setting, the objective minimizes the distributionally robust logistic loss over the full dataset, while the constraint controls the distributionally robust logistic loss over the poisonous subset. This distinction is important because incorrectly classifying a poisonous mushroom as edible is considerably more consequential than classifying an edible mushroom as poisonous. Therefore, placing the poisonous-subset loss in the constraint provides explicit control over this more critical classification error. Moreover, using a DRO model rather than a fully robust optimization (RO) model incorporates distributional information through the ambiguity set and avoids an excessively conservative worst-case formulation that may unnecessarily degrade the overall classification performance. Let \(D\) denote the full dataset and let \(D_p\) denote the poisonous subset. For a given decision vector \(u\), let \(\ell_i(u)\) denote the logistic loss of sample \(i\). The DRO full-data objective is defined as
\begin{equation}
f(u)=\max_{q\in Q_D}\wt{f}(u,q),
\qquad
\wt{f}(u,q):=\tsum_{i\in D}q_i\ell_i(u),
\label{eq:exp2-fu}
\end{equation}
where \(q\) is a probability distribution over the full dataset. We use the moment-based ambiguity set
\begin{equation}
Q_D=\left\{
q\in\mathbb{R}_+^N:
\tsum_{i=1}^{N}q_i=1,\;
A_Dq=C_D
\right\}.
\label{eq:exp2-QD}
\end{equation}
Similarly, for the poisonous subset, we define the DRO function constraint
\begin{equation}
g(u):=\max_{p\in P_{D_p}}\wt{g}(u,p)\leq \bar g,
\qquad
\wt{g}(u,p):=\tsum_{i\in D_p}p_i\ell_i(u),
\label{eq:exp2-gu}
\end{equation}
where
\begin{equation}
P_{D_p}=\left\{
p\in\mathbb{R}_+^n:
\tsum_{i=1}^{n}p_i=1,\;
A_{D_p}p=C_{D_p}
\right\}.
\label{eq:exp2-P}
\end{equation}
Here, \(q_i\) denotes the probability weight assigned to the \(i\)-th sample in the full dataset, and \(p_i\) denotes the probability weight assigned to the \(i\)-th sample in the poisonous subset. The simplex constraints ensure that \(q\) and \(p\) are probability distributions. The matrices \(A_D\in\mathbb{R}^{d\times N}\) and \(A_{D_p}\in\mathbb{R}^{d\times n}\) contain the selected empirical moment features of the samples as columns.In particular,
$
A_D=
\begin{bmatrix}
\psi(\omega_1) & \psi(\omega_2) & \cdots & \psi(\omega_N)
\end{bmatrix}
\in\mathbb R^{d\times N},
$
and
$
A_{D_p}=
\begin{bmatrix}
\psi(\omega_1^p) & \psi(\omega_2^p) & \cdots & \psi(\omega_n^p)
\end{bmatrix}
\in\mathbb R^{d\times n},
$
where \(\psi(\omega_i)\in\mathbb R^d\) and \(\psi(\omega_i^p)\in\mathbb R^d\) denote the prescribed moment-feature vectors of the corresponding samples. The vectors \(C_D\in\mathbb R^d\) and \(C_{D_p}\in\mathbb R^d\) are the corresponding empirical moment vectors, defined by
\[
C_D=\tfrac{1}{N}\tsum_{i=1}^{N}\psi(\omega_i),
\qquad
C_{D_p}=\tfrac{1}{n}\tsum_{i=1}^{n}\psi(\omega_i^p).
\]
Hence, the constraints \(A_Dq=C_D\) and \(A_{D_p}p=C_{D_p}\) require the weighted moments induced by \(q\) and \(p\), respectively, to match the corresponding empirical moments of the full dataset and the poisonous subset.
Thus, the constraints \(A_Dq=C_D\) and \(A_{D_p}p=C_{D_p}\) are vector-valued moment constraints. They require the adversarial distributions \(q\) and \(p\) to preserve the selected empirical moments of the full dataset and the poisonous subset, respectively. 
The value of \(\bar g\) is computed as follows.
First, we solve the full-data DRO problem
\begin{equation}
    \min_{u\in U}\max_{q\in Q_D} f(u),
    \qquad
    \label{eq:exp2-gbar-stage1}
\end{equation}
and this produces \((\bar u,\bar q)\). 
Next, we evaluate the poisonous-subset DRO loss at the solution
\(\bar u\) by solving the linear program 
$ g^\star(\bar u)  =
    \max_{p\in P_{D_p}} \wt g(\bar u,p).$

\begin{figure}[htbp]
\centering

\begin{subfigure}{0.48\textwidth}
\centering
\includegraphics[width=\textwidth]{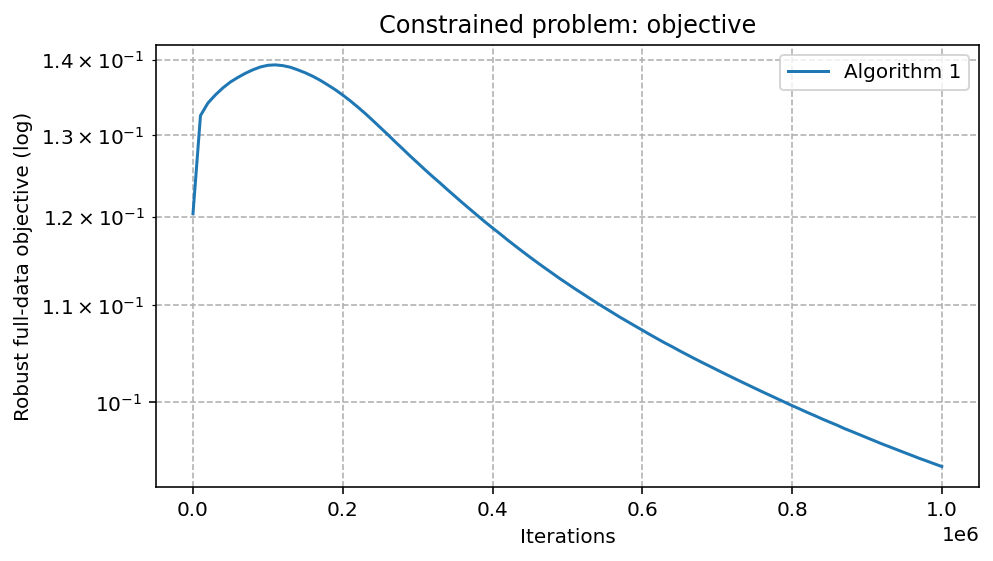}
\caption{Full-data objective value $f(u)$.}
\label{fig:exp2-objective}
\end{subfigure}
\hfill
\begin{subfigure}{0.48\textwidth}
\centering
\includegraphics[width=\textwidth]{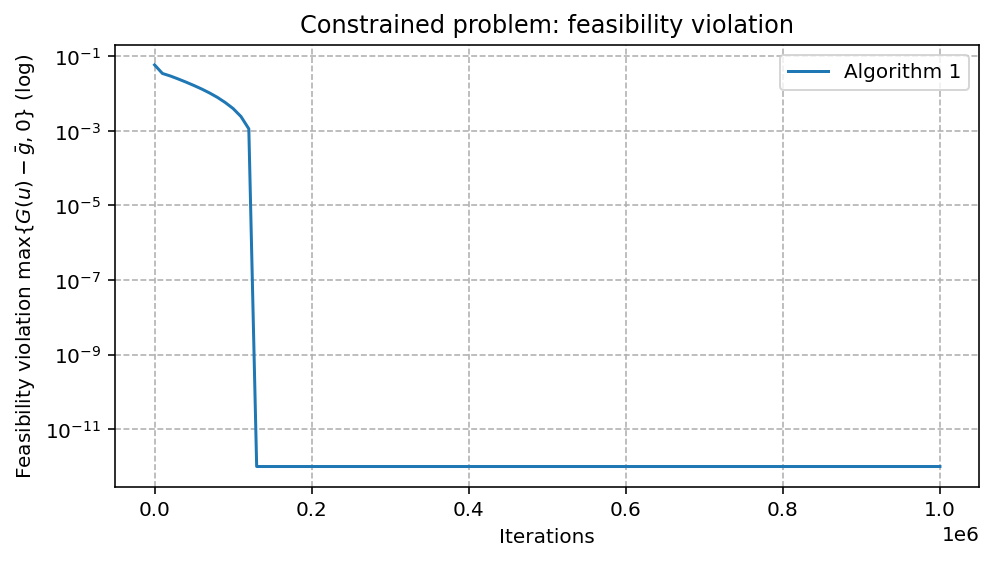}
\caption{Feasibility violation $\max\{g(u) - \bar g, \epsilon\}$.}
\label{fig:exp2-feas}
\end{subfigure}

\vspace{0.3cm}

\begin{subfigure}{0.48\textwidth}
\centering
\includegraphics[width=\textwidth]{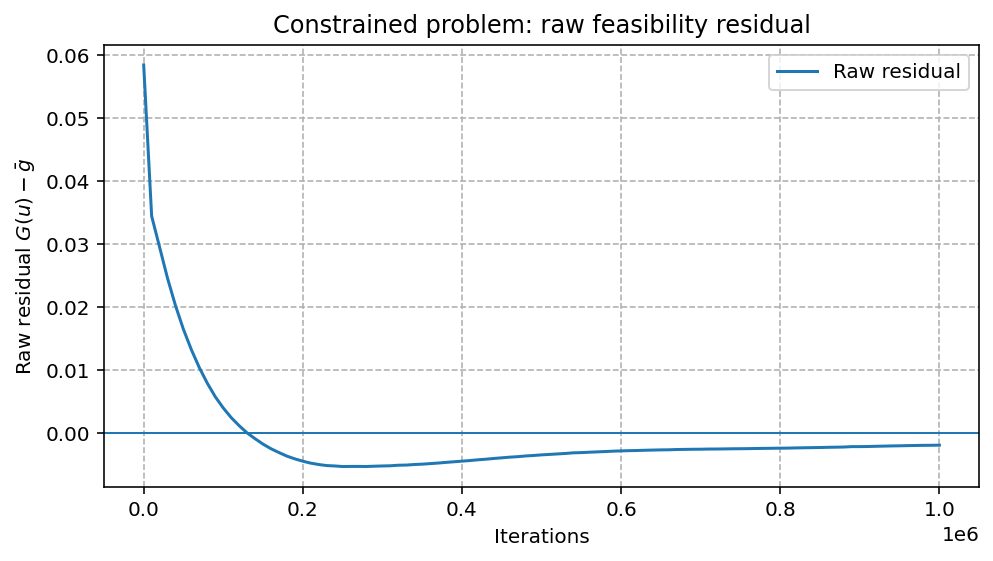}
\caption{The constraint value $g(u) - \bar g$}
\label{fig:exp2-raw}
\end{subfigure}
\hfill
\begin{subfigure}{0.48\textwidth}
\centering
\includegraphics[width=\textwidth]{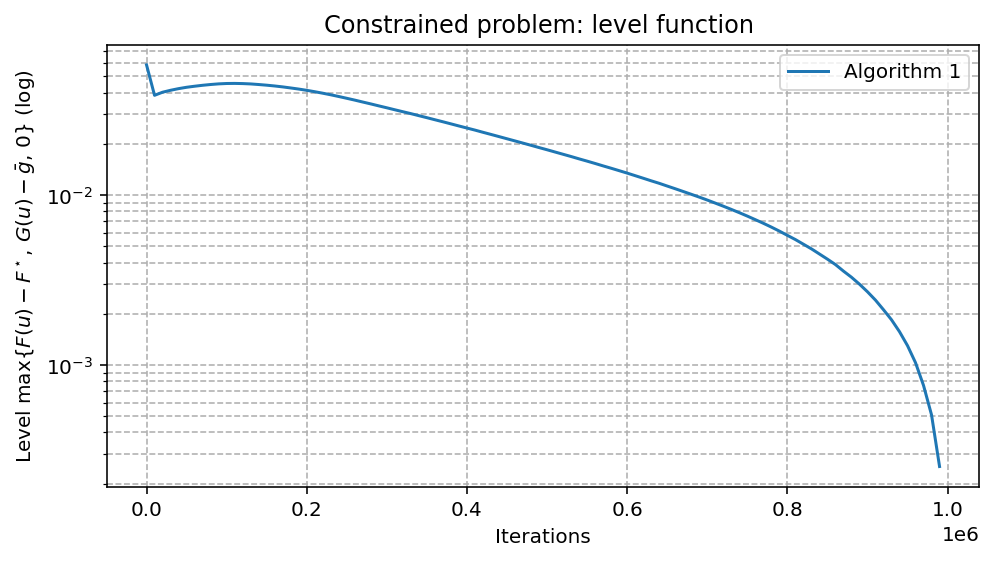}
\caption{Level function $\max\{f(u) -f^*, g(u)-\bar g\}$.}
\label{fig:exp2-level}
\end{subfigure}

\vspace{0.3cm}

\begin{subfigure}{0.48\textwidth}
\centering
\includegraphics[width=\textwidth]{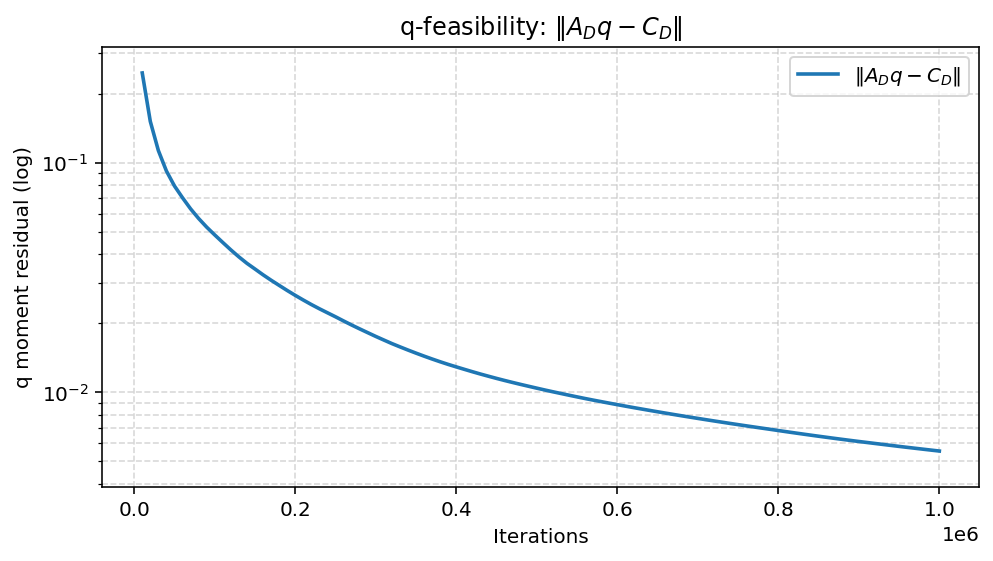}
\caption{\(q\)-feasibility residual \(\|A_Dq-C_D\|\).}
\label{fig:exp2-qres}
\end{subfigure}
\hfill
\begin{subfigure}{0.48\textwidth}
\centering
\includegraphics[width=\textwidth]{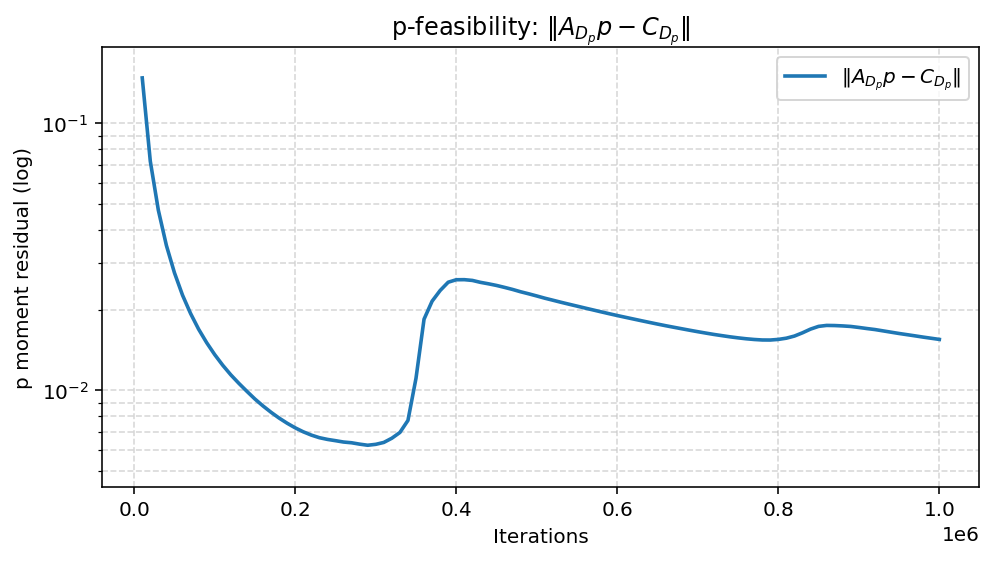}
\caption{\(p\)-feasibility residual \(\|A_{D_p}p-C_{D_p}\|\)}
\label{fig:exp2-pres}
\end{subfigure}

\caption{Performance of the proposed primal-dual method for the constrained moment-based DRO problem on the Mushroom dataset.}
\label{fig:exp2-results}
\end{figure}

One can set $\bar g := \kappa g^\star(\bar u)$ where $\kappa \in (0,1)$ is a prescribed scaling factor. We use $\kappa = 0.5$ in our experiment. 
This construction creates a stricter constraint level and keeps the DRO constraint active in the following problem:
\begin{equation}
\min_{u\in U} f(u) \quad \text{s.t.} \quad g(u)\leq \bar g . 
\label{eq:exp2-constrained-dro} 
\end{equation} 
\textbf{Experimental Setup:}
The Mushroom dataset contains \(8124\) samples. Each sample is represented as \((x_i,y_i)\), where \(x_i\) denotes the feature vector and \(y_i\in\{-1,1\}\) denotes the binary class label. We use the logistic loss 
\begin{equation} \ell_i(u) = \log\left(1+\exp(-y_i x_i^\top u)\right), 
\label{eq:exp2-logistic-loss} 
\end{equation} 
where \(u\) is the classifier parameter. The moment-based DRO problem is implemented through a primal-dual saddle-point reformulation. The probability vectors \(q\) and \(p\) are projected onto the simplexes \(\Delta_N\) and \(\Delta_n\), respectively. The moment constraints \(A_Dq=C_D\) and \(A_{D_p}p=C_{D_p}\) are incorporated through dual variables. A nonnegative multiplier \(\lambda\geq 0\) is used to enforce the DRO constraint \(g(u)\leq \bar g\). The implemented saddle formulation is \begin{equation}
\min_{u,\mu_q,\mu_p}
\max_{\substack{q\in\Delta_N,\;p\in\Delta_n\\ \lambda\geq 0}}
[
\wt{f}(u,q)
+\mu_q^\top(A_Dq-C_D)
+\lambda
(
\wt{g}(u,p)-\bar g
+\mu_p^\top(A_{D_p}p-C_{D_p})
)
],
\label{eq:exp2-saddle}
\end{equation}
Here, \(\mu_q\) and \(\mu_p\) are the dual variables associated with the moment constraints for the full dataset and the poisonous subset, respectively. To evaluate the performance of the PD-ConEx method, we track the feasibility violation and the level function. The feasibility violation is defined as 
\begin{equation} \max\{g(u)-\bar g,0\}, \label{eq:exp2-feas-violation} 
\end{equation} where \(g(u)\) denotes the worst-case poisonous-subset loss. The level function is defined by 
\begin{equation} 
\mathrm{Level}(u) = \max\{f(u)-f^\star,\ g(u)-\bar g,\ 0\}. 
\label{eq:exp2-level} 
\end{equation} This quantity combines the objective residual and the constraint violation into a single measure, and therefore provides a convenient way to monitor both optimality and feasibility during the iterations.\\ 
\textbf{Results:}
We implement our method using the step-size policy presented in Theorem \ref{thm:monotone_conv}. The numerical results are presented in Figure~\ref{fig:exp2-results}. The proposed primal-dual method was run for \(10^6\) iterations on the moment-based DRO problem over the Mushroom dataset. 
All of the quantities are evaluated along the averaged iterates generated by the algorithm. The level function is computed as \eqref{eq:exp2-level}

The objective plot in Figure~\ref{fig:exp2-objective} shows that the full-data objective initially increases during the early stage. Looking at the feasibility plot Figure~\ref{fig:exp2-feas}, the algorithm is making progress to satisfy the poisonous-subset DRO constraint and to the associated dual variables. After this initial phase, the objective decreases steadily over the remaining iterations, indicating progress towards optimality inside the feasible region of the problem. 
The feasibility behavior is shown in Figure~\ref{fig:exp2-feas}. We calculate the feasibility violation as $\max\{g (u) - \bar g, 10^{-12}\}$. The violation is capped at $10^{-12}$ for plotting on the log-scale. As one can see, the feasibility violation decreases rapidly during the early iterations and reaches the numerical plotting floor. After this point, it remains essentially zero for the rest of the computation.
The exact constraint values are reported in Figure~\ref{fig:exp2-raw}. It shows the constraint value goes below zero and then approaches zero from negative side as we make progress towards optimality.

The level function is shown in Figure~\ref{fig:exp2-level}. The estimated best feasible objective value obtained from the run is \(f^\star \approx 0.093879079\). For visualization, the final evaluation point is excluded from the level-function plot, since the level value at the final iterate is zero. The last plotted point corresponds to iteration \(990000\), where the level value is approximately \(2.52\times 10^{-4}\).

Since the feasibility violation reaches zero after the early stage of the run, the later behavior of the level function is mainly governed by the objective residual \(F(u)-F^\star\). The decreasing trend in the plot shows that, after achieving feasibility, the algorithm continues to improve the objective value within the feasible region. Thus, the level function confirms the combined progress of the method in terms of both feasibility restoration and objective reduction.
In addition to the objective and feasibility measures, we also monitor the moment residuals associated with the adversarial probability vectors. Specifically, we plot
\[
\|A_D q-C_D\|
\quad \text{and} \quad
\|A_{D_p}p-C_{D_p}\|,
\]
which measure how well the averaged algorithmic variables \(q\) and \(p\) satisfy the moment constraints defining the ambiguity sets. The \(q\)-residual corresponds to the full-data ambiguity set, while the \(p\)-residual corresponds to the poisonous-subset ambiguity set.

The moment residual plots in Figures~\ref{fig:exp2-qres} and~\ref{fig:exp2-pres} show that both residuals remain small over the course of the algorithm. The \(q\)-residual decreases steadily after an initial transient phase, indicating that the full-data moment constraint is increasingly well enforced. The \(p\)-residual also remains small, with a temporary increase during the middle iterations before stabilizing again.At the final recorded iterate, the $q$- and $p$-moment residuals are approximately
$5.546\times 10^{-3}$ and $1.551\times 10^{-2}$, respectively. These plots confirm that the algorithm not only reduces the objective and feasibility violation, but also maintains approximate feasibility of the moment constraints for the adversarial probability variables.

\section{Conclusion}

In this paper, we studied variational inequality problems with semi-infinite function constraints, where direct projection onto the feasible set is generally difficult due to the presence of infinitely many constraints. To address this challenge, we proposed a primal-dual first-order method that enforces feasibility through Lagrange multipliers, constraint-index variables, linearization, and extrapolation, while requiring projection only onto simple sets. We established convergence guarantees under monotonicity, convexity, and smoothness or nonsmoothness assumptions, and further discussed extensions to semi-infinite saddle-point problems with coupling constraints.

We also showed that moment-based and Wasserstein distributionally robust function-constrained optimization problems can be reformulated within the proposed semi-infinite framework. This connection demonstrates the applicability of the method to robust learning problems. The numerical experiments support the theoretical development. In the robust quadratic examples, the PD-ConEx method showed faster reduction of the level function and stronger feasibility behavior compared with SGSP. In the moment-based DRO experiment on the Mushroom dataset, the method rapidly reduced the feasibility violation, maintained approximate moment feasibility of the adversarial probability variables, and continued to improve the objective value after reaching feasibility.

Overall, the results indicate that the proposed primal-dual framework is effective for solving variational inequality and saddle-point problems with semi-infinite constraints, especially in settings where projection onto the full feasible set is computationally difficult.

\section{Appendix}
\subsection{Proof of 
Proposition~\ref{prop:summed-gap-alg2}}
\label{app:proof-gap-alg2}
\begin{sharedproof}
Let \(x_k=(u_k,v_k)\), and define
$
\Delta F_k^u:=F_u(x_k)-F_u(x_{k-1}),
\quad
\Delta F_k^v:=F_v(x_k)-F_v(x_{k-1}).
$
Recall from \eqref{eq:5.3} that
\[
F(x)=
((F_u(x),F_v(x))
=
(\nabla_u f(u,v),-\nabla_v f(u,v)).
\]
The proof follows the same structure as the proof of Algorithm \ref{alg:vi-sip-agsip}, but now
the primal variable is split into the two blocks \(u\) and \(v\).

First, applying the three-point inequality used in \eqref{eq:4.10} to the
\(u\) update in Algorithm \ref{alg:prob2}, and using the same argument that leads to
\eqref{eq:4.14}, we obtain
\begin{align}
&\left\langle F_u(x),u_{k+1}-u\right\rangle
+
\left\langle
g'_u(x_k,y_{k+1})^\top\lambda_{k+1},
u_{k+1}-u
\right\rangle
\le
\tfrac{\tau_k^u}{2}\|u-u_k\|^2
-
(
\tfrac{\tau_k^u}{2}+\mu_F^u
)\|u-u_{k+1}\|^2
\nonumber\\
&-
\tfrac{\tau_k^u}{2}\|u_{k+1}-u_k\|^2
+
\left\langle \Delta F_{k+1}^u,u_{k+1}-u\right\rangle
-
\theta_k\left\langle \Delta F_k^u,u_k-u\right\rangle
-
\theta_k\left\langle \Delta F_k^u,u_{k+1}-u_k\right\rangle .
\label{eq:7.1}
\end{align}
This is the \(u\)-block analogue of \eqref{eq:4.14}.

Similarly, applying the same argument to the \(v\)-update gives
\begin{align}
&\left\langle F_v(x),v_{k+1}-v\right\rangle
+
\left\langle
g'_v(x_k,y_{k+1})^\top\lambda_{k+1},
v_{k+1}-v
\right\rangle
\le
\tfrac{\tau_k^v}{2}\|v-v_k\|^2
-
(
\tfrac{\tau_k^v}{2}+\mu_F^v
)\|v-v_{k+1}\|^2
\nonumber\\
&-
\tfrac{\tau_k^v}{2}\|v_{k+1}-v_k\|^2
+
\left\langle \Delta F_{k+1}^v,v_{k+1}-v\right\rangle
-
\theta_k\left\langle \Delta F_k^v,v_k-v\right\rangle
-
\theta_k\left\langle \Delta F_k^v,v_{k+1}-v_k\right\rangle .
\label{eq:7.2}
\end{align}

By convexity of \(f(\cdot,v)\) in \(u\) and concavity of \(f(u,\cdot)\)
in \(v\), the two block estimates imply
\begin{align}
&f(u_{k+1},v)-f(u,v_{k+1})
+
\left\langle
g'_u(x_k,y_{k+1})^\top\lambda_{k+1},
u_{k+1}-u
\right\rangle
+
\left\langle
g'_v(x_k,y_{k+1})^\top\lambda_{k+1},
v_{k+1}-v
\right\rangle
\nonumber\\
&\le
\tfrac{\tau_k^u}{2}\|u-u_k\|^2
-
(
\tfrac{\tau_k^u}{2}+\mu_F^u
)\|u-u_{k+1}\|^2
-
\tfrac{\tau_k^u}{2}\|u_{k+1}-u_k\|^2
+
\tfrac{\tau_k^v}{2}\|v-v_k\|^2
-
(
\tfrac{\tau_k^v}{2}+\mu_F^v
)\|v-v_{k+1}\|^2 \nonumber \\
&-
\tfrac{\tau_k^v}{2}\|v_{k+1}-v_k\|^2
+
\left\langle \Delta F_{k+1}^u,u_{k+1}-u\right\rangle
-
\theta_k\left\langle \Delta F_k^u,u_k-u\right\rangle
-
\theta_k\left\langle \Delta F_k^u,u_{k+1}-u_k\right\rangle
\nonumber\\
&
+
\left\langle \Delta F_{k+1}^v,v_{k+1}-v\right\rangle
-
\theta_k\left\langle \Delta F_k^v,v_k-v\right\rangle
-
\theta_k\left\langle \Delta F_k^v,v_{k+1}-v_k\right\rangle .
\label{eq:7.3}
\end{align}

Next, by convexity of \(g(\cdot,\cdot,y_{k+1})\) in \(x=(u,v)\), and
using the same linearization argument as in \eqref{eq:4.22}, we have
\[
\begin{aligned}
&\left\langle
g'_u(x_k,y_{k+1})^\top\lambda_{k+1},
u_{k+1}-u
\right\rangle
+
\left\langle
g'_v(x_k,y_{k+1})^\top\lambda_{k+1},
v_{k+1}-v
\right\rangle
\ge
\left\langle
\lambda_{k+1},
\ell_g(x_{k+1};x_k,y_{k+1})-g(x,y_{k+1})
\right\rangle
\end{aligned}
\]
Substituting this into \eqref{eq:7.3} gives the
Algorithm \ref{alg:prob2} analogue of the one-step inequality used in the proof of Algorithm \ref{alg:vi-sip-agsip} before summing over \(k\)
Since the \(\lambda\) update is unchanged from Algorithm \ref{alg:vi-sip-agsip}, the same \(\lambda\) estimate used in \eqref{eq:4.18} applies.
Thus, we obtain
\begin{align}
&-\left\langle
\ell_g(x_{k+1};x_k,y_{k+1}),
\lambda_{k+1}-\lambda
\right\rangle
\le
\tfrac{\gamma_k}{2}
(
\|\lambda-\lambda_k\|^2
-
\|\lambda-\lambda_{k+1}\|^2
)
-
\tfrac{\gamma_k}{2}\|\lambda_{k+1}-\lambda_k\|^2
\nonumber\\
&
-
\left\langle q_\lambda^{k+1},\lambda_{k+1}-\lambda\right\rangle
+
\theta_k\left\langle q_\lambda^k,\lambda_k-\lambda\right\rangle
+
\theta_k\left\langle q_\lambda^k,\lambda_{k+1}-\lambda_k\right\rangle .
\label{eq:7.4}
\end{align}

The \(y\)-update is also unchanged from Algorithm \ref{alg:vi-sip-agsip}. Hence, using Lemma \ref{lemma3},
we obtain
\begin{align}
&\left\langle
\lambda,
g(x_{k+1},y)-g(x_{k+1},y_{k+1})
\right\rangle
\le
-
\left\langle \Lambda q_y^{k+1},y_{k+1}-y\right\rangle
+
\theta_k\left\langle \Lambda q_y^k,y_k-y\right\rangle
+
\theta_k\left\langle \Lambda q_y^k,y_{k+1}-y_k\right\rangle
\nonumber\\
&
-
\tfrac{\sigma_k}{2}\|y_{k+1}-y_k\|_\lambda^2
+
\tfrac{\sigma_k}{2}\|y-y_k\|_\lambda^2
-
\tfrac{\sigma_k}{2}\|y-y_{k+1}\|_\lambda^2
-
\tfrac{\mu_y}{2}\|y_{k+1}-y\|_\lambda^2 .
\label{eq:7.5}
\end{align}

Adding \eqref{eq:7.3},
\eqref{eq:7.4}, and \eqref{eq:7.5}, and using
the Algorithm \ref{alg:prob2} gap definition gives the
one-step bound for \(Q(z_{k+1},z)\),  multiply by \(t_k\) and sum up \(k=0\) to
\(k=i\).
The remaining cross terms are then bounded exactly as in
\eqref{eq:4.48}--\eqref{eq:4.53}, but now applied separately to the \(u\)- and \(v\)-blocks.

Specifically, the multiplier extrapolation term is bounded using the same
estimate as \eqref{eq:4.48}. The operator extrapolation term is bounded by applying
the analogue of \eqref{eq:4.49} separately as
\[
-t_{k-1}\left\langle
\Delta F_k^u,u_{k+1}-u_k
\right\rangle
\le
\tfrac{t_{k-1}L_F^u}{2}\|u_k-u_{k-1}\|^2
+
\tfrac{t_{k-1}L_F^u}{2}\|u_{k+1}-u_k\|^2,
\]
\[
-t_{k-1}\left\langle
\Delta F_k^v,v_{k+1}-v_k
\right\rangle
\le
\tfrac{t_{k-1}L_F^v}{2}\|v_k-v_{k-1}\|^2
+
\tfrac{t_{k-1}L_F^v}{2}\|v_{k+1}-v_k\|^2.
\]
The \(y\)-gradient extrapolation term is bounded using the same logic as
(2.39)--(2.42),  with 
replaced by the blockwise bound
$
L_g^{yu}\|u_k-u_{k-1}\|
+
L_g^{yv}\|v_k-v_{k-1}\|.
$
The nonsmooth terms are handled by Young's inequality as in (2.44)--(2.45):
\[
t_kH_g^u\|\lambda\|_1\|u_{k+1}-u_k\|
\le
\tfrac{4t_k(H_g^u\|\lambda\|_1)^2}{\tau_k^u}
+
\tfrac{t_k\tau_k^u}{16}\|u_{k+1}-u_k\|^2,
\]
\[
t_kH_g^v\|\lambda\|_1\|v_{k+1}-v_k\|
\le
\tfrac{4t_k(H_g^v\|\lambda\|_1)^2}{\tau_k^v}
+
\tfrac{t_k\tau_k^v}{16}\|v_{k+1}-v_k\|^2,
\]
\[
t_{k-1}H_g^y\|\lambda\|_1^{1/2}
\|y_{k+1}-y_k\|_\lambda
\le
\tfrac{4t_{k-1}(H_g^y)^2\|\lambda\|_1}{\sigma_k}
+
\tfrac{t_{k-1}\sigma_k}{16}
\|y_{k+1}-y_k\|_\lambda^2 .
\]
By the definition of \(\ell_g\) and the
block smoothness of \(g\),
\begin{align}
&\left\langle
\lambda,
g(u_{k+1},v_{k+1},y_{k+1})
-
\ell_g(u_{k+1},v_{k+1};u_k,v_k,y_{k+1})
\right\rangle
\le
(
L_g^{uu}\|\lambda\|_1
+
\tfrac{L_g^{uv}\|\lambda\|_1}{2}
)
\|u_{k+1}-u_k\|^2
\nonumber\\
&
+
(
L_g^{vv}\|\lambda\|_1
+
\tfrac{L_g^{uv}\|\lambda\|_1}{2}
)
\|v_{k+1}-v_k\|^2
+
H_g^u\|\lambda\|_1\|u_{k+1}-u_k\|
+
H_g^v\|\lambda\|_1\|v_{k+1}-v_k\|.
\label{eq:7.6}
\end{align}
The y-gradient extrapolation term is bounded using \eqref{eq:4.50}–\eqref{eq:4.53}, substituting the cross-term estimates \eqref{eq:4.48}, \eqref{eq:4.49}, and \eqref{eq:4.53}, together with the Young-inequality bounds, and rearranging as in \eqref{eq:4.63},  gives the gap bound and set \(i=K-1\). Recall
\[
\Gamma_K:=\tsum_{k=0}^{K-1}t_k,
\qquad
\bar u_K=\tfrac{1}{\Gamma_K}\tsum_{k=0}^{K-1}t_k u_{k+1},
\qquad
\bar v_K=\tfrac{1}{\Gamma_K}\tsum_{k=0}^{K-1}t_k v_{k+1}.
\]
Using the convexity of \(f(\cdot,v)\) in \(u\) and concavity of  \(f(u,\cdot)\) in \(v\) we can write 
\begin{align}
\Gamma_K
[
f(\bar u_K,v)-f(u,\bar v_K)
]
&\le
\tsum_{k=0}^{K-1}t_k
[
f(u_{k+1},v)-f(u,v_{k+1})
].
\label{eq:7.7}
\end{align}

Similarly, by convexity of \(g(\cdot,\cdot,y)\) in \(x=(u,v)\) and
because \(\lambda\ge0\),
\begin{align}
\Gamma_K
\left\langle
g(\bar u_K,\bar v_K,y),\lambda
\right\rangle
&\le
\tsum_{k=0}^{K-1}t_k
\left\langle
g(u_{k+1},v_{k+1},y),\lambda
\right\rangle .
\label{eq:7.8}
\end{align}
Combining \eqref{eq:7.7} and
\eqref{eq:7.8} with the gap bound , we obtain the following
complete summed bound before applying the step-size conditions:
\begin{align}
&\Gamma_K
[
f(\bar u_K,v)-f(u,\bar v_K)
]
+
\Gamma_K
\left\langle
g(\bar u_K,\bar v_K,y),\lambda
\right\rangle
-
\tsum_{k=0}^{K-1}t_k
\left\langle
g(u,v,y_{k+1}),\lambda_{k+1}
\right\rangle
\nonumber\\
&\le
\tfrac{t_0\tau_0^u}{2}\|u-u_0\|^2
+
\tfrac{t_0\tau_0^v}{2}\|v-v_0\|^2
+
\tfrac{t_0\gamma_0}{2}\|\lambda-\lambda_0\|^2
+
\tfrac{t_0\sigma_0}{2}\|y-y_0\|_\lambda^2
\nonumber\\
&
+
[
\tfrac{t_{K-1}L_F^u}{2}
-
t_{K-1}
(
\tfrac{\tau_{K-1}^u}{2}+\mu_F^u
)
]\|u_K-u\|^2
+
[
\tfrac{t_{K-1}L_F^v}{2}
-
t_{K-1}
(
\tfrac{\tau_{K-1}^v}{2}+\mu_F^v
)
]\|v_K-v\|^2
\nonumber\\
&
+
[
\tfrac{t_{K-1}\gamma_K}{4}
-
\tfrac{t_{K-1}\gamma_{K-1}}{2}
]\|\lambda_K-\lambda\|^2
+
[
\tfrac{3t_{K-1}\sigma_K}{16}
-
\tfrac{t_{K-1}(\sigma_{K-1}+\mu_y)}{2}
]\|y_K-y\|_\lambda^2
\nonumber\\
&
+
\tsum_{k=1}^{K-1}
[
\tfrac{t_{k-1}\gamma_k}{8}
-
\tfrac{t_k\gamma_k}{2}
]\|\lambda_{k+1}-\lambda_k\|^2
\nonumber\\
&
+
\tsum_{k=1}^{K-1}
[
\tfrac{6t_kM_g^2}{\gamma_{k+1}}
+
t_kL_F^u
+
\tfrac{t_{k-1}L_F^u}{2}
+
\tfrac{24t_{k+1}M_g^2}{\gamma_{k+2}}
-
\tfrac{7t_k\tau_k^u}{32}
]\|u_{k+1}-u_k\|^2
\nonumber\\
&
+
\tsum_{k=1}^{K-1}
[
\tfrac{6t_kM_g^2}{\gamma_{k+1}}
+
t_kL_F^v
+
\tfrac{t_{k-1}L_F^v}{2}
+
\tfrac{24t_{k+1}M_g^2}{\gamma_{k+2}}
-
\tfrac{7t_k\tau_k^v}{32}
]\|v_{k+1}-v_k\|^2
\nonumber\\
&
+
\Psi_u(\|\lambda\|_1,K)
+
\Psi_v(\|\lambda\|_1,K)
+
\tsum_{k=1}^{K-1}
[
\tfrac{8t_k(L_g^{yy})^2}{\sigma_{k+1}}
+
\tfrac{3t_{k-1}\sigma_k}{16}
-
\tfrac{t_k\sigma_k}{2}
]\|y_{k+1}-y_k\|_\lambda^2
\nonumber\\
&
+
\tsum_{k=0}^{K-1}
[
\tfrac{4t_k(H_g^u\|\lambda\|_1)^2}{\tau_k^u}
+
\tfrac{4t_k(H_g^v\|\lambda\|_1)^2}{\tau_k^v}
]
+
4
(
\tsum_{k=1}^{K-1}\tfrac{t_k}{\sigma_{k+1}}
+
\tsum_{k=1}^{K-1}\tfrac{t_{k-1}}{\sigma_k}
)
(H_g^y)^2\|\lambda\|_1 .
\label{eq:7.9}
\end{align}
Then we can obtained the step size conditions stated as in Propositions~\ref{prop:summed-gap-alg2}.
Therefore, all the above distance and increment terms can be dropped from the upper bound in \eqref{eq:7.9}. The only
remaining terms are the initial distance terms, the two residual terms
\(\Psi_u\) and \(\Psi_v\), and the nonsmoothness terms. Hence,
\begin{align}
&\Gamma_K
[
f(\bar u_K,v)-f(u,\bar v_K)
]
+
\Gamma_K
\left\langle
g(\bar u_K,\bar v_K,y),\lambda
\right\rangle
-
\tsum_{k=0}^{K-1}t_k
\left\langle
g(u,v,y_{k+1}),\lambda_{k+1}
\right\rangle
\nonumber\\
&\le
\tfrac{t_0\tau_0^u}{2}\|u-u_0\|^2
+
\tfrac{t_0\tau_0^v}{2}\|v-v_0\|^2
+
\tfrac{t_0\gamma_0}{2}\|\lambda-\lambda_0\|^2
+
\tfrac{t_0\sigma_0}{2}\|y-y_0\|_\lambda^2
\nonumber\\
&
+
\Psi_u(\|\lambda\|_1,K)
+
\Psi_v(\|\lambda\|_1,K)
\nonumber\\
&
+
\tsum_{k=0}^{K-1}
[
\tfrac{4t_k(H_g^u\|\lambda\|_1)^2}{\tau_k^u}
+
\tfrac{4t_k(H_g^v\|\lambda\|_1)^2}{\tau_k^v}
]
+
4
(
\tsum_{k=1}^{K-1}\tfrac{t_k}{\sigma_{k+1}}
+
\tsum_{k=1}^{K-1}\tfrac{t_{k-1}}{\sigma_k}
)
(H_g^y)^2\|\lambda\|_1 .
\label{eq:7.10}
\end{align}
Finally, since \(\|\lambda\|_1\le\Lambda\), we replace
\(\|\lambda\|_1\) by \(\Lambda\) on the right-hand side. This gives
exactly \eqref{eq:gap_prop3}, where \(\Psi_u(\Lambda,K)\) and \(\Psi_v(\Lambda,K)\)
are defined in \eqref{eq:T-Lambda-K_u}--\eqref{eq:T-Lambda-K_v}. Therefore, Propositions~\ref{prop:summed-gap-alg2} follows.
\end{sharedproof}
\subsection{Adaptive PD-Lagrangian Extrapolation Theorem~\ref{thm:adaptive-pdconex}}
\begin{sharedproof}
\label{app:adaptive-proof}
Let
$z=(x,\lambda,y)$ and
$z^{k+1}=(x_{k+1},\lambda_{k+1},y_{k+1})$.
Recall $
Q(z^{k+1},z)
=
\langle F(x),x_{k+1}-x\rangle
+\lambda^\top g(x_{k+1},y)
-(\lambda_{k+1})^\top g(x,y_{k+1}).
$
We decompose$
Q(z^{k+1},z)=Q_0+Q_1+Q_2,
$
where
\begin{align*}
Q_0
&:=
\lambda^\top
[g(x_{k+1},y)-g(x_{k+1},y_{k+1})],\\
Q_1
&:=
(\lambda-\lambda_{k+1})^\top
g(x_{k+1},y_{k+1}),\\
Q_2
&:=
\langle F(x),x_{k+1}-x\rangle
+(\lambda_{k+1})^\top
[g(x_{k+1},y_{k+1})-g(x,y_{k+1})].
\end{align*}
Define
$
q_y^k
:=
g_y'(x_k,y_k)-g_y'(x_{k-1},y_{k-1}),
q_\lambda^k
:=
g(x_k,y_k)-g(x_{k-1},y_k),
q_x^k
:=
P(x_k,y_k,\lambda_k)
-P(x_{k-1},y_k,\lambda_{k-1}),
$
with
$q_y^0=q_\lambda^0=q_x^0=0$.
\paragraph{Blockwise estimates.}
The optimality condition of the $y$-update and the
$\mu_y$-strong concavity of $g_i(x,\cdot)$ give
\begin{align}
Q_0
\le&
-\langle\Lambda q_y^{k+1},y_{k+1}-y\rangle
+\theta_k\langle\Lambda q_y^k,y_k-y\rangle
+\theta_k\langle\Lambda q_y^k,y_{k+1}-y_k\rangle
\nonumber\\
&+
\tfrac{\sigma_k}{2}
[
\|y-y_k\|_\lambda^2-\|y-y_{k+1}\|_\lambda^2
]
-\tfrac{\sigma_k}{2}\|y_{k+1}-y_k\|_\lambda^2
-\tfrac{\mu_y}{2}\|y_{k+1}-y\|_\lambda^2 .
\label{eq:app-adap-Q0}
\end{align}
Similarly, optimality of the multiplier update gives
\begin{align}
Q_1
\le&
\tfrac{\gamma_k}{2}
[
\|\lambda-\lambda_k\|_2^2
-\|\lambda-\lambda_{k+1}\|_2^2
]
-\tfrac{\gamma_k}{2}
\|\lambda_{k+1}-\lambda_k\|_2^2
\nonumber\\
&-
\langle q_\lambda^{k+1},\lambda_{k+1}-\lambda\rangle
+\theta_k
\langle q_\lambda^k,\lambda_k-\lambda\rangle
+\theta_k
\langle q_\lambda^k,\lambda_{k+1}-\lambda_k\rangle .
\label{eq:app-adap-Q1}
\end{align}
Finally, monotonicity of $F$, convexity of
$g_i(\cdot,y_i)$, and the primal proximal optimality condition imply
\begin{align}
Q_2
\le&
\tfrac{\tau_k}{2}
[\|x-x_k\|^2-\|x-x_{k+1}\|^2]
-\tfrac{\tau_k}{2}
\|x_{k+1}-x_k\|^2+
\langle q_x^{k+1},x_{k+1}-x\rangle\nonumber\\&
-\theta_k\langle q_x^k,x_k-x\rangle
-\theta_k\langle q_x^k,x_{k+1}-x_k\rangle .
\label{eq:app-adap-Q2}
\end{align}
The smoothness assumptions further imply
\begin{align}
\|q_\lambda^k\|_2
&\le
M_g\|x_k-x_{k-1}\|,
\nonumber\\
\|q_x^k\|
&\le
[L_F+L_g^{xx}\|\lambda_{k-1}\|_2]
\|x_k-x_{k-1}\|
+
M_g\|\lambda_k-\lambda_{k-1}\|_2,
\nonumber\\
\|q_y^k\|_\lambda^2
&\le
2(L_g^{yx})^2\|\lambda\|_1
\|x_k-x_{k-1}\|^2
+
2(L_g^{yy})^2
\|y_k-y_{k-1}\|_\lambda^2.\nonumber
\end{align}
Let
$
d_x^k:=x_k-x_{k-1},\quad
d_\lambda^k:=\lambda_k-\lambda_{k-1},\quad
d_y^k:=y_k-y_{k-1}.
$
Multiplying
\eqref{eq:app-adap-Q0}--\eqref{eq:app-adap-Q2}
by $t_k$ and using
$t_k\tau_k=a_0,\quad
t_k\gamma_k=\beta a_0,\quad
t_k\sigma_k=b_k,\quad
t_k\theta_k=t_{k-1},
$
gives
\begin{align}
&t_kQ(z^{k+1},z)
\le
\tfrac{a_0}{2}
[
\|x-x_k\|^2-\|x-x_{k+1}\|^2
]
+
\tfrac{\beta a_0}{2}
[
\|\lambda-\lambda_k\|_2^2
-\|\lambda-\lambda_{k+1}\|_2^2
]
\nonumber\\
&+
\tfrac{b_k}{2}
[
\|y-y_k\|_\lambda^2
-\|y-y_{k+1}\|_\lambda^2
]
-\tfrac{\mu_yt_k}{2}
\|y-y_{k+1}\|_\lambda^2-
\tfrac{a_0}{2}\|d_x^{k+1}\|^2
-\tfrac{\beta a_0}{2}\|d_\lambda^{k+1}\|_2^2
-\tfrac{b_k}{2}\|d_y^{k+1}\|_\lambda^2
\nonumber\\
&+
\mathcal R_x^k+\mathcal R_\lambda^k+\mathcal R_y^k
+\mathcal T_x^k+\mathcal T_\lambda^k+\mathcal T_y^k .
\label{eq:app-adap-one-step}
\end{align}
Here the $\mathcal R$-terms collect the local cross terms and the
$\mathcal T$-terms collect the extrapolation terms.
Using \eqref{eq:adap-step-control},
Cauchy--Schwarz, and Young's inequality,
\begin{align}
|\mathcal R_x^k|
\le&
\tfrac{a_0}{2}
(\tfrac1{c_1}+\tfrac1{c_2})
(\|d_x^k\|^2+\|d_x^{k+1}\|^2)+
\tfrac{a_0}{2c_1}
(\|d_\lambda^k\|_2^2+\|d_x^{k+1}\|^2),
\label{eq:app-adap-Rx}\\
|\mathcal R_\lambda^k|
\le&
\tfrac{a_0}{2c_1}
(\|d_x^k\|^2+\|d_\lambda^{k+1}\|_2^2).
\label{eq:app-adap-Rlambda}\\
|\mathcal R_y^k|
\le&
\tfrac{b_k}{4}\|d_y^{k+1}\|_\lambda^2
+
\tfrac{
2(L_g^{yx})^2\|\lambda\|_1t_{k-1}^2
}{b_k}
\|d_x^k\|^2+
\tfrac{
2(L_g^{yy})^2t_{k-1}^2
}{b_k}
\|d_y^k\|_\lambda^2.
\label{eq:app-adap-Ry}
\end{align}
Because $b_0\ge4L_g^{yy}$ and $b_k$ is nondecreasing,
$
\tfrac{
2(L_g^{yy})^2t_{k-1}^2
}{b_k}
\le
\tfrac{b_{k-1}}8.
$\\
Let
\[
A_\lambda
:=
2(L_g^{yx})^2\|\lambda\|_1,
\qquad
E_k(\lambda)
:=
[
A_\lambda\tfrac{t_{k-1}^2}{b_k}
-\rho_xa_0
]_+.
\]
If $E_k(\lambda)>0$, then
$b_k<A_\lambda/(\rho_xa_0)$. Since
$b_k-b_{k-1}=\mu_yt_{k-1}$ and
$t_{k-1}^2\le t_{k-1}$,
\begin{align*}
\tsum_{\{k:E_k(\lambda)>0\}}
A_\lambda\tfrac{t_{k-1}^2}{b_k}
&\le
\tfrac{A_\lambda}{\mu_y}
\tsum_{\{k:E_k(\lambda)>0\}}
\tfrac{b_k-b_{k-1}}{b_k}
\le
\tfrac{A_\lambda}{\mu_y}
[
\log(
\tfrac{A_\lambda}{\rho_xa_0b_0}
)
]_+ .
\end{align*}
Using $\|d_x^k\|\le D_X$, we obtain
\begin{equation}
\tsum_{k=1}^{K-1}
E_k(\lambda)\|d_x^k\|^2
\le
C_{yx}(\lambda),
\label{eq:app-adap-Cyx}
\end{equation}
where $C_{yx}(0):=0$ and, for $\lambda\neq0$,
$
C_{yx}(\lambda)
:=
\tfrac{A_\lambda D_X^2}{\mu_y}
[
\log(
\tfrac{A_\lambda}{\rho_xa_0b_0}
)
]_+ .
$\\
Let
$A_k:=\|y-y_k\|_\lambda^2$.
Since $b_{k+1}=b_k+\mu_yt_k$,
\[
\tsum_{k=0}^{K-1}
[
\tfrac{b_k}{2}(A_k-A_{k+1})
-\tfrac{\mu_yt_k}{2}A_{k+1}
]
=
\tfrac{b_0}{2}A_0-\tfrac{b_K}{2}A_K.
\]
The extrapolation terms telescope to
$
\mathcal B_K
:=
t_{K-1}\langle q_x^K,x_K-x\rangle
-t_{K-1}\langle q_\lambda^K,\lambda_K-\lambda\rangle
-t_{K-1}\langle\Lambda q_y^K,y_K-y\rangle.
$
Bounding these terminal terms and combining the multiplier term with
\(-\tfrac{\beta a_0}{2}\|\lambda-\lambda_K\|_2^2\) yields the boundary
constant
\[
C_{\rm b}(\lambda)
:=
a_0(\tfrac1{c_1}+\tfrac1{c_2})D_X^2
+\tfrac{4a_0D_X^2}{\beta c_1^2}
+\|\lambda\|_1
 (L_g^{yx}D_XD_y+L_g^{yy}D_y^2).
\]
Consequently,
\begin{align}
\tsum_{k=0}^{K-1}t_kQ(z^{k+1},z)
\le&
\tfrac{a_0}{2}\|x-x_0\|^2
+
\tfrac{\beta a_0}{2}
\|\lambda-\lambda_0\|_2^2
+
\tfrac{b_0}{2}\|y-y_0\|_\lambda^2+
C_{yx}(\lambda)
+
C_{\rm b}(\lambda)
-
\tfrac{3\beta a_0}{8}
\|\lambda-\lambda_K\|_2^2 .
\label{eq:app-adap-master}
\end{align}
For
$z^*=(x^*,\lambda^*,y^*)$,
convexity, complementarity, feasibility, and the KKT condition yield
$
Q(z^{k+1},z^*)\ge0.
$\\
Let, $\qquad
\lambda^*
:=
\|\lambda^*\|_1+1,
\Delta_B
:=
[\lambda^*-B]_+.
$
Then
$\|\lambda^*\|_2\le B+\Delta_B$ and
$\|y^*-y_0\|_{\lambda^*}^2
\le
(B+\Delta_B)D_y^2.$
Define
\begin{align*}
V_*(B,\Delta_B)
:={}&
\tfrac{a_0}{2}D_X^2
+
\tfrac{\beta a_0}{2}(B+\Delta_B)^2
+
\tfrac{b_0}{2}(B+\Delta_B)D_y^2+
\tfrac{
2(L_g^{yx})^2(B+\Delta_B)D_X^2
}{\mu_y}
[
\log(
\tfrac{
2(L_g^{yx})^2(B+\Delta_B)
}{
\rho_xa_0b_0
}
)
]_+\\
&+
a_0
(
\tfrac1{c_1}+\tfrac1{c_2}
)D_X^2
+
\tfrac{4a_0D_X^2}{\beta c_1^2}+
(B+\Delta_B)
(
L_g^{yx}D_XD_y
+
L_g^{yy}D_y^2
).
\end{align*}
Applying \eqref{eq:app-adap-master} at $z=z^*$ gives
\[
\tfrac{3\beta a_0}{8}
\|\lambda_K-\lambda^*\|_2^2
\le
V_*(B,\Delta_B).
\]
Hence
$\|\lambda_K\|_2
\le
R_*(B,\Delta_B)
:=
B+\Delta_B
+
\sqrt{
\tfrac{
8V_*(B,\Delta_B)
}{
3\beta a_0
}
}.$
Thus the generated multiplier sequence is uniformly bounded.
Since
$
a_k
=
c_1\bar L
+
c_2L_g^{xx}\|\lambda_k\|_2,
$
we have
\[
a_k
\le
a_*(B,\Delta_B)
:=
c_1\bar L
+
c_2L_g^{xx}R_*(B,\Delta_B),
\]
and therefore
$
\Gamma_K
=
\tsum_{k=0}^{K-1}t_k
\ge
\tfrac{a_0}{a_*(B,\Delta_B)}K.
\label{eq:app-adap-Gamma}$
\paragraph{VI gap and feasibility.}
For any $x\in\widetilde X$, choose $\lambda=0$ in
\eqref{eq:app-adap-master}. Since
\[
Q(z^{k+1},(x,0,y))
\ge
\langle F(x),x_{k+1}-x\rangle,
\]
we obtain
$
\sup_{x\in\widetilde X}
\langle F(x),\bar x_K-x\rangle
\le
\tfrac{C_{\rm gap}}{\Gamma_K},
$
where
$
C_{\rm gap}
:=
\tfrac{a_0}{2}D_X^2+C_{\rm b}(0).
$
Using \eqref{eq:app-adap-Gamma} proves
\eqref{eq:adaptive-gap}.
For feasibility, choose
\[
j_K
\in
\arg\max_i[\phi_i(\bar x_K)]_+,
\qquad
\widehat y_{K,i}
\in
\arg\max_{y_i\in Y_i}
g_i(\bar x_K,y_i),
\]
and let $
\widehat\lambda_K
:=
\lambda^*+ 1.$
Then
$
\|\widehat\lambda_K\|_1
=
\lambda^*
\le
B+\Delta_B.
$
By convexity and the KKT conditions,
\[
\tsum_{k=0}^{K-1}
t_k
Q
(
z^{k+1},
(x^*,\widehat\lambda_K,\widehat y_K)
)
\ge
\Gamma_K
\|[\phi(\bar x_K)]_+\|_\infty.
\]
Applying \eqref{eq:app-adap-master} with this comparison point gives
\[
\|[\phi(\bar x_K)]_+\|_\infty
\le
\tfrac{C_{\rm feas}(B,\Delta_B)}{\Gamma_K},
\]
where
\begin{align*}
C_{\rm feas}(B,\Delta_B)
:={}&
\tfrac{a_0}{2}D_X^2
+
\tfrac{\beta a_0}{2}(B+\Delta_B)^2
+
\tfrac{b_0}{2}(B+\Delta_B)D_y^2+
\tfrac{
2(L_g^{yx})^2(B+\Delta_B)D_X^2
}{\mu_y}
[
\log(
\tfrac{
2(L_g^{yx})^2(B+\Delta_B)
}{
\rho_xa_0b_0
}
)
]_+\\
&+
a_0
(
\tfrac1{c_1}+\tfrac1{c_2}
)D_X^2
+
\tfrac{4a_0D_X^2}{\beta c_1^2}+
(B+\Delta_B)
(
L_g^{yx}D_XD_y
+
L_g^{yy}D_y^2
).
\end{align*}
Combining this with \eqref{eq:app-adap-Gamma} proves
\eqref{eq:adaptive-feas} and concludes the proof.
\end{sharedproof}

\bibliographystyle{abbrvnat}
\bibliography{ref.bib}

\end{document}